\documentclass[11pt,a4paper,reqno]{article}

\usepackage{amsmath,amssymb,exscale,mathrsfs}
\usepackage{graphics,graphicx,color}
\usepackage{float}

\usepackage{enumerate}

\usepackage{color}
\definecolor{marin}{rgb}   {0.,   0.3,   0.7} 
\definecolor{rouge}{rgb}   {0.8,   0.,   0.} 
\definecolor{sepia}{rgb}   {0.8,   0.5,   0.} 
\usepackage[colorlinks,citecolor=marin,linkcolor=rouge,
            bookmarksopen,
            bookmarksnumbered
           ]{hyperref}

\newtheorem{lemma}{Lemma}[section]

\newtheorem{proposition}[lemma]{Proposition}

\newtheorem{remark}[lemma]{Remark}
\newtheorem{example}[lemma]{Example}

\newtheorem{notation}[lemma]{Notation}
\newtheorem{definition}[lemma]{Definition}
\newtheorem{conclusion}[lemma]{Conclusion}
\numberwithin{equation}{section}
\newcommand{\QED}{\mbox{}\hfill \raisebox{-0.2pt}{\rule{5.6pt}{6pt}\rule{0pt}{0pt}} 
          \medskip\par}             

\newenvironment{Proofof}[1]{\noindent
    \parindent=0pt\abovedisplayskip = 0.5\abovedisplayskip
    \belowdisplayskip=\abovedisplayskip{\bfseries Proof of #1. }}{\QED}

\newcommand{\E}{\mathbb{E}}

\newcommand{\N}{\mathbb{N}}

\newcommand{\R}{\mathbb{R}}

\newcommand{\T}{\mathbb{T}}
\newcommand{\Z}{\mathbb{Z}}

\newcommand{\be}{\begin{equation}}
\newcommand{\ee}{\end{equation}}
\newcommand{\bea}{\begin{eqnarray}}
\newcommand{\eea}{\end{eqnarray}}
\newcommand{\bee}{\begin{eqnarray*}}
\newcommand{\eee}{\end{eqnarray*}}
\def\ds{\displaystyle}
\def\ni{\noindent}
\def\bs{\bigskip}
\def\ms{\medskip}

\def\eps{\varepsilon}

\def\fref#1{{\rm (\ref{#1})}}

\def\I{\mbox{\rm I}}
\def\pa{\partial}
\def\na{\nabla}

\def\im{{\mathcal I} \hspace{-2pt}{\textit m}\,}

\author{Nicolas Crouseilles \thanks{INRIA-Rennes Bretagne Atlantique, IPSO Project} 
  \and 
Mohammed Lemou
 \thanks{CNRS and IRMAR, Universit\'e de Rennes 1 and INRIA-Rennes Bretagne Atlantique, IPSO Project}
  \and 
Florian M\'ehats 
 \thanks{IRMAR, Universit\'e de Rennes 1 and INRIA-Rennes Bretagne Atlantique, IPSO Project}}
\title{Asymptotic preserving schemes for highly oscillatory kinetic equations}        
\begin{document}
\maketitle


\begin{abstract}
This work is devoted to the numerical simulation of a Vlasov-Poisson model describing a charged particle beam under the action 
of a rapidly oscillating external electric field. We construct an Asymptotic Preserving numerical scheme for this kinetic equation in the highly oscillatory limit. This scheme enables to simulate the problem without using any time step refinement technique. Moreover, since our numerical method is not based on the derivation of the simulation of asymptotic models, it works in the regime where the solution does not oscillate rapidly, and in the highly oscillatory regime as well. Our method is based on a "double-scale" reformulation of the initial equation, with the introduction of an additional periodic variable.
\end{abstract}


\section{Introduction} 


In this article, we are interested in the construction of numerical schemes for collisionless kinetic equations which involve rapid oscillations in time. Our study is done in the framework of a specific physical application, the case of a charged particle beam in the paraxial approximation, but our strategy can be applied to other highly oscillatory kinetic models, for instance in the physics of magnetized plasmas \cite{golsesaintraymond,frenodcg,frenodflr,bostan1,bostan2} for the guiding-center limit or the finite Larmor radius limit.

Let us first present our model. The paraxial approximation of the Vlasov-Maxwell equations concerns stationary, non collisional, charged particle beams which display a predominant length scale, called the longitudinal direction, such that the transverse width of the beam is very small compared to the typical longitudinal length. The paraxial model is obtained by expanding the Vlasov-Maxwell model with respect to the ratio $\eps>0$ between the characteristic lengths in the transverse and in the longitudinal directions, we refer to \cite{degond,filbet} for a derivation of this model. Here, following \cite{besse,frenod,mouton}, we consider the simpler case of an axisymmetric beam (with zero angular momentum). The paraxial Vlasov-Poisson model takes then the following form, in dimensionless variables,
\begin{equation}
\label{eqf1}
\partial_t f^\eps + \frac{v}{\varepsilon}\partial_r f^\eps +(E_{f^\eps}+E_{\rm app}) \partial_{v} f^\eps =0, 
\end{equation}
where $f^\eps(t, r, v)$ is the distribution function of the particles, $t\geq 0$ corresponds to the longitudinal position coordinate (the direction of propagation of the beam, denoted as a time here), $r\in \R_+$ is the radial component of the position in the transverse plane, and $v\in \R$ is the radial velocity in this plane. The total electric field has two contributions, the self-consistent electric field $E_{f^\eps}=E_{f^\eps}(t, r)$ satisfying the Poisson equation in the transverse plane, written in cylindrical symmetry as
\begin{equation}
\label{poisson1}
\frac{1}{r}\partial_r (r E_{f^\eps}) = \int_\mathbb{R} f^\eps dv
\end{equation}
and an applied electric field $E_{\rm app}$, chosen as in \cite{frenod} under the following form 
\begin{equation}
\label{eapp1}
E_{\rm app}(t, r) = -\frac{r}{\varepsilon} + a\left(\frac{t}{\varepsilon}\right) r,  
\end{equation}
where $a$ is a given $2\pi$-periodic function (the so-called tension function). This system is initially defined for $r\geq 0$ but can be extended to $r\in \mathbb{R}$ by using the conventions $f^\eps(t, -r, -v) = f^\eps(t, r, v)$ and $E(t, -r)=-E(t, r)$.

To summarize, in this paper we consider the following one-dimensional Vlasov-Poisson system satisfied by $f^\eps (t, r, v)$, where $r\in \R$ and $v\in \R$, 
 \begin{equation}
\label{eqf}
\partial_t f^\eps + \frac{v}{\varepsilon}\partial_r f^\eps +\left(E_{f^\eps} - \frac{r}{\varepsilon}+a\left(\frac{t}{\varepsilon}\right) r \right) \partial_{v} f^\eps =0,\qquad f^\eps(t=0,r,v)=f_0(r,v),
\end{equation}
\begin{equation}
\label{poisson}
E_{f^\eps}(t,r)=\frac{1}{r}\int_0^rs\rho^\eps(t,s)ds\quad \mbox{with}\quad \rho^\eps(t,r)=\int_\R f^\eps(t,r,v)dv. 
\end{equation}
The initial data $f_0$ is a given smooth function. When there is no confusion we shall omit the subscript $\eps$  to ease notations.

\bs
The main purpose of this work is the construction of robust numerical methods for stiff transport equations of type (\ref{eqf}) in the limit 
$\varepsilon\rightarrow 0$. We seek a method that is able to capture the properties of the various scales in the considered system, while the numerical parameters may be kept independent of the stiffness degree of these scales. Contrary to collisional kinetic equations in hydrodynamic or diffusion asymptotics, collisionless equations like (\ref{eqf}) involve time oscillations. In this context, the notion of two-scale convergence \cite{allaire, frenod, nguetseng, cfhm} is well-adapted in order to derive asymptotic models. However, these asymptotic models are valid only when $\varepsilon$ is small. In this paper, we develop numerical schemes that are able to deal with a wide range of values for $\eps$. We construct a numerical method in the so-called {\em Asymptotic Preserving} (AP) class \cite{jin}: such schemes are consistent with the kinetic model for all positive value of $\eps$, and degenerate into consistent schemes with the asymptotic model when $\eps \to 0$.

\bs
To do this, let us first put the stiff equation \fref{eqf} into a filtered form by rewriting it in the adapted rotating frame. The characteristic equations associated with \fref{eqf} read
$$\frac{d}{dt}\left(\begin{array}{ll}r \\v\end{array}\right)=\frac{1}{\varepsilon}J\left(\begin{array}{ll}r \\v\end{array}\right)+
\left(\begin{array}{ll}0 \\E_f(t,r)+a\left(t/\eps\right)r\end{array}\right),$$
where the matrix $J$ is defined by
$$J=
\left( 
\begin{array}{llcc}
0 & 1 \\
-1 & 0
\end{array}
\right).
$$
Hence, introducing the oscillatory variable $\xi\in \R^2$ defined by
\begin{equation}
\label{changevar}
\left(
\begin{array}{ll}
\xi_1 \\
\xi_2
\end{array}
\right)
=
e^{-Jt/\varepsilon}\left(\begin{array}{ll}r \\v\end{array}\right)
=\left(\begin{array}{ll}
\cos (t/\eps)&-\sin(t/\eps) \\ \sin(t/\eps)&\cos(t/\eps)
\end{array}\right)
\left(\begin{array}{ll}r \\v\end{array}\right),
\end{equation}
the associated filtered distribution function $\tilde{f^\eps}(t, \xi_1,\xi_2) = f^\eps(t, r, v)$ satisfies
\begin{equation}
\label{eqftilde}
\partial_t \tilde{f^\eps}(t, \xi) + \left(\widetilde E_{\tilde f^\eps}(t, t/\varepsilon,  \xi) +\widetilde E_{\rm app}(t/\varepsilon,  \xi) \right)\cdot \nabla_\xi \tilde f^\eps(t, \xi)=0,\qquad \tilde f^\eps(t=0,\cdot)=f_0,
\end{equation}
where the vector field is the sum of the applied field
\be
\label{Eapp}
\widetilde E_{\rm app}( \tau,  \xi) = a(\tau)(\xi_1 \cos\tau+ \xi_2 \sin\tau)
\left( \begin{array}{ll}
-\sin \tau \\
\cos \tau
\end{array}\right)
\ee
and of the self-consistent field defined by
$$
\widetilde E_{\tilde f}(t,\tau,\xi)=\left( \begin{array}{ll}
-\sin \tau \\
\cos \tau
\end{array}\right)\frac{1}{r(\tau,\xi)}
\int_0^{r(\tau,\xi)}\int_{-\infty}^{+\infty}s\tilde f\left(t,s \cos\tau- v \sin\tau,s\sin \tau+v\cos\tau\right)dsdv
$$
with $r(\tau,\xi)=\xi_1 \cos\tau+ \xi_2 \sin\tau$.

Let us briefly describe the strategy we propose to deal with equations like \fref{eqftilde}. As a matter of fact, we embed the function $\tilde f^\eps(t,\xi)$ into the family of solutions $F^\eps(t,\tau,\xi)$ of an "augmented" kinetic equation, where we separate the two scales $t/\eps$ and $t$. Assume indeed that $F^\eps$ solves the equation
\be
\label{eqFintro}
\pa_t F^\eps+\left(\widetilde E_{F^\eps}(t, \tau,  \xi) +\widetilde E_{\rm app}(\tau,  \xi) \right)\cdot \nabla_\xi F^\eps=-\frac{1}{\eps}\pa_\tau F^\eps,
\ee
and that, additionnally, we have 
\be
\label{initFintro}
\forall \xi\in \R^2,\qquad F^\eps(0,0,\xi)=f_0(\xi),
\ee
then it is readily seen that $F^\eps(t,t/\eps,\xi)$ satisfies the initial-value problem \fref{eqftilde}, so we recover $\tilde f^\eps (t,\xi)=F^\eps(t,t/\eps,\xi)$. The point is, in this double-scale formulation \fref{eqFintro} of \fref{eqftilde}, the stiffness is confined in the sole term $-\frac{1}{\eps}\pa_\tau F^\eps$ in the right-hand side. Reinterpreting this singularly perturbed term as a "collision" operator in this collisionless context, we can obtain the asymptotic behavior of $F^\eps$ (then of $\tilde f^\eps$) by a Chapman-Enskog expansion. In turn, this suggests a systematic method to construct Asymptotic Preserving numerical schemes, based on a micro-macro decomposition of $F$, see \cite{ml,bennoune,clemou}.

\bs

This paper is organized as follows. In Section \ref{gene}, we present the double-scale formulation in a general framework, and perform in subsection \ref{chapman} the Chapman-Enskog expansion of $F$. We then discuss in subsection \ref{discuinit} the crucial question of the choice of the initial data  $F^\eps(0,\tau,\xi)$ for this augmented kinetic equation \fref{eqFintro}. In this subsection, we state the main (formal) theoretical result  of this paper, in Proposition \ref{mainth}. In subsection \ref{sectionlinear}, we compute explicitely the averaged equations for our problem in the linear setting (when self-consistent interactions are neglected). Then, in Section \ref{micro-macro}, we present our AP numerical scheme. In subsection \ref{sectionscheme}, we introduce the scheme, which is a second order (in time and space) Eulerian numerical scheme. In subsection \ref{sectionAP}, we prove formally that this scheme is Asymptotic Preserving at the limit $\eps\to 0$. In subsection \ref{sectiondiffusion}, we show how the micro-macro decomposition method enables to construct AP schemes in more complicated situations, such as the diffusion limit.
Finally, the last Section \ref{numerics} is devoted to a series of numerical tests which characterize the properties of our scheme.  

\section{Double-scale formulation of the oscillatory equation}

\label{gene}

In this section, we introduce a general strategy in order to deal with highly oscillatory problems under the form
\begin{equation}
\label{eqftildegene}
\partial_t \tilde{f^\eps} + A(t, t/\varepsilon,  \xi,\tilde f^\eps)=0,\qquad \tilde f^\eps(t=0,\cdot)=f_0,
\end{equation}
where the unknown is the distribution function $(t,\xi)\in \R_+\times \R^d\mapsto \tilde f^\eps(t,\xi)\in \R$ and the vector-field $(t,\tau,\xi,f)\mapsto A(t,\tau,\xi,f)\in \R$ is a functional which is $P$-periodic with respect to the variable $\tau\in\T$ ($\T$ denotes the torus $\R/P\Z$).
Our target equation \fref{eqftilde} is under the form \fref{eqftildegene}, with $d=2$ and
$$A (t,\tau,\xi,f)=\left(\widetilde E_{f}(t, \tau,  \xi) +\widetilde E_{\rm app}(\tau,  \xi) \right) \cdot \nabla_\xi f.$$

\bs
We now introduce the following "double-scale formulation"
\be
\label{eqF}
\pa_t F^\eps+A(t, \tau,  \xi,F^\eps)=-\frac{1}{\eps}\pa_\tau F^\eps,
\ee
where the unknown is the function $(t,\tau,\xi)\in \R_+\times \T\times\R^d\mapsto F^\eps(t,\tau,\xi)$. This problem is an augmented version of \fref{eqftildegene}. Indeed if a function $F^\eps(t,\tau,\xi)$ solves \fref{eqF} and satisfies additionally
\be
\label{initF}
\forall \xi\in \R^d,\qquad F^\eps(0,0,\xi)=f_0(\xi),
\ee
then by differentiating $F^\eps(t,t/\eps,\xi)$ we obtain that $\tilde f^\eps(t,\xi):=F^\eps(t,t/\eps,\xi)$ satisfies the initial-value problem \fref{eqftildegene}.
 
It is important to note that \fref{eqF}, \fref{initF} is not sufficient to  uniquely determine the function $F^\eps$. Indeed, \fref{initF} is not a Cauchy condition for \fref{eqF}. The question of choosing a "good" initial condition $F(0,\tau,\xi)=F_0(\tau,\xi)$ for all $(\tau,\xi)\in \T\times \R^d$ is a delicate issue and is discussed in subsection \ref{discuinit}. In fact, we will see --\,in a formal setting\,-- that there is a unique way (up to order $\mathcal O(\eps^2)$ terms) to define $F_0$ in order to get a {\em smooth function}
$$(t,\tau,\xi,\eps)\in [0,t_{final}]\times \T\times \R^d\times [0,\eps_0[ \,\longmapsto \,F^\eps(t,\tau,\xi)$$
that satisfies \fref{eqF} and \fref{initF}. Here $t_{final}>0$ is a fixed final observation time and $\eps_0>0$ is arbitrary. The important point here is the assumed regularity with respect to $\eps$ when this parameter goes to zero.

More precisely, our aim is to ensure that the function $F^\eps$ and its two first derivatives $\pa_t F^\eps$ and $\pa^2_{t} F^\eps$ are bounded. Roughly speaking, this regularity is a constraint that prevents a dependency of $F^\eps$ in the fast variable $t/\eps$  (up to order $\mathcal O(\eps^2)$ terms), and $F^\eps$ will "only depend" on $t$, $\tau$ and $\xi$. Under this condition, one can pretend that we have succeeded in separating (up to order $\mathcal O(\eps^2)$ terms) the two scales $t$ and $\tau=t/\eps$ that were initially in \fref{eqftildegene}. The main result of this section is Proposition \ref{mainth}.

\subsection{Chapman-Enskog expansion}
\label{chapman}
In this subsection, we analyze {\em formally} the behavior of \fref{eqF} when $\eps\to 0$, assuming that its solution $F^\eps$ is smooth enough. To this aim, we carry out the Chapman-Enskog expansion of this function. Consider the following linear operator, defined for all periodic (regular) function $\tau\in \T\mapsto h(\tau)$ by
$$Lh=\pa_\tau h.$$
This operator is skew-adjoint with respect to the $L^2(\T)$ scalar product and \fref{eqF} can be rewritten
\be
\label{eqFL}
\pa_t F^\eps+A(t, \tau,  \xi,F^\eps)=-\frac{1}{\eps}LF^\eps.
\ee
The kernel of $L$ is the set of constant functions and the $L^2$ projector on this kernel is the average
$$\Pi h:=\frac{1}{|\T|}\int_{\T} h(\tau) d\tau,$$ 
where $|\T| = P$ is the measure of $\T$.

Moreover, $L$ is invertible in the set of functions with zero average and, if $\int_\T h(\tau)d\tau=0$, we have
$$
(L^{-1}h)(\tau)=(\I-\Pi)\int_0^\tau h(\sigma)d\sigma=\int_0^\tau h(\sigma)d\sigma+\frac{1}{|\T|}\int_\T\sigma h(\sigma)d\sigma.
$$
Performing the Chapman-Enskog expansion of $F^\eps(t,\tau,\xi)$ consists in writing
\be
\label{decomp}
F^\eps(t,\tau,\xi)=G^\eps(t,\xi)+h^\eps(t,\tau,\xi)\qquad \mbox{with}\quad G^\eps(t,\xi)=\Pi \left(F^\eps(t,\tau,\xi)\right)
\ee
and deriving asymptotic equations for $G^\eps$ and $h^\eps$ when $\eps\to 0$. As we said, we proceed at a formal level, and the rule that we follow in this analysis is that $F^\eps$ is assumed to be smooth with respect to all its variables (in particular with respect to the parameter $\eps$ which can be very small).

Inserting the decomposition \fref{decomp} into \fref{eqFL} leads to 
\be
\label{eqinsert}
\partial_t G^\eps+\partial_t h^\eps + A(t, \tau,\xi,G^\eps+h^\eps)= -\frac{1}{\varepsilon}L h^\eps. 
\ee
Averaging this last equation with respect to $\tau$ (i.e. applying $\Pi$) yields, since $\Pi h^\eps=0$,
\be
\label{eqmoy}
\partial_t G^\eps+ \Pi \left(A(t, \tau,\xi,G^\eps+h^\eps)\right) = 0. 
\ee
Then, from \fref{eqinsert} and \fref{eqmoy} we deduce that $h^\eps$ satisfies
\be
\label{eqh}
\partial_t h^\eps +(\I-\Pi)\left(A(t, \tau,\xi,G^\eps+h^\eps)\right)= -\frac{1}{\varepsilon}L h^\eps.  
\ee
Now, from \fref{eqh} and the fact that $h^\eps$ belongs to the range of $L$, we deduce that
\be
\label{h-esti1}
h^\eps=-\eps L^{-1}\left(\partial_t h^\eps +(\I-\Pi)\left(A(t, \tau,\xi,G^\eps+h^\eps)\right)\right).
\ee
Hence, using our smoothness assumption and in particular that we have $\pa_t F^\eps=\mathcal O(1)$, $\pa^2_t F^\eps=\mathcal O(1)$ (hence $G^\eps$ and $h^\eps$ have also bounded derivatives), we deduce from \fref{h-esti1} that
$$h^\eps=\mathcal O(\eps)\quad\mbox{and}\quad  \pa_t h=\mathcal O(\eps).$$
From these estimates and \fref{eqmoy}, we deduce a first approximate equation satisfied by $G^\eps$:
\be
\label{eq0G}
\partial_t G^\eps+ \Pi A(t,\tau,\xi, G^\eps) = \mathcal O(\eps).
\ee
Next, using again \fref{h-esti1}, we obtain an expression of $h^\eps$ in terms of $G^\eps$, up to a small remainder:
\be
\label{h-fonction-G}
h^\eps=-\eps L^{-1}(\I-\Pi)A(t,\tau,\xi,G^\eps)+\mathcal O(\eps^2)
\ee
and this expression, together with \fref{eqmoy}, enables to derive the following equation satisfied by $G^\eps$ up to second order terms: 
\be
\label{eq1G}
\partial_t G^\eps+\Pi A(t,\tau,\xi, G^\eps)-\eps \Pi \left(\pa_f A(t,\tau,\xi,G^\eps)\left(L^{-1}(\I-\Pi)A(t,\tau,\xi, G^\eps)\vphantom{\widetilde G}\right)\right)=\mathcal O(\eps^2).
\ee
Finally, the function $F^\eps$ can be deduced from $G^\eps$, up to second order terms, by using \fref{decomp} and \fref{h-fonction-G}:
\be
\label{Fordre2v0}
F^\eps=G^\eps-\eps L^{-1}(\I-\Pi)A(t,\tau,\xi,G^\eps)+\mathcal O(\eps^2).
\ee

\subsection{Discussion on the initial data and main result}
\label{discuinit}

In the previous subsection, the Chapman-Enskog expansion was performed formally under a regularity assumption on $F^\eps$. In this subsection, we reverse the argument and deduce from these expansions a Cauchy data for \fref{eqF} that ensures that $F^\eps$ is regular enough (up to order $\mathcal O(\eps^2)$ terms).

A natural initial condition for \fref{eqF} can be deduced from \fref{Fordre2v0}. Indeed, by evaluating \fref{Fordre2v0} at $t=0$, one gets 
\be
\label{condinitv0}
F^\eps(0, \tau,\xi) = G^\eps(0,\xi) - \eps (I-\Pi)\int_0^\tau (I-\Pi)A(0, s, \xi,G^\eps(0,\xi))ds+\mathcal O(\eps^2)
\ee
and then, by taking this equation at $\tau=0$ and by using \fref{initF},
\be
\label{condinitv1}
f_0(\xi) = G^\eps(0,\xi) + \eps\Pi \int_0^\tau (I-\Pi)A(0, s, \xi,G^\eps(0,\xi))ds+\mathcal O(\eps^2).
\ee
By substracting these two identities \fref{condinitv0} and \fref{condinitv1}, one gets
\be
\label{condinitv2}
F^\eps(0, \tau,\xi) = f_0(\xi)- \eps \int_0^\tau (I-\Pi)A(0, s, \xi,G^\eps(0,\xi))ds+\mathcal O(\eps^2).
\ee
Moreover, from \fref{condinitv1}, one deduces $G^\eps(0,\xi)=f_0(\xi)+\mathcal O(\eps)$, which can finally be inserted into \fref{condinitv2} and yields
\be
\label{condinitv3}
F^\eps(0, \tau,\xi) = f_0(\xi)- \eps \int_0^\tau (I-\Pi)A(0, s, \xi,f_0(\xi))ds+\mathcal O(\eps^2).
\ee
The correction term in $\eps$ is important here and, as we show further, will guarantee that $F^\eps(t,\tau,\xi)$ does not oscillate in time. By analogy with boundary value problems in collisional kinetic theory (see \cite{golsejin}), one can interpret this term as "boundary corrector" (where the boundary is the initial time $t=0$). The interesting point in our case is that we {\em do not have} to assume that the initial data is well-prepared since, as we said in the introduction of this section, we have a degree of freedom on $F_0$ which is not totally prescribed. We have then the possibility to enforce that \fref{condinitv3} is satisfied (see \fref{defF0}).

\bs
Let us formulate in the following proposition the main result of this section.
\begin{proposition}[formal]
\label{mainth}
Let $F^\eps(t,\tau,\xi)$ be the unique solution of \fref{eqF} subject to the initial condition
\be
\label{condinit}
\forall (\tau,\xi)\in\T\times \R^2,\qquad  F^\eps(0, \tau,\xi) = F_0^\eps(\tau,\xi)
\ee
with $F_0^\eps$ defined for all $\tau\in\T$ and $\xi\in\R^2$ by
\be
\label{defF0}
F_0^\eps(\tau,\xi)=f_0(\xi)- \eps \int_0^\tau (I-\Pi)A(0, s, \xi,f_0(\xi))ds.
\ee
Then we have
\be
\label{Fordre2}
F^\eps(t,\tau,\xi)=\widetilde G^\eps(t,\xi)-\eps L^{-1}(\I-\Pi)A(t,\tau,\xi,\widetilde G^\eps(t,\xi))+\mathcal O(\eps^2),
\ee
where $\widetilde G^\eps(t,\xi)$ is the solution of the initial-value problem
\bea
\label{eqGtilde}
\hspace*{-1.8cm}&&\partial_t \widetilde G^\eps+\Pi A(t,\tau,\xi,\widetilde G^\eps)-\eps \Pi \left(\pa_f A(t,\tau,\xi,\widetilde G^\eps)\left(L^{-1}(\I-\Pi)A(t,\tau,\xi, \widetilde G^\eps)\right)\right)=0,\\
\label{condinitGtilde}
\hspace*{-1.8cm}&&\widetilde G^\eps(0,\xi)=\Pi F_0^\eps(\tau,\xi) = f_0(\xi)- \eps \Pi \int_0^\tau (I-\Pi)A(0, s, \xi,f_0(\xi))ds.
\eea
\end{proposition}
\begin{remark}
Since \fref{condinit} and \fref{defF0} imply \fref{initF}, one can recover the solution $\widetilde f^\eps$ to the oscillatory equation \fref{eqftildegene} by setting
$$\widetilde f(t,\xi)=F^\eps(t,t/\eps,\xi).$$
Moreover, whereas the term of order  $\eps$ in $\widetilde f$ varies rapidly in time (it depends on $t$ and $t/\eps$), the corresponding term in the function $F^\eps$ is smooth --\,since $\tau$ replaces the variable $t/\eps$\,-- and then easier to compute numerically. Our Asymptotic Preserving numerical method is constructed on the double-scale formulation \fref{eqF} instead of \fref{eqftildegene}.
\end{remark}
\begin{Proofof}{Proposition \ref{mainth}}
Let $F^\eps$ be the solution of \fref{eqF}, \fref{condinit} and let $\widetilde G^\eps$ be the solution of \fref{eqGtilde}, \fref{condinitGtilde}. Denote
$$\widetilde F^\eps=\widetilde G^\eps+\widetilde h^\eps+\eps^2 \chi^\eps$$
with
$$\widetilde h^\eps=-\eps L^{-1}(\I-\Pi)A(t,\tau,\xi,\widetilde G^\eps)$$
and where $\chi^\eps$ is a bounded corrector that is defined below (see \fref{corrector}). Proving the Proposition amounts to proving that $$F^\eps(t,\tau,\xi)-\widetilde F^\eps(t,\tau,\xi)=\mathcal O(\eps^2).$$
By substracting \fref{condinit} and \fref{condinitGtilde}, one gets
\bee
F^\eps(0,\tau,\xi)-\widetilde G^\eps(0,\xi)&=&-\eps (\I-\Pi)\int_0^\tau (I-\Pi)A(0, s, \xi,f_0(\xi))ds\\
&=&-\eps (\I-\Pi)\int_0^\tau (I-\Pi)A(0, s, \xi,\widetilde G^\eps(0,\xi))ds+\mathcal O(\eps^2)\\
&=&-\eps L^{-1}(I-\Pi)A(0, s, \xi,\widetilde G^\eps(0,\xi))ds+\mathcal O(\eps^2)\\
&=&\widetilde h^\eps(0,\tau,\xi)+\mathcal O(\eps^2).
\eee
This gives
\be
\label{initw}
F^\eps(0,\tau,\xi)-\widetilde F^\eps(0,\tau,\xi)=-\eps^2\chi^\eps(0,\tau,\xi)+\mathcal O(\eps^2)=\mathcal O(\eps^2).
\ee
Let us now derive an approximate equation satisfied by $\widetilde F^\eps(t,\tau,\xi)$. By inserting $\widetilde F^\eps$ in the equation \fref{eqFL}, one gets
\bee
&&\hspace*{-1.5cm}\pa_t \widetilde F^\eps+\frac{1}{\eps}L\widetilde F^\eps+A(t,\tau,\xi,\widetilde F^\eps)\\
&=&\pa_t \widetilde G^\eps-\eps L^{-1}(\I-\Pi)\pa_fA(t,\tau,\xi,\widetilde G^\eps)\pa_t \widetilde G^\eps-\eps L^{-1}(\I-\Pi)\pa_tA(t,\tau,\xi,\widetilde G^\eps)\\
&&-(\I-\Pi)A(t,\tau,\xi,\widetilde G^\eps)+\eps L\chi^\eps\\
&&+A(t,\tau,\xi,\widetilde G^\eps)-\eps \pa_f A(t,\tau,\xi,\widetilde G^\eps)\left(L^{-1}(\I-\Pi)A(t,\tau,\xi,\widetilde G^\eps)\right)+\mathcal O(\eps^2)\\
&=&-\Pi A(t,\tau,\xi,\widetilde G^\eps)+\eps \Pi \left(\pa_f A(t,\tau,\xi,\widetilde G^\eps)\left(L^{-1}(\I-\Pi)A(t,\tau,\xi, \widetilde G^\eps)\right)\right)\\
&&-\eps L^{-1}(\I-\Pi)\pa_fA(t,\tau,\xi,\widetilde G^\eps)\Pi A(t,\tau,\xi,\widetilde G^\eps)
-\eps L^{-1}(\I-\Pi)\pa_tA(t,\tau,\xi,\widetilde G^\eps)\\
&&-(\I-\Pi)A(t,\tau,\xi,\widetilde G^\eps)+\eps L\chi^\eps\\
&&+A(t,\tau,\xi,\widetilde G^\eps)-\eps \pa_f A(t,\tau,\xi,\widetilde G^\eps)\left(L^{-1}(\I-\Pi)A(t,\tau,\xi,\widetilde G^\eps\right)+\mathcal O(\eps^2)\\
&=&\eps L\chi^\eps-\eps (\I-\Pi)\left(\pa_f A(t,\tau,\xi,\widetilde G^\eps)\left(L^{-1}(\I-\Pi)A(t,\tau,\xi,\widetilde G^\eps\right)\right)\\
&&-\eps L^{-1}(\I-\Pi)\pa_tA(t,\tau,\xi,\widetilde G^\eps)\\
&&-\eps L^{-1}(\I-\Pi)\pa_fA(t,\tau,\xi,\widetilde G^\eps)\Pi A(t,\tau,\xi,\widetilde G^\eps)+\mathcal O(\eps^2)
\eee
where we used \fref{eqGtilde} in the second equality. Hence, by defining the corrector $\chi^\eps$ as
\bea
\chi^\eps&=&L^{-1}\left[(\I-\Pi)\left(\pa_f A(t,\tau,\xi,\widetilde G^\eps)\left(L^{-1}(\I-\Pi)A(t,\tau,\xi,\widetilde G^\eps\right)\right)\right.\nonumber\\
&&\hspace*{-10mm}\left.+L^{-1}(\I-\Pi)\pa_tA(t,\tau,\xi,\widetilde G^\eps)+L^{-1}(\I-\Pi)\pa_fA(t,\tau,\xi,\widetilde G^\eps)\Pi A(t,\tau,\xi,\widetilde G^\eps)\right]\label{corrector}
\eea
one gets finally
\be
\label{eqappFtilde}
\pa_t \widetilde F^\eps+\frac{1}{\eps}\pa_\tau \widetilde F^\eps+A(t,\tau,\xi,\widetilde F^\eps)=\mathcal O(\eps^2).
\ee
We shall now conclude by integrating the characteristics associated to this equation. Let
$$w(t,\tau,\xi)=(F^\eps-\widetilde F^\eps)(t,\tau+t/\eps,\xi).$$
From \fref{eqF} and \fref{eqappFtilde}, one deduces
$$\pa_t w=-A(t,\tau+t/\eps,\xi,F^\eps(t,\tau+t/\eps,\xi))+A(t,\tau+t/\eps,\xi,\widetilde F^\eps(t,\tau+t/\eps,\xi))+\mathcal O(\eps^2).$$
Hence, using the estimate \fref{initw} at the initial time, a Gronwall lemma yields formally (recall that $A$ is periodic with respect to $\tau$)
$$w(t,\tau,\xi)=\mathcal O(\eps^2)$$
for $t\in[0,T]$, $\tau\in \T$, $\xi\in \R^2$,
and, finally, one has proved that
$$\left(F^\eps-\widetilde F^\eps\right)(t,\tau,\xi)=w (t,\tau-t/\eps,\xi)=\mathcal O(\eps^2).$$
The proof of Proposition \ref{mainth} is complete.
\end{Proofof}

\begin{remark}
In fact, this averaging procedure can be pushed forward to higher orders in $\eps$ by iterating further the Chapman-Enskog procedure. 
Other approaches may be used to obtain formally higher order averaged models for $F^\eps$, under a higher order initial condition, see for instance the approach developed in \cite{perko,SV} which is widely used in the context of ODEs. However, the purpose of this paper being to build an AP numerical method for our problem, we stop this construction at order $\mathcal O(\eps^2)$. We also refer to \cite{strobo} for a presentation of the so-called stroboscopic averaging in a way which is very close to the method introduced here. Indeed, in \cite{strobo}, a systematic construction of high order averaged models for oscillatory equations such as \fref{eqftilde} is based on the transport equation \fref{eqF}. It is proved in this paper that, for any fixed integer $N>0$, the solution of \fref{eqftilde} can be written under the form (omitting the dependencies in $\xi$ for simplicity and assuming that $A(t,\tau,f)$ does not depend on $t$)
\be
\label{strobo}
\tilde f^\eps(t)=\Phi^{\eps,N}\left(t/\eps,G^{\eps,N}(t)\right)+\mathcal O(\eps^{N+1}),
\ee
where $G^{\eps,N}(t)$ satisfies an autonomous averaged equation of the form $\pa_t G=A^{\eps,N}_{\rm av}\left(G\right)$ with  $G^{\eps,N}(0)=f_0$ and where $(\tau,f)\mapsto \Phi^{\eps,N}(\tau,f)$ is a close-to-identity mapping which is $2\pi$-periodic with respect to $\tau$ and satisfies $\Phi^{\eps,N}(0,f)=f$. The link with our construction is the following. If we choose
$F_0(\tau)=\Phi^{\eps,N}\left(\tau,f_0\right)$ as initial data for \fref{eqF}, then the stroboscopic averaging result says that $F(\tau,t)=\Phi^{\eps,N}\left(\tau,G^{\eps,N}(t)\right)+\mathcal O(\eps^{N+1})$, i.e. $F(\tau,t)$ is smooth, up to $\mathcal O(\eps^{N+1})$ terms. This gives the natural generalization of our initial data \fref{defF0} in order to get higher order estimates. Of course, one can check that, for $N=1$,
$$\Phi^{\eps,1}\left(\tau,f_0\right)=f_0- \eps \int_0^\tau (I-\Pi)A(s,f_0(\xi))ds.$$

\end{remark}

\subsection{The case of a linear transport equation}
\label{sectionlinear}

In this subsection, we compute explicitely the initial condition $F_0$ and the averaged system in the special situation of the following linear transport equation in dimension $d=2$:
\begin{equation}
\label{eqftildelinear}
\partial_t \tilde{f^\eps} + E(\tau,  \xi) \cdot \nabla_\xi \tilde f^\eps=0,\qquad \tilde f^\eps(t=0,\cdot)=f_0,
\end{equation}
where the field $E(\tau,\xi)=\left(\begin{array}{c}E_1(\tau,\xi)\\ E_2(\tau,\xi)\end{array}\right)$ is given and divergence-free. This equation is under the form \fref{eqftildegene} with
\be
\label{Alinear}
A (\tau,\xi,f)=E(\tau,  \xi) \cdot \nabla_\xi f.
\ee
In particular, when the self-consistent Poisson field $E_f$ is neglected, the filtered equation \fref{eqftilde} associated to the paraxial beam model \fref{eqf} is under this form, with $E(\tau,\xi)=\widetilde E_{\rm app}(\tau,\xi)$ defined by \fref{Eapp} (it is a divergence-free vector field).

In this linear case, the following proposition is a variant of Proposition \ref{mainth}.
\begin{proposition}[formal]
\label{theo2}
Assume that $A$ takes the form \fref{Alinear}. Let $F^\eps(t,\tau,\xi)$ be the unique solution of \fref{eqF} subject to the initial condition $F^\eps(0, \tau,\xi) = F_0(\tau,\xi)$
with $F_0$ defined for all $\tau\in\T$ and $\xi\in\R^2$ by
\be
\label{defF0linear}
F_0(\tau,\xi)=f_0\left(\xi- \eps \int_0^\tau (I-\Pi)E(s,\xi)ds\right).
\ee
Then we have
\be
\label{Fordre2linear}
F^\eps(t,\tau,\xi)=\widetilde G^\eps\left(t,\xi-\eps L^{-1}(\I-\Pi)E(\tau,\xi) \right)+\mathcal O(\eps^2),
\ee
where $\widetilde G^\eps(t,\xi)$ is the solution of the averaged transport equation
\bea
\label{eqGtildelinear}
\hspace*{-1.8cm}&&\partial_t \widetilde G^\eps+\left(E^{(0)}+\eps E^{(1)}\right)\cdot\na_\xi \widetilde G^\eps=0,\\
\label{condinitGtildelinear}
\hspace*{-1.8cm}&&\widetilde G^\eps(0,\xi)=f_0\left(\xi- \eps \Pi \int_0^\tau (I-\Pi)E(s,\xi)ds\right),
\eea
and where $E^{(0)}=\Pi E$ and $E^{(1)}=J^{-1}\na_\xi\mathcal D$ is the vector-field associated with the Hamiltonian
$$
\mathcal D(\xi)= \frac{1}{|\T|} \int_{\T} \left[  (I-\Pi)E_2\right](\tau,\xi)\int_0^\tau (I-\Pi)E_1(s,\xi)ds d\tau. 
$$
\end{proposition}
\begin{remark}
This result is the Eulerian version of an averaging theorem formulated directly in terms of the characteristics equations associated to the vector field $E(\tau,\xi)$. Indeed, consider the flow $\Xi$ associated to the averaged vector field: $\Xi(t,t_0,\xi_0)$ solves
$$\frac{d\Xi}{dt}=E^{(0)}(\Xi)+\eps E^{(1)}(\Xi),\qquad \Xi(t_0,t_0,\xi_0)=\xi_0.$$
Then we have $F^\eps(t,\tau,\xi)=f_0(\widetilde\Xi(t,\tau,\xi))+\mathcal O(\eps^2)$, where $\widetilde \Xi$ is defined by
$$\widetilde\Xi(t,\tau,\xi)=\left(\I- \eps \Pi \int_0^\tau (I-\Pi)E(s,\cdot)ds\right)\left(\Xi\left(0,t,\xi-\eps \int_0^\tau (I-\Pi)E(s,\xi)ds\right)\right).$$
\end{remark}
\begin{remark}The averaged equation (up to the the order $\mathcal O(\eps^2)$) shares the geometric structure of the initial equation \fref{eqftildelinear}. Indeed, since $E$ is divergence-free, so is $E^{(0)}$ and if $E$ is Hamiltonian, with Hamiltonian $H(\tau,\xi)$, then $E^{(0)}$ is Hamiltonian, with Hamiltonian given by $H^{(0)}=\Pi H$. Moreover, the correction $\eps E^{(1)}$ is always divergence-free and Hamiltonian.
\end{remark}

\begin{Proofof}{Proposition \ref{theo2}}
The initial data \fref{defF0linear} and \fref{condinitGtildelinear} can be deduced from \fref{defF0} and \fref{condinitGtilde} by a Taylor expansion, up to order $\mathcal O(\eps^2)$ terms: one has indeed
$$f_0\left(\xi+ \eps B(\tau,\xi)\right)=f_0(\xi)+ \eps B(\tau,\xi)\cdot\na_\xi f_0(\xi)+\mathcal O(\eps^2).
$$
Similarly, the change of variable \fref{Fordre2linear} can be deduced from \fref{Fordre2} by a Taylor expansion. Moreover, we have clearly $\Pi A(t,\tau,\xi,\widetilde G^\eps)=E^{(0)}\cdot \na_\xi \widetilde G^\eps$. Hence, to end the proof of the proposition, we simply have to compute the first order correction in the equation of $\widetilde G^\eps$ given by Proposition \ref{mainth}, i.e. the operator
\bee
G&\mapsto &-\Pi \left(\pa_f A(t,\tau,\xi,G)\left(L^{-1}(\I-\Pi)A(t,\tau,\xi, G)\right)\right)\\
&&\qquad =-\frac{1 }{|\T|}\int_{\T} E\cdot\na_\xi \left(L^{-1}(\I-\Pi)E\cdot\na_\xi G\right)d\tau= \na_\xi\cdot \left(\mathbb D\na_\xi G\right),
\eee
where we used that $E$ is divergence-free and where $\mathbb D$ is the $2\times 2$ "diffusion" matrix of components  
$$
\mathbb D_{i,j}= - \frac{1 }{|\T|}\int_{\T} E_i  L^{-1} [(I-\Pi) E_j]  d\tau, \;\;\; i, j=1, 2.  
$$
In fact, this matrix $\mathbb D$ inherits the skew-symmetry property of $L$. Indeed, for all $i,j$, we have
 \begin{eqnarray*}
\mathbb D_{i,j}&=& -\frac{1 }{|\T|}\int_{\T} E_i  (I-\Pi) L^{-1} [(I-\Pi) E_j]  d\tau\\
&=& -\frac{1 }{|\T|}\int_{\T} \left[(I-\Pi) E_i\right]  L^{-1} [(I-\Pi) E_j]  d\tau\\
&=&\frac{1 }{|\T|}\int_{\T} L^{-1} [(I-\Pi) E_i] \;  (I-\Pi) E_j  d\tau= -\mathbb D_{j,i}. 
\end{eqnarray*}
Hence, setting $\mathcal D=\mathbb D_{1,2}=-\mathbb D_{2,1}$, the "diffusion" term $\nabla_\xi \cdot(\mathbb D \nabla_\xi G)$ can be simplified as
$$
\nabla_\xi\cdot (\mathbb D \nabla_\xi G) =\partial_{\xi_1} (\mathcal D\partial_{\xi_2} G)-\partial_{\xi_2} (\mathcal D\partial_{\xi_1} G)=(\partial_{\xi_1} \mathcal D)\partial_{\xi_2} G-(\partial_{\xi_2} \mathcal D)\partial_{\xi_1} G,
$$
which is the desired result. Note that the first order model is a pure transport equation and does not include second order derivative.
\end{Proofof}

\bs
\ni
{\bf Explicit calculations in our example.}

\nopagebreak
\ms
\ni
Let us compute explicitely the approximate model in terms of the Fourier coefficients of $E$. 
From now, the period is taken as $|\T|=P=2\pi$. Introduce the decomposition of the two (real-valued) components of the vector field $E$ on the Fourier basis:
$$E_j(\tau,\xi)=\sum_{k\in \Z}A_{k,j}(\xi)e^{ik\tau}\quad \mbox{for }j=1,2,$$
with $A_{-k,j}=\overline{A_{k,j}}$ for all $k\in \Z$ and $j=1,2$.
Then direct calculations yield
\be
\label{expliciteFourier}
E^{(0)}_j=A_{0,j}\quad \mbox{for}\,j=1,2\quad \mbox{and}\quad \mathcal D=2\im\sum_{k\in\N^*}\frac{1}{k}A_{k,1}\overline{A_{k,2}}\,.
\ee
We now calculate the quantities defined in Proposition \ref{theo2} in a specific example that we use later for numerical experiments. In the beam model \fref{eqftilde}, if we neglect the Poisson field, then we have $E=\widetilde E_{\rm app}$ defined by \fref{Eapp}. 
Choosing $a(\tau)=\cos^2(2\tau)$, one computes from \fref{Eapp} the Fourier coefficients of $E_1$ and $E_2$:
\bee
E_1&=&\frac{1}{16}\left(-4\xi_2+(3\xi_2+i\xi_1)e^{2i\tau}-2\xi_2e^{4i\tau}+(\xi_2+i\xi_1)e^{6i\tau}+c.c.\right)\\
E_2&=&\frac{1}{16}\left(4\xi_1+(3\xi_1-i\xi_2)e^{2i\tau}+2\xi_1e^{4i\tau}+(\xi_1-i\xi_2)e^{6i\tau}+c.c.\right).
\eee
Hence, we obtain by simple integrations
\be
\label{d0d1}
\Pi\int_0^\tau(\I-\Pi)E=D_0\xi\quad \mbox{and}\quad L^{-1}(\I-\Pi)E=D_1(\tau)\xi,
\ee
with
\bee
D_0&=&\frac{1}{12}\left(\begin{array}{cc}-1&0\\0&1 \end{array}\right),\\
D_1&=&\frac{1}{48}\left(\begin{array}{cc}3\cos(2\tau)+\cos(6\tau)&9\sin(2\tau)-3\sin(4\tau)+\sin(6\tau)\\9\sin(2\tau)+3\sin(4\tau)+\sin(6\tau)&-3\cos(2\tau)-\cos(6\tau) \end{array}\right),
\eee
and also, from \fref{expliciteFourier}, we obtain that the averaged vector field (up to order $\mathcal O(\eps^2)$ terms) is the following Hamiltonian vector field:
$$
E^{(0)}+\eps E^{(1)}=J^{-1}\na_\xi \mathcal H, 
$$
with 
\be
\label{omega}
\mathcal H=\frac{\omega}{2}(\xi_1^2+\xi_2^2),\quad\omega=\omega_0+\eps\omega_1 = \left(\frac{1}{4}+\frac{5\eps}{192}\right).
\ee
The averaged equation \fref{eqGtildelinear}, \fref{condinitGtildelinear} for $\widetilde G^\eps$ is thus the equation of a rotation in the phase space and has an explicit solution:
\be
\label{Gexplicit}\widetilde G^\eps(t,\xi)=\widetilde G^\eps(0,e^{t\omega J}\xi)=f_0\left((\I-\eps D_0)e^{t\omega J}\xi\right).
\ee
We have thus analytic expressions for the solution of the limit model as $\eps\to 0$ and also for the solution of a next order approximation, which are then easy to implement numerically. The solution of the limit model reads (see \cite{frenod})
\be
\label{lim1}
F_{\rm limit}(t,\tau,\xi)=f_0\left(e^{t\omega_0 J}\xi\right)
\ee
and the solution of the second order model will be
\be
\label{lim2}
F_{\rm second order}(t,\tau,\xi)=f_0\left((\I-\eps D_0)e^{t\omega J}(\I-\eps D_1(\tau))\xi\right).  
\ee
This last relation is obtained by using successively \fref{Fordre2linear}, \fref{Gexplicit} and recalling that $D_0$ and $D_1$ are given by (\ref{d0d1}).  Indeed, one has
\bee
F^\eps(t, \tau, \xi) &=& \widetilde G^\eps(t,(I-\varepsilon D_1)\xi)+{\cal O}(\varepsilon^2)\\
&=& \widetilde G^\eps(0,e^{t\omega J}(I-\varepsilon D_1)\xi)+{\cal O}(\varepsilon^2)\\ 
&=&f_0\left((\I-\eps D_0)e^{t\omega J}(\I-\eps D_1(\tau))\xi\right)+{\cal O}(\varepsilon^2).
\eee

\section{Asymptotic Preserving numerical schemes}
\label{micro-macro}

In this section, we construct some Asymptotic Preserving numerical schemes for \fref{eqF}, hence for the original problem \fref{eqf}. Let us insist on the fact that we do {\em not} base the construction of our numerical method on the approximate models derived in the previous section, since we want a method which is efficient for the regimes where $\eps$ small {\em and} where $\eps=\mathcal O(1)$.

Recall that, in order to solve the filtered equation \fref{eqftilde}, we have introduced the augmented equation
\be
\label{eqFnum}
\pa_t F^\eps+E(t,\tau,\xi)\cdot \nabla_\xi F^\eps=-\frac{1}{\eps}\pa_\tau F^\eps,
\ee
where we denote for simplicity the field (which depends on the unknown $F^\eps$) by
$$E(t,\tau,\xi)=\widetilde E_{F^\eps}(t, \tau,  \xi) +\widetilde E_{\rm app}(\tau,  \xi).$$
After the asymptotic analysis in the previous section, and according to Proposition \ref{mainth} (see also Proposition \ref{theo2}), we know (see (\ref{defF0linear})) that a suitable initial condition for this problem is $F(0,\tau,\xi)=F_0(\tau,\xi)$ with
\be
\label{defF0bis}
F_0(\tau,\xi)=f_0\left(\xi- \eps \int_0^\tau (I-\Pi)E(0,s,\xi)ds\right).
\ee
Note that this choice is asymptotically close to \fref{defF0}, up to order $\mathcal O(\eps^2)$ terms, but is preferable since it garantees the positivity of the initial distribution function. Under this choice, we know two important facts:
\begin{itemize}
\item[--] one recovers the solution of \fref{eqftilde} by $\tilde f(t,\xi)=F^\eps(t,t/\eps,\xi)$,
\item[--] the function $F^\eps$ is smooth and, up to terms of order $\mathcal O(\eps^2)$, does not oscillate. In particular, its derivatives $\pa_t F^\eps$ and $\pa^2_t F^\eps$ are bounded when $\eps\to 0$.
\end{itemize}
In order to emphasize the role of the choice of the initial condition $F_0$, in our numerical experiments we will also test the most simple choice:
\be
\label{defF0ter}
F_0(\tau,\xi)=f_0(\xi).
\ee
This choice only garantees that $F^\eps=G^\eps+\mathcal O(\eps)$: we show below that, with this initial data, the numerical method will capture the right limit, but not the details of order $\mathcal O(\eps)$. In the sequel, the initial condition \fref{defF0bis} is referred to as "with correction", and the initial condition \fref{defF0ter} is referred to as "without correction".

\subsection{The numerical scheme}
\label{sectionscheme}
In this subsection, we present our AP numerical scheme. Due to the lack of relaxation 
or diffusion operator, very fine structures and filamentations can be observed which forbid the use 
of low order numerical methods. Hence a second order (in time $t$ and phase space $\xi$) finite difference discretization 
is applied to (\ref{eqFnum}), which is based on a Lax-Wendroff-Richtmyer numerical scheme (see \cite{richtmyer, gourlay}). 

First, we introduce the time discretization $t^n=n\Delta t$ with $n\in \mathbb{N}$ and the time step $\Delta t$. The phase space discretization is uniform so that the domain $[-\xi_{\max}, \xi_{\max}]^2$ 
is meshed by $\xi_{1,i}= -\xi_{\max} + i\Delta \xi$ and $\xi_{2,j} = -\xi_{\max} + j\Delta \xi$ for 
$i,j=0, \dots,  N-1$ 
and $\Delta \xi = 2\xi_{\max}/N$, $N$ being the number of points per direction. For the direction $\tau$, 
we also use a uniform mesh of size $\Delta \tau$, so that $\tau_\ell=\ell\Delta \tau$, for 
$\ell=0, \dots, N_\tau-1$, $\Delta \tau = 2\pi/N_\tau$. 
Denoting $\xi_{i,j}= (\xi_{1,i}, \xi_{2,j})$, the discrete unknown is then 
$F^n_{i,j,\ell} \approx F^\eps(t^n, \tau_\ell,  \xi_{i,j})$. 
In the following description, we keep the $\tau$ variable continuous in order to focus on the discretization 
in the $\xi_1$ and $\xi_2$ directions. In practice, since periodic boundary conditions are considered in this direction $\tau$, the fast Fourier transform is very efficient for this variable. At the boundary of the phase space domain in $\xi$, zero inflow boundary conditions are prescribed.

We then introduce the flux in $\xi$ which approximates $(E^n\cdot \nabla_\xi)F^n_{i,j}$ by centered finite differences:
$$
{\Phi}^n_{i,j}(F^n) = \frac{  E^n_{1,i+1,j} F^n_{i+1,j} - E^n_{1,i-1,j} F^n_{i-1,j}}{2\Delta \xi}+\frac{E^n_{2,i,j+1} F^n_{i,j+1} - E^n_{2,i,j-1} F^n_{i,j-1}}{2\Delta \xi},$$
and we also consider the following four-points average 
$$
\overline{F}^n_{i,j} =\left( F^n_{i+1,j}+F^n_{i-1,j}+F^n_{i,j+1}+F^n_{i,j-1}\right)/4.
$$

\noindent A first step on $\Delta t/2$ is performed to get intermediate unknowns $F^{n+1/2}_{i,j}$ 
\begin{equation}
\label{lf}
F^{n+1/2}_{i,j} = \overline{F}^n_{i,j} - \frac{\Delta t}{2}  {\Phi}^n_{i,j}(F^n)- \frac{\Delta t}{2\varepsilon}\partial_\tau F^{n+1/2}_{i,j}\,.
\end{equation}
The second step reads 
\begin{equation}
\label{lw}
F^{n+1}_{i,j} =  F^n_{i,j} - \Delta t 
{\Phi}^{n+1/2}_{i,j}(F^{n+1/2}) - \frac{\Delta t}{2\varepsilon}\partial_\tau (F^{n}_{i,j}+F^{n+1}_{i,j}).
\end{equation}
Standard results (see \cite{gourlay, richtmyer}) say that this numerical scheme is second order in time and phase space $\xi$ for all fixed $\varepsilon > 0$.

Recall now that the model is nonlinear due to the presence of the self-consistent electric field $\widetilde E_F$. Let us explain how we update the field $\widetilde E^{n+1}_F$, once $F^{n+1}$ is known. The inversion of the Poisson equation is easier in the original variables $(r,v)$ than in the variables $\xi$, since it takes the simple form \fref{poisson} of an ODE in the $r$ variable. At the continuous level, coming back to $(r,v)$ can be done easily by introducing the function
\be
\label{comeback}
f(t,\tau, r, v) = F(t, \tau,\xi_1,\xi_2),\qquad \mbox{with}\quad \left(\begin{array}{c}\xi_1\\ \xi_2\end{array}\right)=e^{-\tau J}\left(\begin{array}{c}r\\ v\end{array}\right).
\ee
It is not that simple at the discrete level. Indeed, if $(r_i,v_j)$ is the mesh in $(r,v)$, then, for all given $\tau_\ell$, the points $e^{-\tau_\ell J}\left(\begin{array}{c}r\\ v\end{array}\right)$ do not necessarily coincide with mesh points $\xi_{i,j}$. To evaluate $f^{n+1}(\tau_\ell, r_i, v_j) $, we thus need an interpolation algorithm in dimension 2. Since this interpolation is done at each step, we choose a simple linear interpolation algorithm. Then, once we known the values $f^{n+1}_{i,j}$ for each $\tau_\ell$, it is easy to compute the Poisson field $E^{n+1}_f$ by integrating \fref{poisson}. To deduce $\widetilde E^{n+1}_F$ on the $\xi$ mesh, another interpolation step is required. Finally, we also remark that our algorithm in two steps imposes to predict the advection field $E$ at time $t^{n+1/2}$, so a Poisson field evaluation is needed also before computing the flux $\Phi^{n+1/2}_{i,j}$.
 
At the final time $t_{final}$ of the simulation, we come back to the solution of our initial problem \fref{eqf} by setting $f(t_{final}, r, v)=F(t_{final},t_{final}/\eps,\xi)$, so a last interpolation algorithm in the two-dimensional $(r,v)$ variable is needed, as well as in the $\tau$ variable (since $t_{final}/\eps$ does not necessarily coincide with a discrete $\tau_\ell$).

\subsection{Asymptotic Preserving property}
\label{sectionAP}
In this subsection, we check formally that the numerical scheme presented above is Asymptotic Preserving, as announced. Thanks to the implicitation of the stiff term $\frac{1}{\eps}\pa_\tau$, the only stability condition will be a standard CFL condition of the form $\Delta t\leq C\Delta \xi$. In the sequel, we consider for simplicity that $\Delta t\sim \Delta \xi$. We have already seen that, for fixed $\eps>0$, this scheme is consistent (and of order 2) with the equation \fref{eqFnum}. We now have to examinate its behavior when $\eps\to 0$. 

It is convenient to analyse the asymptotics of numerical schemes written with the micro-macro decomposition technique, which was developed in \cite{ml,bennoune} as a flexible method in order to construct Asymptotic Preserving numerical schemes for collisional kinetic equations. Remark that, here, we have rewritten \fref{eqftilde} under the "collisional form" \fref{eqFintro} (the operator $\pa_\tau$ plays the role of the collision operator). The micro-macro method consists in mimicking the Chapman-Enskog expansion and decomposing the unknown $F^\eps$ into a macro part $G^\eps=\Pi F^\eps$ and the remainding micro part $h^\eps=(\I-\Pi)F^\eps$. This micro part is small when $\eps$ is small (but plays an important role when $\eps$ is not small, ensuring the AP property). In fact, our scheme \fref{lf}, \fref{lw} is already under a "micro-macro" form, thanks to the simple form of the operator $L=\pa_\tau$. Indeed, it suffices to set
$$G^n_{i,j}=\Pi F^n_{i,j},\qquad h^n_{i,j}=(\I-\Pi)F^n_{i,j}$$
to realize that our scheme is reformulated as follows:
\begin{equation}
\label{lfmm}
\begin{cases}
 G^{n+1/2}_{i,j} &= \overline{G}_{i,j}^n - \frac{\Delta t}{2} \Pi {\Phi}^n_{i,j}(G^n+ h^n), \\[3mm]
  h^{n+1/2}_{i,j} &= \overline{h}_{i,j}^n
- \frac{\Delta t}{2} (I-\Pi) {\Phi}^n_{i,j}(G^{n}+ h^n) 
- \frac{\Delta t}{2\varepsilon}\partial_\tau h^{n+1/2}_{i,j}\,,
\end{cases}
\end{equation}
\begin{equation}
\label{lwmm}
\begin{cases}
G^{n+1}_{i,j} &= G^n_{i,j} - \Delta t \, \Pi {\Phi}^{n+1/2}_{i,j}(G^{n+1/2}+ h^{n+1/2}).\\[3mm]
h^{n+1}_{i,j} &=  h^n_{i,j} - \Delta t 
(I-\Pi) {\Phi}^{n+1/2}_{i,j}(G^{n+1/2}+ h^{n+1/2}) - \frac{\Delta t}{2\varepsilon}\partial_\tau (h^{n}_{i,j}+h^{n+1}_{i,j})
\end{cases}
\end{equation}

We will proceed by an induction argument. From our choice \fref{defF0bis} of initial data, we deduce that
$$G^0=\mathcal O(1),\qquad h^0= -\varepsilon  L^{-1} (I-\Pi){\Phi}^0_{i,j}(G^0) + {\cal O}(\varepsilon^2).$$
Now, assume that we have proved that
$$
G^n=\mathcal O(1),\qquad h^n= -\varepsilon  L^{-1} (I-\Pi){\Phi}^n_{i,j}(G^n) + {\cal O}(\varepsilon^2+\eps\Delta t).
$$
On the one side, the micro part of the first step (\ref{lfmm}) gives 
$$
\left(I+\frac{\Delta t}{2\varepsilon}L\right) h^{n+1/2}_{i,j} = \overline{h}^n_{i,j} -\frac{\Delta t}{2} (I-\Pi) {\Phi}^n_{i,j}(G^{n}+h^n), 
$$
from which we deduce that 
\bea
h^{n+1/2}_{i,j} &=& -\varepsilon  L^{-1} (I-\Pi){\Phi}^n_{i,j}(G^n) + {\cal O}(\varepsilon^2)\nonumber\\
&=&
 - \varepsilon  L^{-1} (I-\Pi){\Phi}^{n+1/2}_{i,j}(G^{n+1/2}) + {\cal O}(\varepsilon^2+\eps \Delta t),\label{hnp12}
\eea
since $G^{n+1/2}=G^n+\mathcal O(\Delta t)$ and $E^{n+1/2}=E^n+\mathcal O(\Delta t)$.
On the other side, the micro part of \fref{lwmm} leads to 
\bea
h^{n+1}_{i,j} &=& -2\varepsilon  L^{-1} (I-\Pi){\Phi}^{n+1/2}_{i,j}(G^{n+1/2}) -h^n+ {\cal O}(\varepsilon^2)\nonumber\\
&=&-2\varepsilon  L^{-1} (I-\Pi){\Phi}^{n+1/2}_{i,j}(G^{n+1/2}) +\varepsilon  L^{-1} (I-\Pi){\Phi}^n_{i,j}(G^n)+ {\cal O}(\varepsilon^2)\nonumber\\
&=&-\varepsilon  L^{-1} (I-\Pi){\Phi}^{n+1}_{i,j}(G^{n+1})+ {\cal O}(\varepsilon^2+\eps\Delta t),\label{approxh}
\eea
which ends the induction proof. 

\medskip 

Let us now focus on the AP property. The macro part of (\ref{lfmm}) gives  
\bea
G^{n+1/2}_{i,j} &=& \overline{G}^n_{i,j} - \frac{\Delta t}{2} \Pi {\Phi}^n_{i,j}(G^n+h^n)\nonumber\\
&=&\overline{G}^n_{i,j} - \frac{\Delta t}{2} \Pi {\Phi}^n_{i,j}(G^n-\varepsilon  L^{-1} (I-\Pi){\Phi}^n_{i,j}(G^n)) +{\cal O}(\eps^2\Delta t+\eps\Delta t^2).\;\;\;\;\;\; \label{lflimit}
\eea
If we now insert (\ref{hnp12}) into the second equation of (\ref{lwmm}), we then obtain
\begin{eqnarray}
G^{n+1}_{i,j} \!\!\!\! &=& \!\!\! G^n_{i,j} -\Delta t \Pi {\Phi}_{i,j}^{n+1/2} \left(G^{n+1/2}+ h^{n+1/2}\right) \nonumber\\
&=&  \!\!\! G^n_{i,j} -\Delta t \Pi {\Phi}^{n+1/2}_{i,j} \left( G^{n+1/2}-  \varepsilon  L^{-1} (I-\Pi){\Phi}^{n+1/2}_{i,j}(G^{n+1/2})\right)+{\cal O}(\varepsilon^2\Delta t+\eps\Delta t^2). \nonumber\\
\label{lwlimit}
\end{eqnarray}
Passing to the limit as $\eps\to0$ (for fixed $\Delta \xi$, $\Delta t$) in \fref{lflimit}, \fref{lwlimit} yields
$$
\begin{cases} 
G^{n+1/2}_{i,j} &\ds =\overline{G}^n_{i,j} - \frac{\Delta t}{2} \Pi {\Phi}^n_{i,j}(G^n),\\[3mm]
 G^{n+1}_{i,j} &\ds = G^n_{i,j} -\Delta t \Pi {\Phi}^{n+1/2}_{i,j} \left( G^{n+1/2}\right),
\end{cases}
$$
which is a Lax-Wendroff-Richtmyer numerical discretization of the limit equation
$$\pa_t G+\Pi E\cdot \na_\xi G=0.$$
This proves that our scheme is Asymptotic Preserving. Furthermore, we observe that when $\varepsilon$ is small but not zero, up to a ${\cal O}(\varepsilon^2\Delta t+\eps\Delta t^2)$ remainder, the numerical scheme (\ref{lflimit}), (\ref{lwlimit}) is nothing but a second order Lax-Wendroff-Richtmyer numerical discretization for the approximate asymptotic equation \fref{eqGtilde} of $\widetilde G^\eps$. Hence, accumulating the errors will yield $\|G^n-\widetilde G^\eps\|_\infty\leq C\eps^2 + C\eps \Delta t\leq C\eps^2+C\Delta t^2$ (here $C$ denotes a generic constant independent of $\eps$, $\Delta t$ and $\Delta \xi$).

We have then, for all $n$,
\bee
h^{n}_{i,j} &=&- \varepsilon  L^{-1} (I-\Pi){\Phi}^{n}_{i,j}(G^{n})+ {\cal O}(\varepsilon^2+\Delta t^2)\\
&=&- \varepsilon  L^{-1} (I-\Pi)E\cdot \na_\xi \widetilde G^\eps+ {\cal O}(\varepsilon^2+\Delta t^2),
\eee
where we recall that we have assumed $\Delta t\sim \Delta \xi$.
Finally, in view of \fref{Fordre2}, we have $F^n=G^n+h^n=F^\eps+{\cal O}(\varepsilon^2+\Delta t^2)$. 
So far, this analysis concerns the asymptotics $\eps \to 0$. For a fixed $\eps>0$, we already know that our scheme is 
 of order two in time and space, which means that there exists a constant $K(\eps) >0$ only depending on $\eps$ and not on $\Delta t$ such that
 $\|F^n-F^\eps\| \leq K(\eps) \Delta t^2.$ These two behaviors can be summarized in the following estimate
 $$\|F^n-F^\eps\|_\infty \leq C \min \left(K(\eps) \Delta t^2, \eps^2+ \Delta t^2\right).$$
 
This means that our scheme is in fact a {\em second order} Asymptotic Preserving in the following sense:
\begin{itemize}
\item[--] for all fixed $\eps$, this scheme provides a second order approximation of the original equation \fref{eqF}; 
\item[--] when $\eps\to0$, this scheme degenerates into a second order approximation of the system \fref{Fordre2}, \fref{eqGtilde}, which itself approximates the original equation \fref{eqF} up to $\mathcal O(\eps^2)$ terms.
\end{itemize}

\subsection{Extension to the diffusion limit}
\label{sectiondiffusion}

The micro-macro decomposition is not only a tool to analyze the limit $\eps\to 0$ (as in \cite{ml,liu-mieussens}), it is also a practical method that allows to extend the construction of AP schemes to more complicated situations, see \cite{bennoune,relaxation,clemou,bord} for instance for collisional kinetic problems. Let us briefly present another oscillatory example that will be developed in a future work. For simplicity, we present this example in the linear setting of subsection \ref{sectionlinear}. We still consider \fref{eqftildelinear}, \fref{Alinear} but assume now that the average of $E(\tau,\xi)$ in $\tau$ vanishes: $\Pi E\equiv 0$ (this is the case for the paraxial beam model if the forcing term $a(\tau)$ has no Fourier component in the frequencies $0$, $2$ or $-2$). Then, the limit field $E^{(0)}$ in \fref{eqGtildelinear} vanishes and it is convenient to rescale the time variable in order to get a non trivial model at the limit. This amounts to considering, from the beginning, the so-called "diffusion scaling" of \fref{eqf} (even if the final model here will not contain any second order derivative):
$$
\partial_t f^\eps + \frac{v}{\varepsilon^2}\partial_x f^\eps +\left(a\left(\frac{t}{\varepsilon}\right) \frac{r}{\eps}- \frac{r}{\varepsilon^2} \right) \partial_{v} f^\eps =0.
$$
In this case, the associated equation in $F$ takes the following form, where the variable $\tau$ stands for $t/\eps^2$:
\be
\label{eqFdiffusion}
\pa_t F^\eps+\frac{1}{\eps}E(\tau,\xi)\cdot \nabla_\xi F^\eps=-\frac{1}{\eps^2}\pa_\tau F^\eps.
\ee
Our micro-macro scheme for \fref{eqFdiffusion} will consist in decomposing the discrete unknown as $F^n_{i,j}=G^n_{i,j}+h^n_{i,j}$, where the macro part $G^n_{i,j}=\Pi F^n_{i,j}$ and the micro part $h^n_{i,j}=(\I-\Pi)F^n_{i,j}$ are calculated by
\begin{equation}
\label{lfdiffusion}
\begin{cases}
 G^{n+1/2}_{i,j} &= \overline{G}^n_{i,j} - \frac{\Delta t}{2\eps} \Pi {\Phi}_{i,j}(h^n), \\[3mm]
  h^{n+1/2}_{i,j} &= \overline{h}^n_{i,j}
- \frac{\Delta t}{2\eps} (I-\Pi) {\Phi}_{i,j}(G^{n+1/2}+ h^n) 
- \frac{\Delta t}{2\varepsilon^2}\partial_\tau h^{n+1/2}_{i,j}\,,
\end{cases}
\end{equation}
\begin{equation}
\label{lwdiffusion}
\begin{cases}
G^{n+1}_{i,j} &= G^n_{i,j} - \frac{\Delta t}{\eps} \, \Pi {\Phi}_{i,j}(h^{n+1/2}),\\[3mm]
h^{n+1}_{i,j} &=  h^n_{i,j} - \frac{\Delta t}{\eps} 
(I-\Pi) {\Phi}_{i,j}(\frac{1}{2}(G^{n+1}+G^{n})+ h^{n+1/2}) - \frac{\Delta t}{2\varepsilon^2}\partial_\tau (h^{n}_{i,j}+h^{n+1}_{i,j}).
\end{cases}
\end{equation}
Let us briefly discuss the limit of this scheme as $\eps\to 0$. Since, initially, one has $h^0=\mathcal O(\eps)$ (see the discussion in section \ref{discuinit}), it is readily seen that our semi-implicit scheme will propagate this property. For all $n$, one has $h^n=\mathcal O(\eps)$, so the flux terms in the equations for $G$ in \fref{lfdiffusion} and in \fref{lwdiffusion} are not singular. The last equation implies that 
$$h^{n+1}_{i,j}= - \varepsilon  L^{-1} (I-\Pi){\Phi}_{i,j}(G^{n+1})+\mathcal O(\eps^2)$$
if this property holds true at step $n$. Hence, since it is true at step $n=0$, it holds true for all $n$. Consequently, one deduces successively from the three first equations of our scheme \fref{lfdiffusion}, \fref{lwdiffusion} that
\be
\label{diff1}
G^{n+1/2}_{i,j} = \overline{G}_{i,j}^n + \frac{\Delta t}{2} \Pi {\Phi}_{i,j}(L^{-1} (I-\Pi){\Phi}(G^{n}))+\mathcal O(\eps),
\ee
$$h^{n+1/2}_{i,j}= - \varepsilon  L^{-1} (I-\Pi){\Phi}_{i,j}(G^{n+1/2})+\mathcal O(\eps^2)$$
and
\be
\label{diff2}
G^{n+1}_{i,j}= G^n_{i,j} + \Delta t \, \Pi {\Phi}_{i,j}(L^{-1} (I-\Pi){\Phi}(G^{n+1/2}))+\mathcal O(\eps).
\ee
Finally, if we disgard the remainders $\mathcal O(\eps)$, the limit scheme \fref{diff1}, \fref{diff2} is a Lax-Wendroff-Richtmyer scheme for the limit equation for $G$:
$$\pa_tG-\Pi \left(E\cdot \na_\xi(L^{-1} (I-\Pi)E\cdot \na_\xi G)\right)=0.$$
The scheme \fref{lfdiffusion}, \fref{lwdiffusion} is thus Asymptotic Preserving in the diffusion limit.


\section{Numerical results}
\label{numerics}

In this section, we present some numerical results for the paraxial beam model \fref{eqf}, \fref{poisson} described in the introduction. In particular, our aim is to validate the Asymptotic Preserving property of our scheme. For all the simulations, the function $a$ in the applied electric field $E_{\rm app}$ defined by \fref{eapp1} is chosen as $a(\tau)=\cos^2(2\tau)$. In the first series of tests, in subsection \ref{nonlinear}, we solve the complete Vlasov-Poisson model. Then, in subsection \ref{linear}, we restrict our study to the linear case when the Poisson field is set to zero, and where the asymptotic models (the limit model and its $\eps$-correction) are explicit and can be solved analytically, which provide some additional reference solutions for small $\eps$'s. 

\begin{figure}[!htbp]
\begin{center}
\begin{tabular}{@{}c@{}c@{}}
\includegraphics[width=7cm]{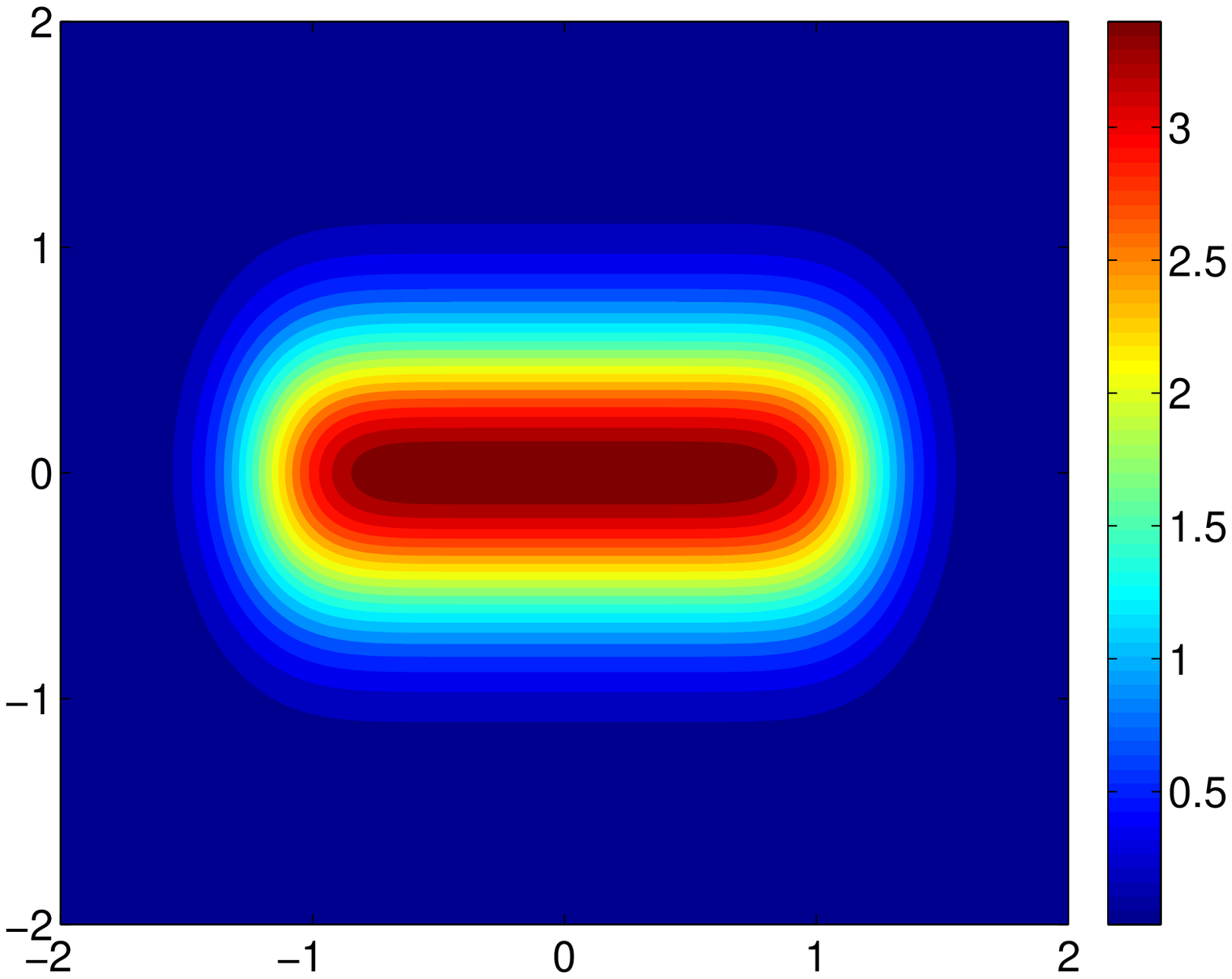} &     
\includegraphics[width=7cm]{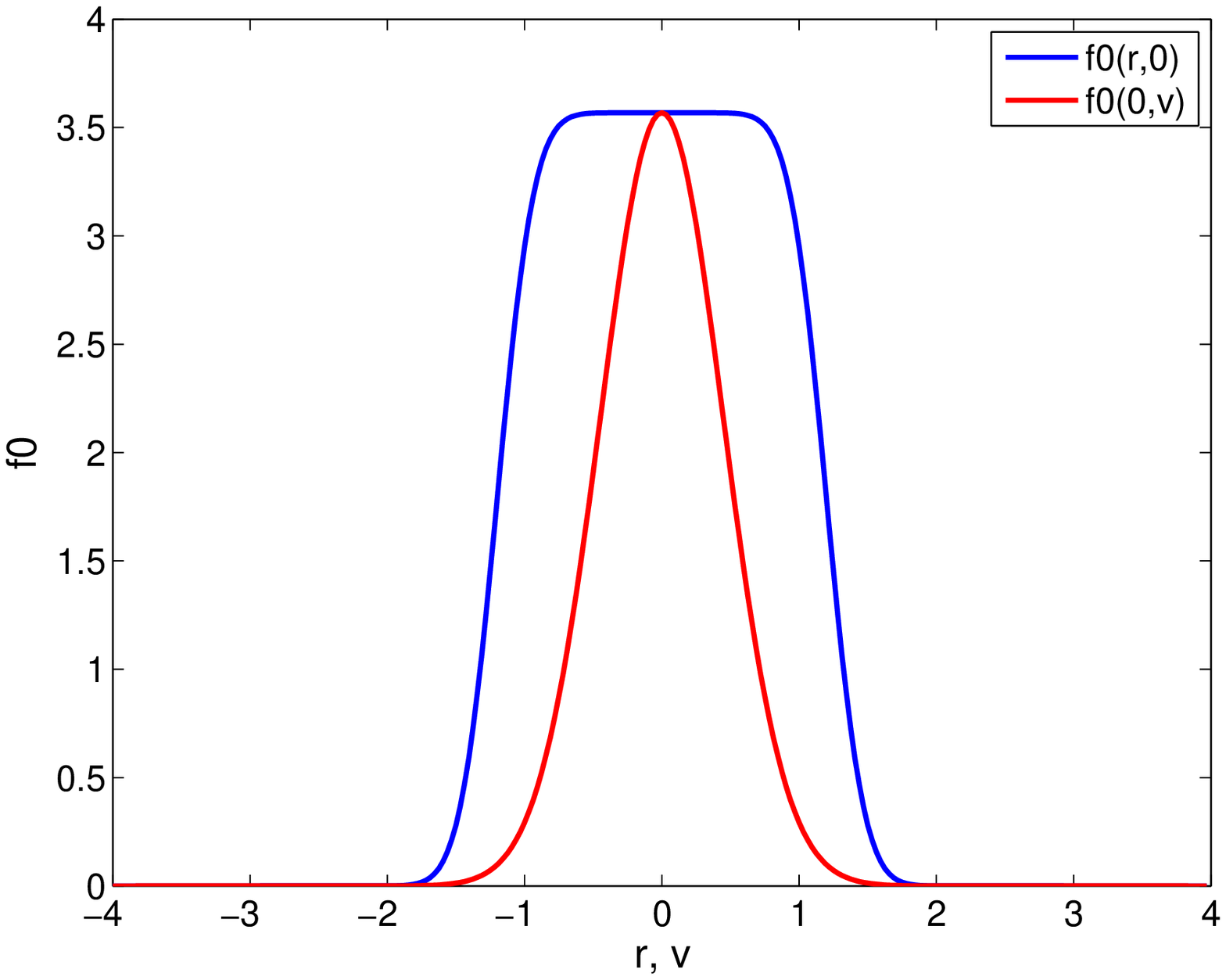}\end{tabular}
\caption{Plot of the initial data $f_0$. Left: 2D plot of the function in the $(r,v)$ space (zoomed for $(r,v)\in [-2,2]^2$). Right: the two curves $r\mapsto f_0(r,0)$ and $v\mapsto f_0(0,v)$.}
\label{figf0}
\end{center}
\end{figure}

The initial condition for \fref{eqf} is the same for all the simulations. It is taken as a Gaussian in velocity multiplied by a regularized step function in $r$:
\begin{equation}
\label{init0}
f_0(r, v) = \frac{4}{\sqrt{2\pi \alpha}} \chi(r)\,\exp\left(-\frac{v^2}{2\alpha} \right), \quad \chi(r)=\frac{1}{2}{\rm erf}\left(\frac{r+1.2}{0.3}\right)-\frac{1}{2}{\rm erf}\left(\frac{r-1.2}{0.3}\right)
\end{equation}
with $\alpha = 0.2$. For all the simulations, the space-velocity domain is $(r,v)\in [-4,4]^2$. We represent on Figure \ref{figf0} this initial data for $(r,v)\in [-2,2]^2$.

Let us list the numerical methods which are tested below:
\begin{itemize}
\item[--] our numerical scheme \fref{lf}, \fref{lw} with the initial data $F_0$ given by \fref{defF0bis}, containing the $\mathcal O(\eps)$ correction term, will be referred to as {\sl AP with correction};
\item[--] the same numerical scheme \fref{lf}, \fref{lw}, but with the initial data \fref{defF0ter}, without the correction term, will be referred to as {\sl AP without correction};
\item[--] a splitting method for the initial, non filtered equation, \fref{eqf}, \fref{poisson}: we apply a second order time-splitting method (Strang splitting) for \fref{eqf}, that we split into
$$\partial_t f^\eps + \frac{v}{\varepsilon}\partial_r f^\eps =0\qquad \mbox{and}\qquad \partial_t f^\eps +\left(E_{f^\eps} - \frac{r}{\varepsilon}+a\left(\frac{t}{\varepsilon}\right) r \right) \partial_{v} f^\eps =0,$$
each split equation being solved by a spectral method based on fast Fourier transform; this method will be referred to as the {\sl splitting scheme};
\item[--] in the linear case (see subsection \ref{linear}), we have the analytic expression \fref{lim1} for the exact solution $F_{\rm limit}$ of the limit model as $\eps\to 0$ --\,referred to as {\sl limit model}\,-- and we have \fref{lim2} for the solution $F_{\rm second order}$ of the limit model with the first correction in $\eps$ --\,referred to as {\sl second order model}.
\end{itemize}
For all the simulations, the number of discretization points in the $\tau$ direction is $N_\tau= 64$, hence the derivative $\pa_\tau$ and the integrals $\int_0^\tau$ are calculated with a spectral accuracy. The strategy for the choice of the time step is the following. For the two AP schemes, the time step is taken independently of $\eps$, it only has to satisfy the stability CFL condition related to our Lax-Wendroff-Richtmyer scheme, i.e. we always choose $\Delta t = \Delta \xi/\xi_{\max}\max{|E|}$, with $\xi_{\max}=4$ and $\Delta \xi=2\xi_{\max}/N$, $N$ being the number of points in the $\xi_1$ (or in the $\xi_2$) direction. For the {\sl splitting scheme}, we have to adapt $\Delta t$ proportionally to $\eps$. The {\sl limit model} and {\sl second order model} are analytic and do not require any time discretization.

\subsection{The Vlasov-Poisson model for the beam}
\label{nonlinear}

Our first series of simulations concern the full model \fref{eqf}, \fref{poisson} or its filtered equivalent version \fref{eqftilde}. 

\begin{figure}[!htbp]
\begin{center}
\begin{tabular}{@{}c@{}c@{}}
\includegraphics[width=7cm]{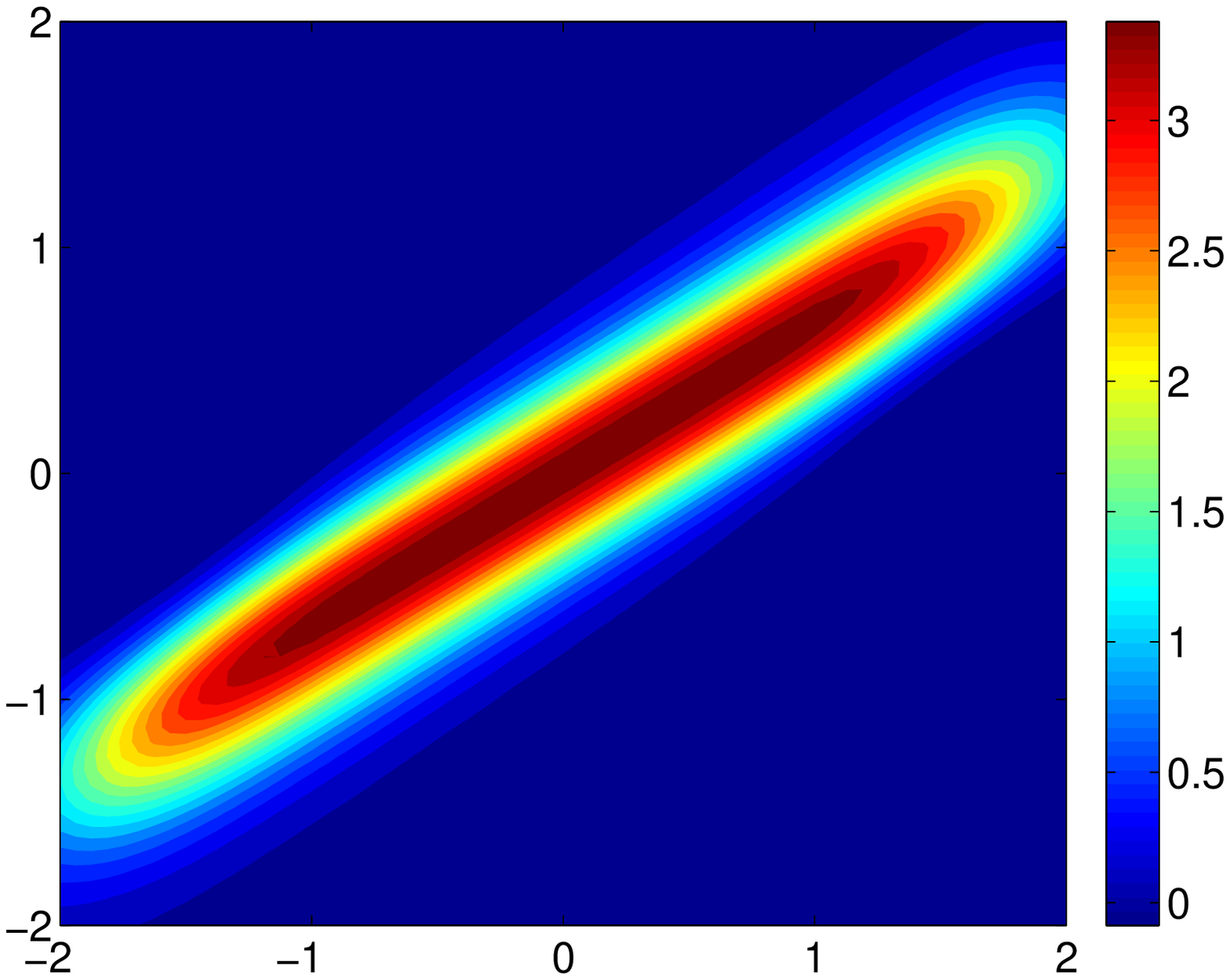} &     
\includegraphics[width=7cm]{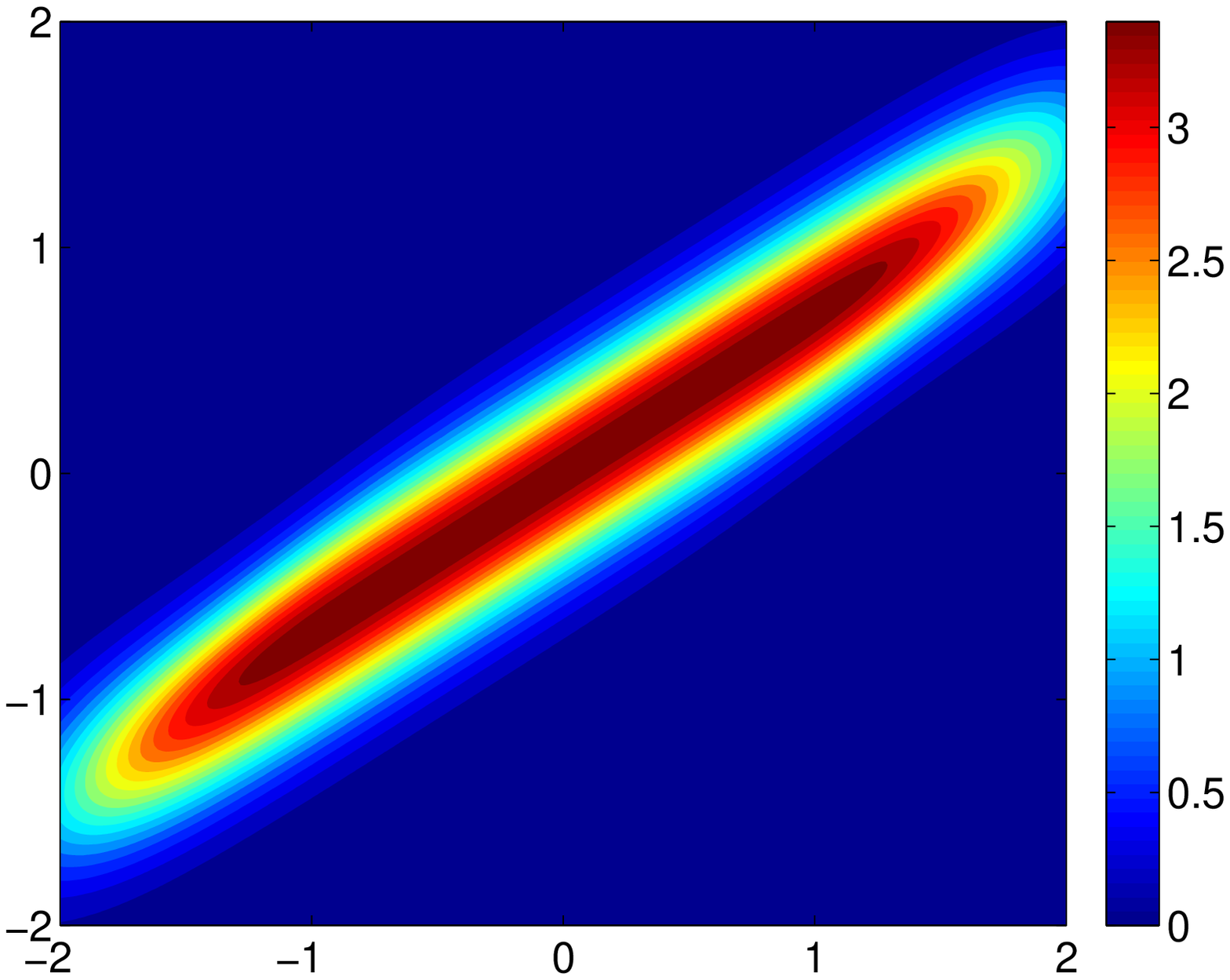}\\[-3mm]
\footnotesize $\eps=1$, {\sl AP with correction}&\footnotesize $ \eps=1$, {\sl splitting scheme}\\
\includegraphics[width=7cm]{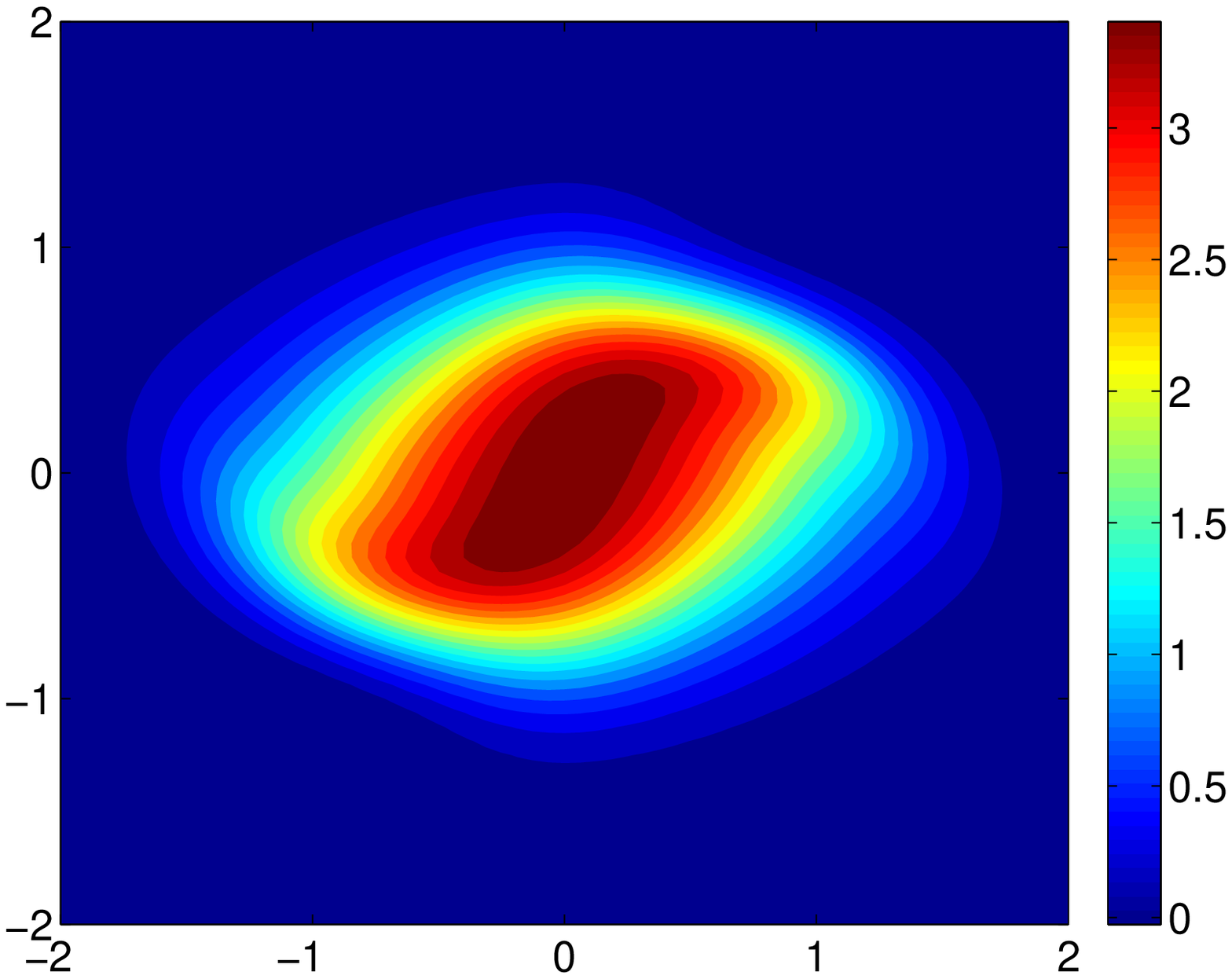} &     
\includegraphics[width=7cm]{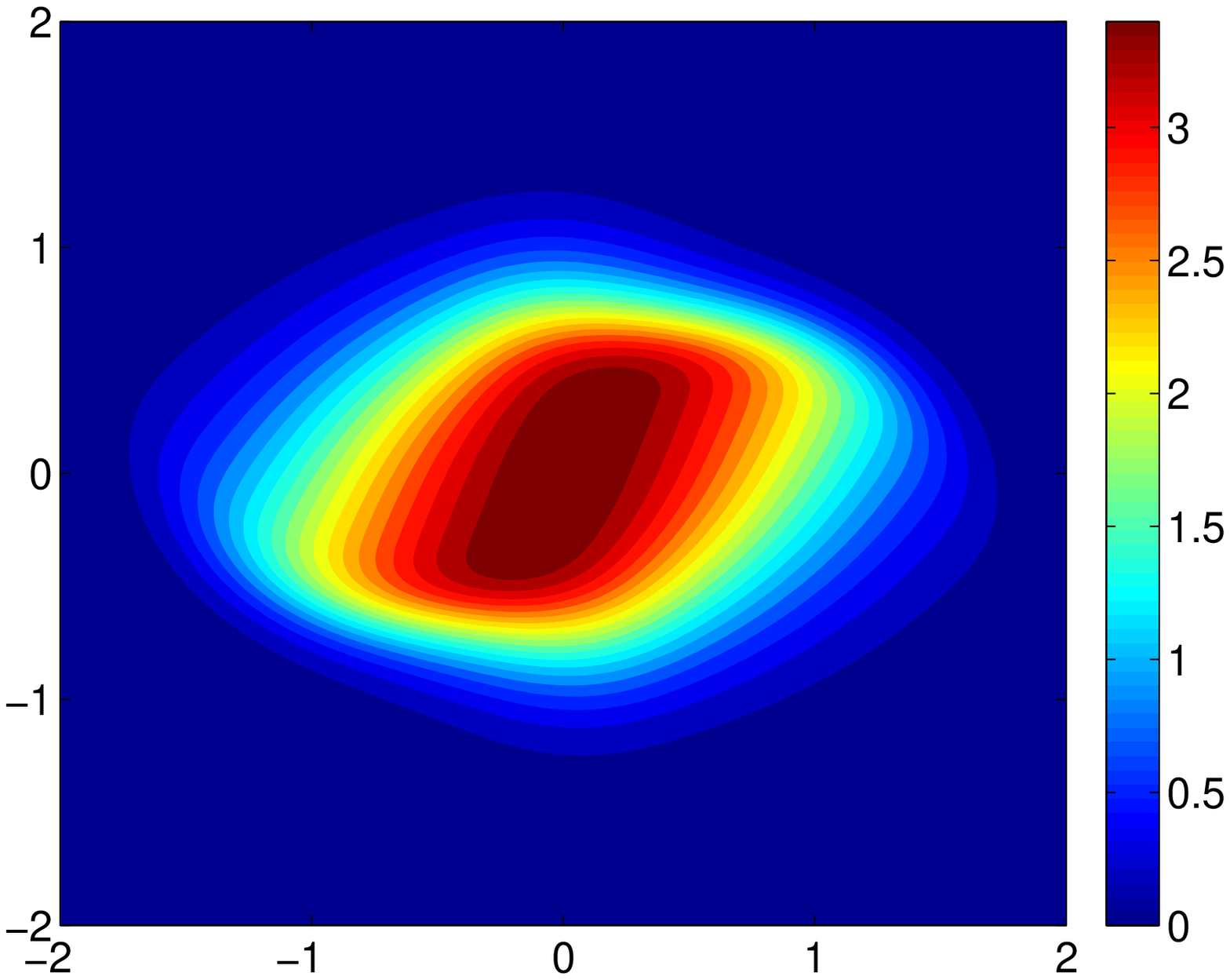}\\[-3mm]
\footnotesize $\eps=0.25$, {\sl AP with correction}&\footnotesize $\eps=0.25$, {\sl splitting scheme}\\
\includegraphics[width=7cm]{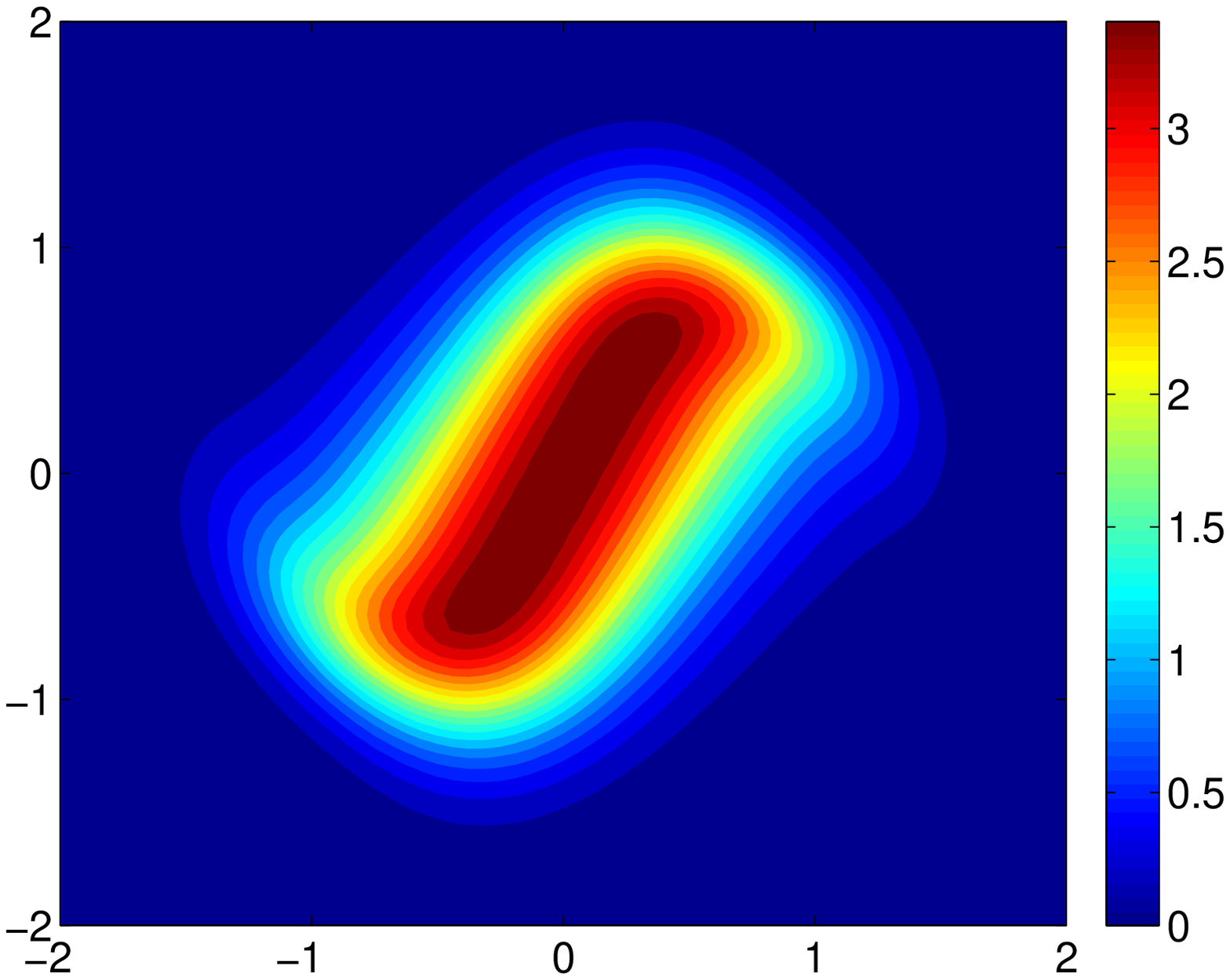} &     
\includegraphics[width=7cm]{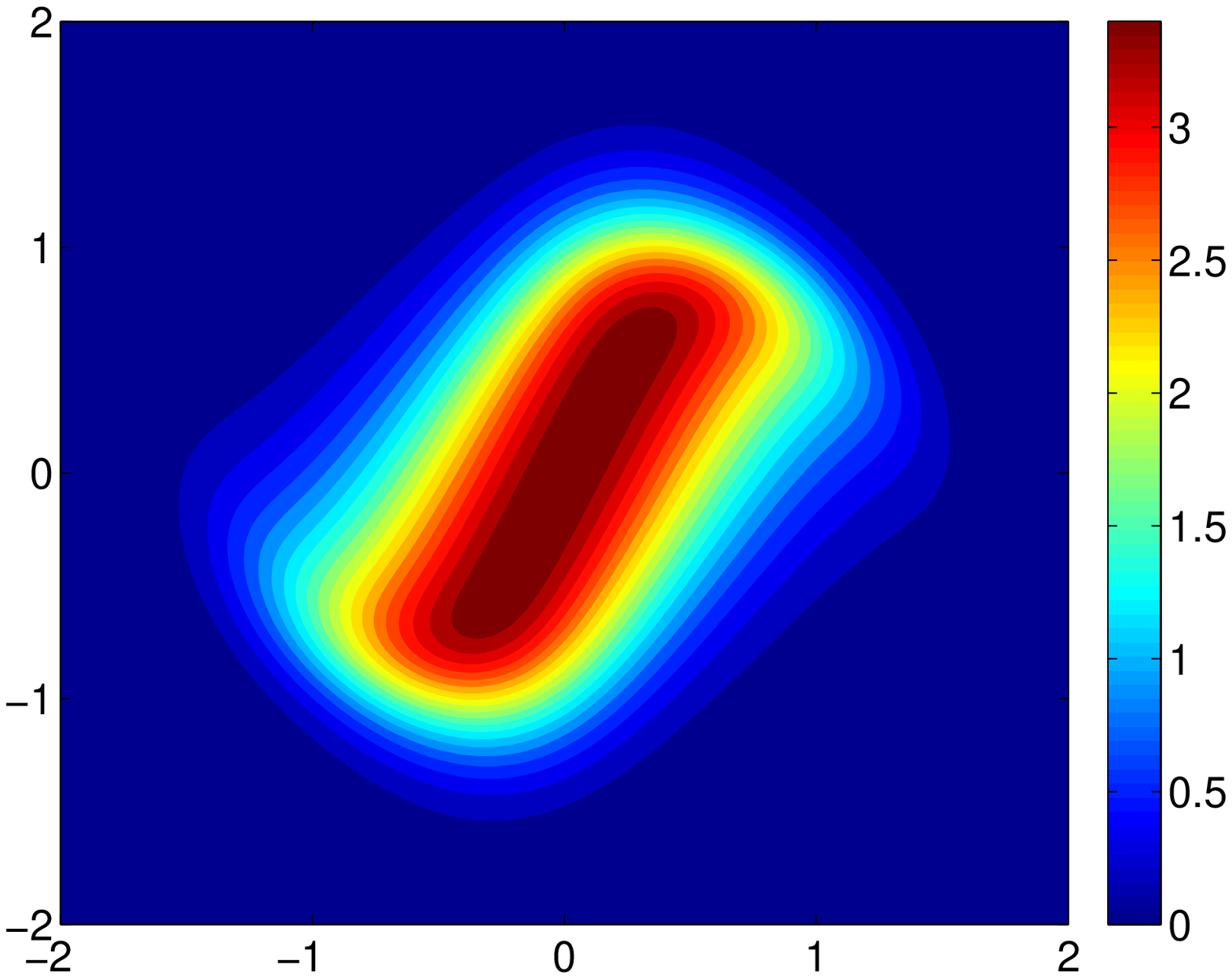}\\[-3mm]
\footnotesize $\eps=0.01$, {\sl AP with correction}&\footnotesize $\eps=0.01$, {\sl splitting scheme}
\end{tabular}
\caption{2D plots  for $(r,v)\in [-2,2]^2$ of the numerical solutions $f^\eps(t,r,v)$ at time $t=\pi/4$. Left column: computed with {\sl AP with correction}. Right column: computed with the {\sl splitting scheme}. Top line: $\eps=1$. Middle line: $\eps=0.25$. Bottom line: $\eps=0.01$.}
\label{figcontourPoi1}
\end{center}
\end{figure}

\bs
\ni
{\bf Qualitative results for different regimes in $\eps$}: Figure \ref{figcontourPoi1} 
\nopagebreak

\ms
\ni
Let us start with a few qualitative results. We first show some 2D plots of the function at the same final time $t_{final}=\pi/4$, for the three values $\eps=1$, $\eps=0.25$ and $\eps=0.01$. We compare in Figure \ref{figcontourPoi1} the numerical solution obtained by {\sl AP with correction} (here $N=128$), to the reference solution computed with the {\sl splitting scheme} with an adapted small time step. The time step for our AP scheme is $\Delta t=0.02$ for the three values of $\eps$. These plots show a good agreement between our solution and the reference solution: the scheme {\sl AP with correction} is able to capture all the regimes in $\eps$.

\begin{figure}[!htbp]
\begin{center}
\begin{tabular}{@{}c@{}c@{}}
\includegraphics[width=7cm]{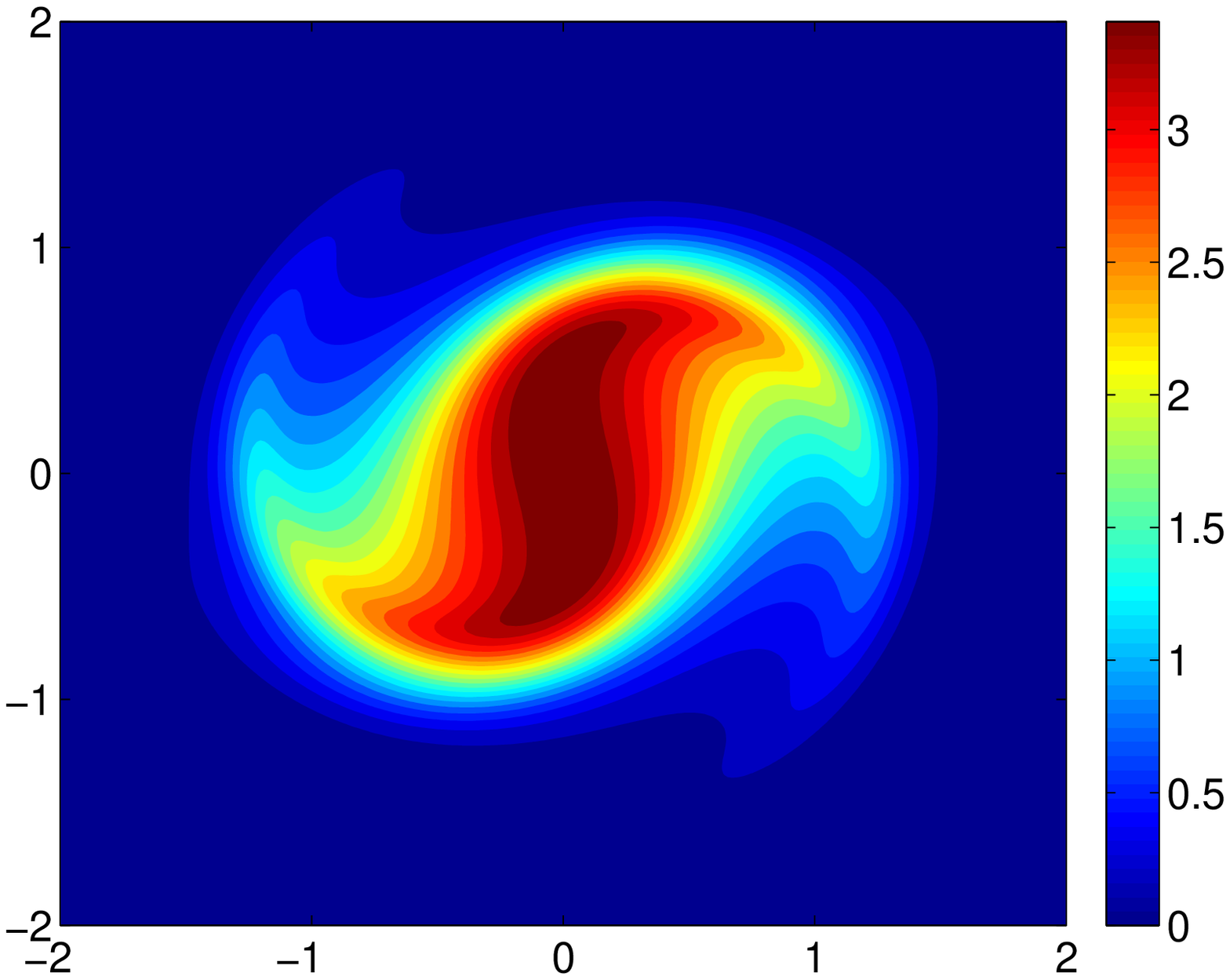} &     
\includegraphics[width=7cm]{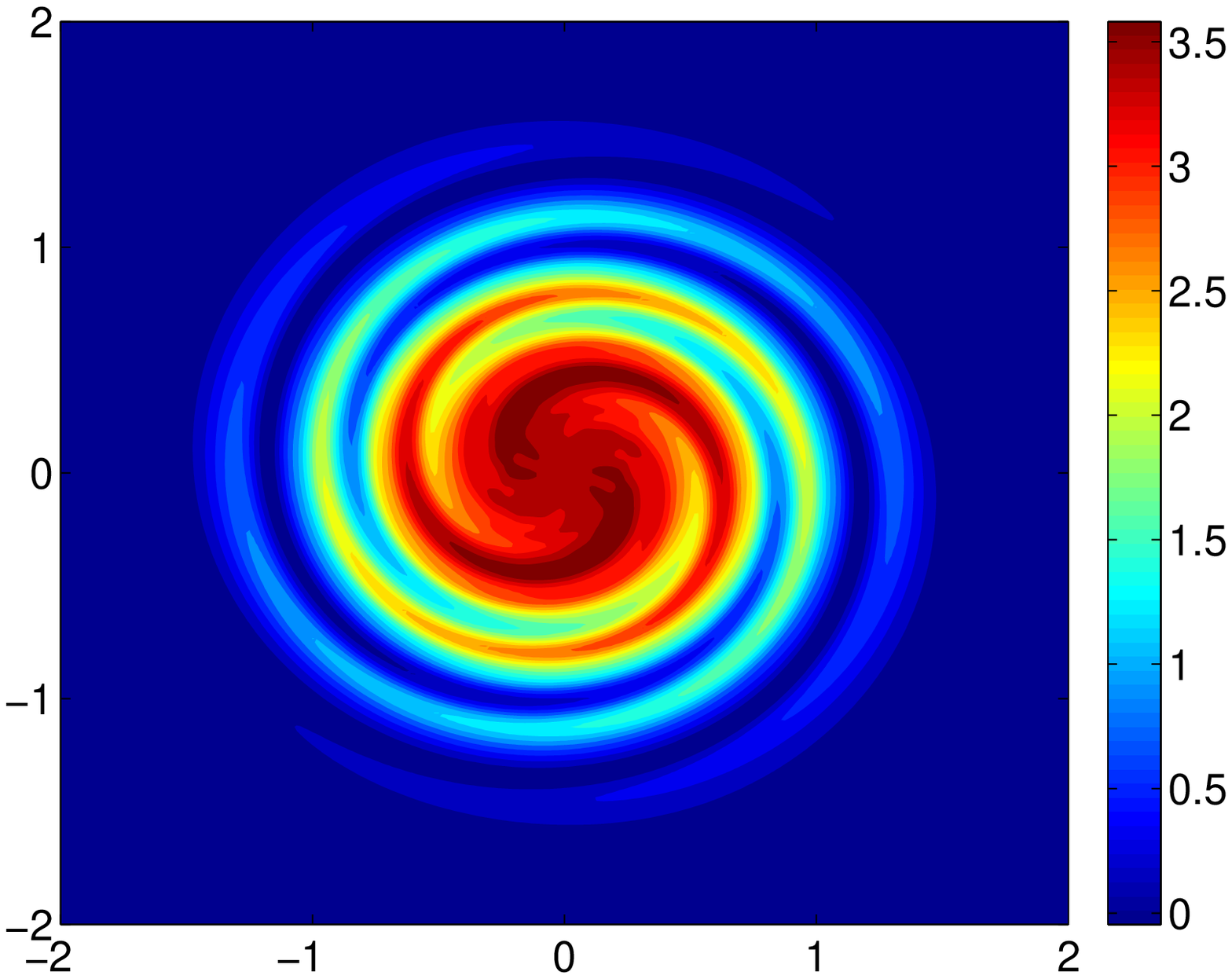}\\[-3mm]
\footnotesize $t=3.14$&\footnotesize $ t=12.6$\\
\includegraphics[width=7cm]{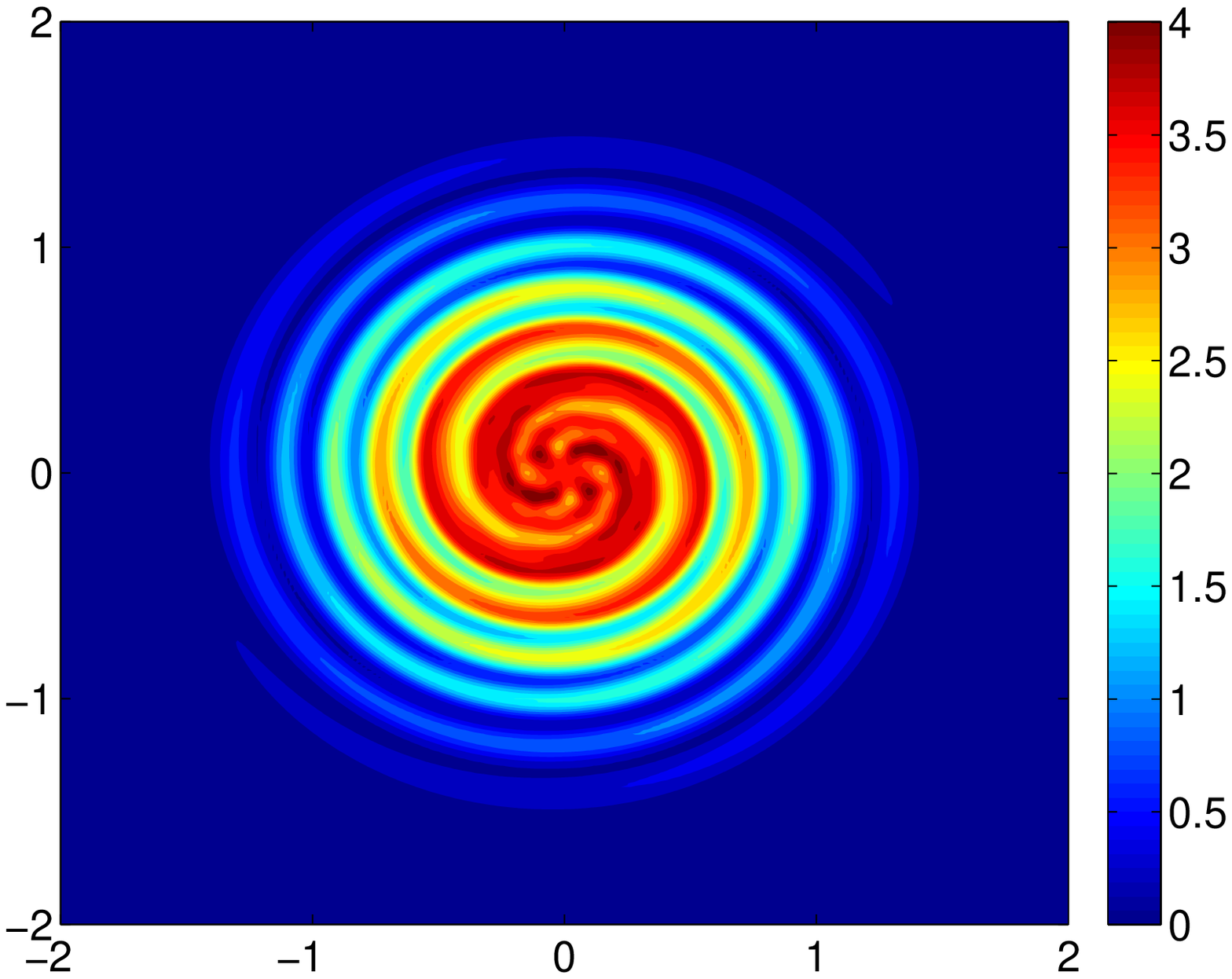} &     
\includegraphics[width=7cm]{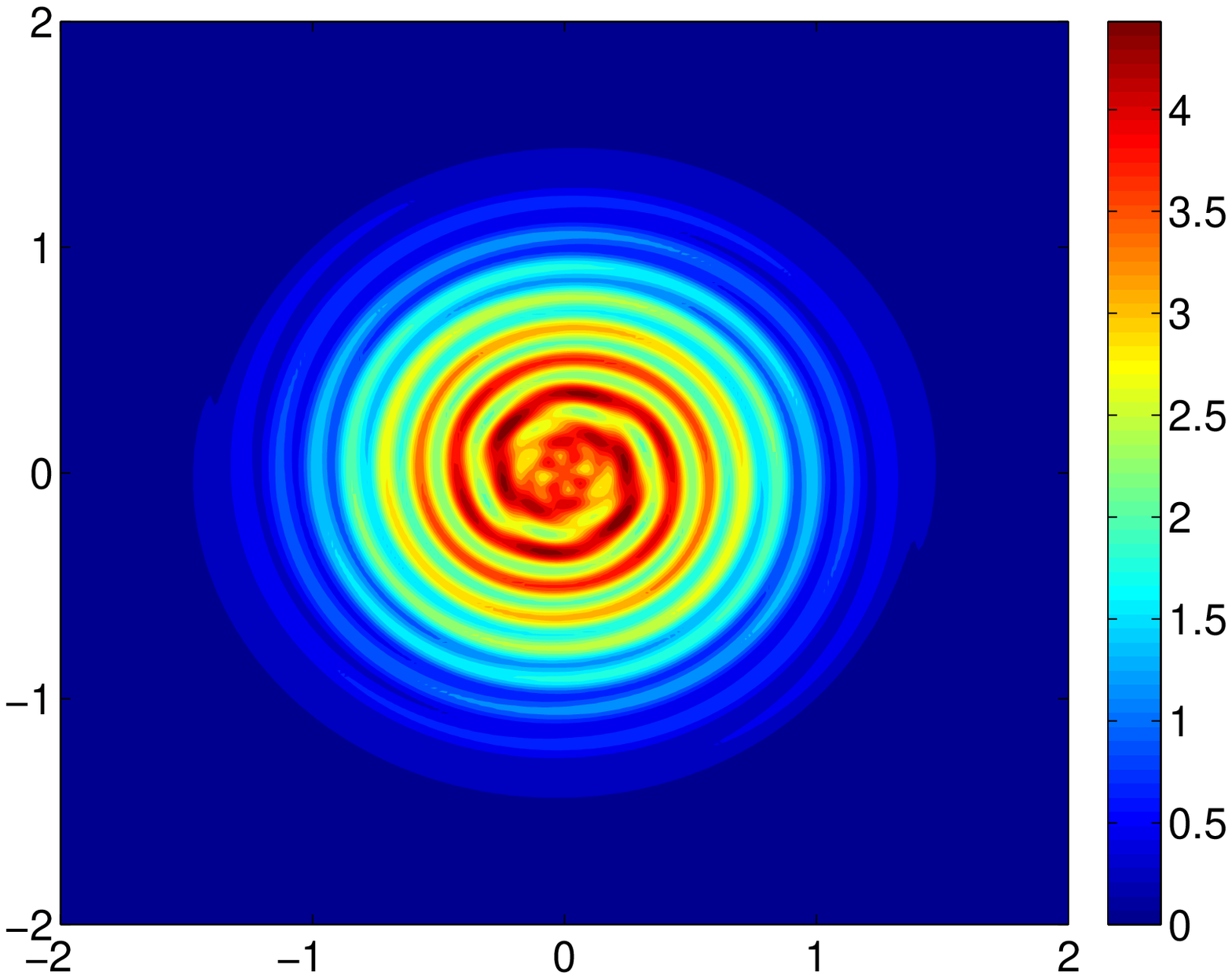}\\[-3mm]
\footnotesize $t=22$&\footnotesize $ t=31.4$\\
\end{tabular}
\caption{2D plots for $(r,v)\in [-2,2]^2$ of the numerical solutions $f^\eps(t,r,v)$ with the scheme {\sl AP with correction} for $\eps=0.001$, at times $t=\pi$, $t=4\pi$, $t=7\pi$ and $t=10\pi$. }
\label{figcontourPoi2}
\end{center}
\end{figure}

\bs
\ni
{\bf Long time behavior and filamentation}: Figure \ref{figcontourPoi2}
\nopagebreak

\ms
\ni
Now, we show that our AP scheme is able to capture very thin structures, with a numerical cost independent of $\eps$. On Figure \ref{figcontourPoi2} we plot the numerical solution obtained with the scheme {\sl AP with correction} (with $N=512$), for a very small $\eps=0.001$ and for different times $t=\pi$, $t=4\pi$, $t=7\pi$ and $t=10\pi$. We observe the filamentation due to the self-consistent Poisson field effect (compare to Figure \ref{figcontoursansPoi} below, obtained at $t=2\pi$ without the Poisson field).

\begin{figure}[!htbp]
\begin{center}
\hspace*{-15mm}
\begin{tabular}{@{}c@{}c@{}}
\includegraphics[width=10cm]{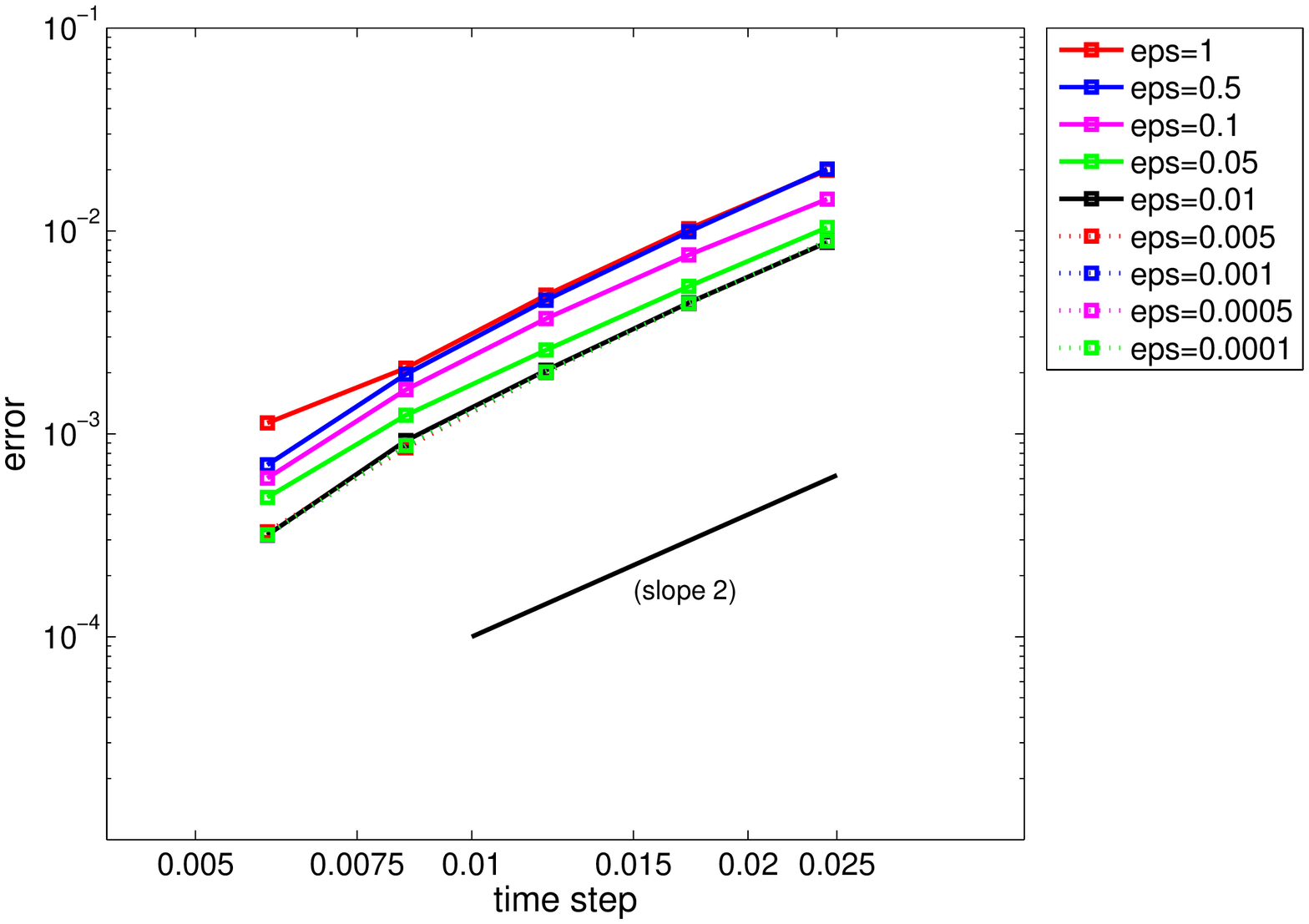}&
\hspace*{-7mm}\includegraphics[width=10cm]{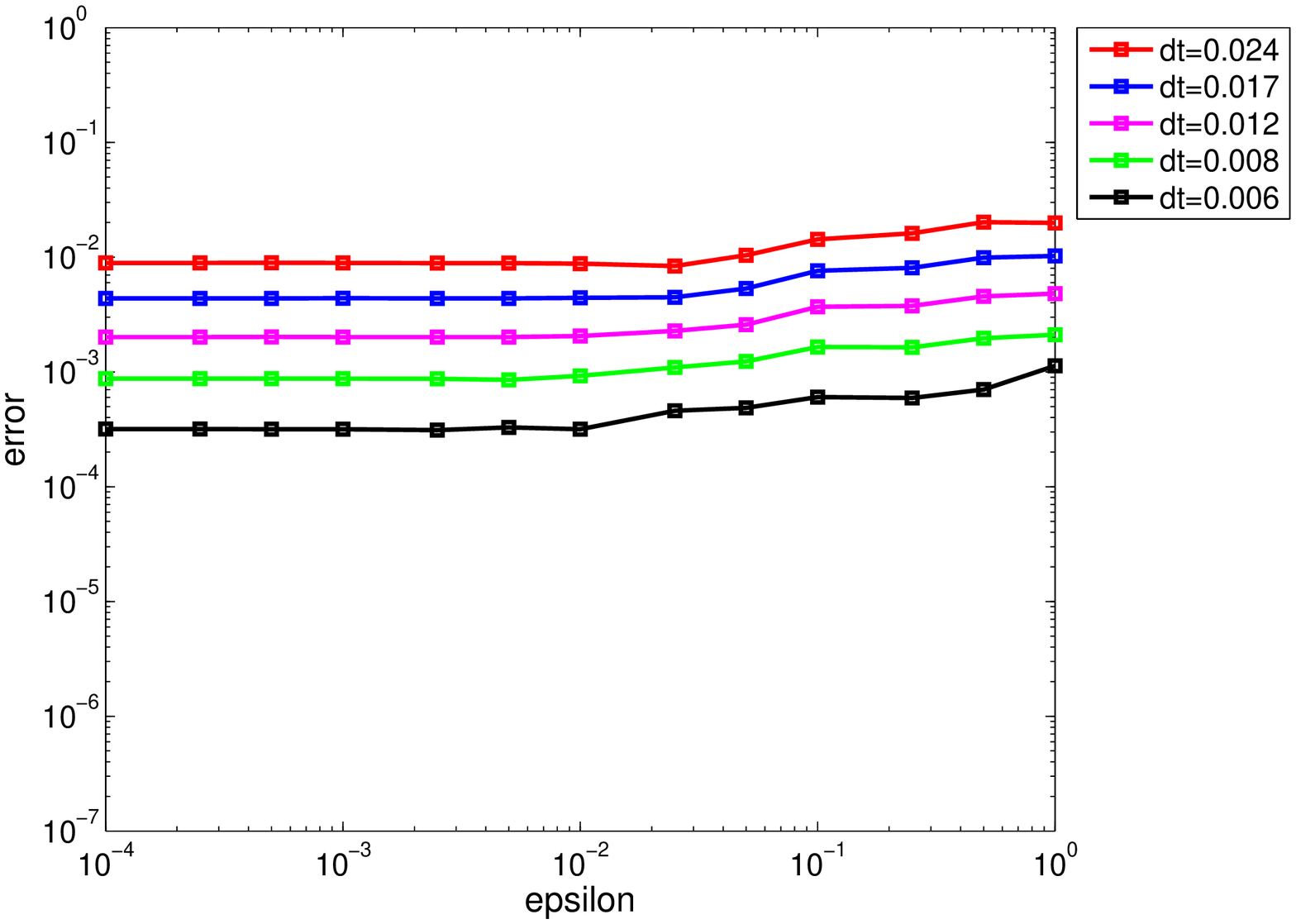}
\end{tabular}
\caption{Plot of the relative $L^2$ error for the scheme {\sl AP with correction}. Left: error as a function of $\Delta t$ for different $\eps$. Right: error as a function of $\eps$ for different $\Delta t$. Conclusion: the scheme is of order 2 and the error (nearly) does not depend on $\eps$.}
\label{figordrecorr}
\end{center}
\end{figure}

\begin{figure}[!htbp]
\begin{center}
\hspace*{-15mm}
\begin{tabular}{@{}c@{}c@{}}
\includegraphics[width=10cm]{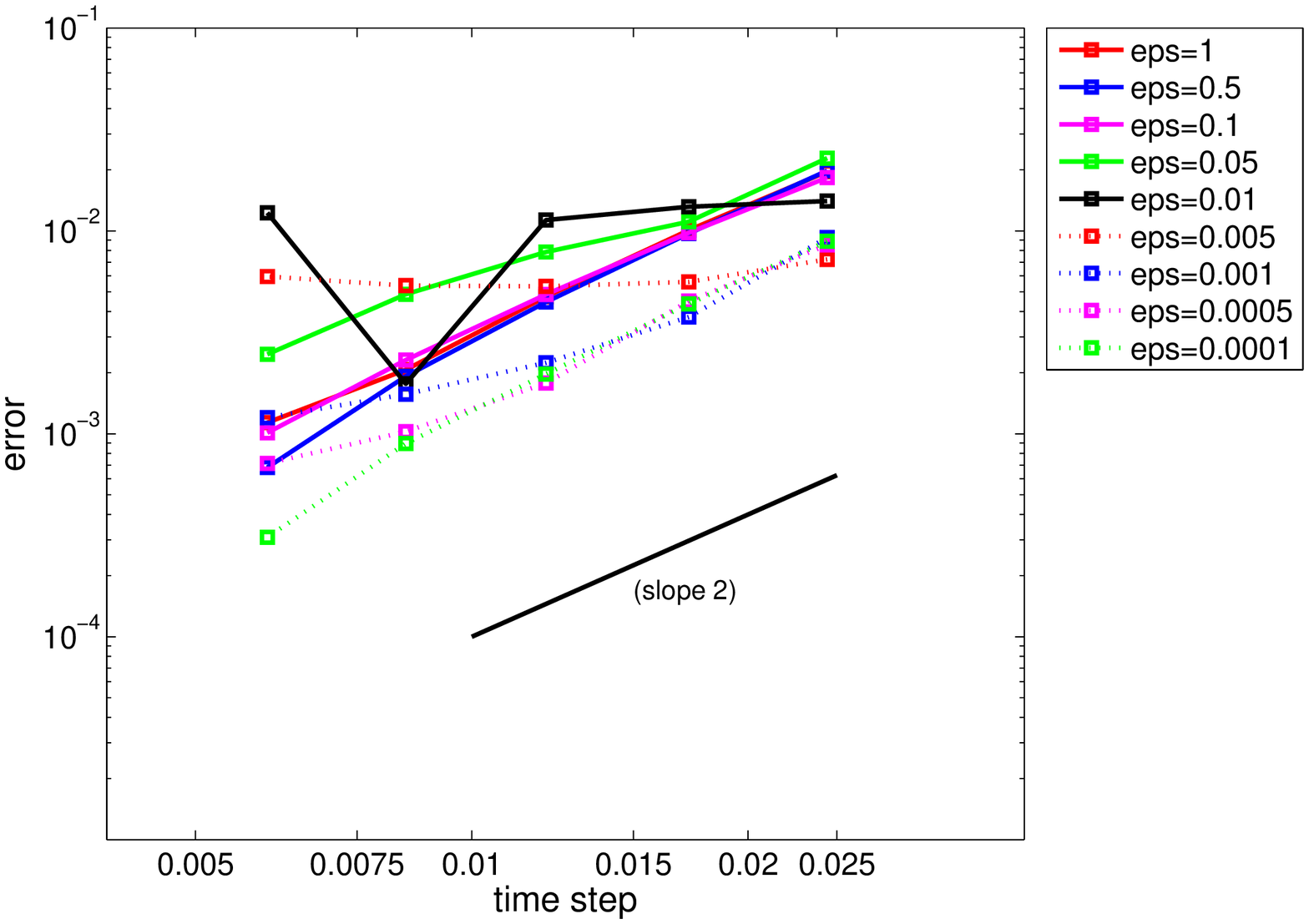}&
\hspace*{-7mm}\includegraphics[width=10cm]{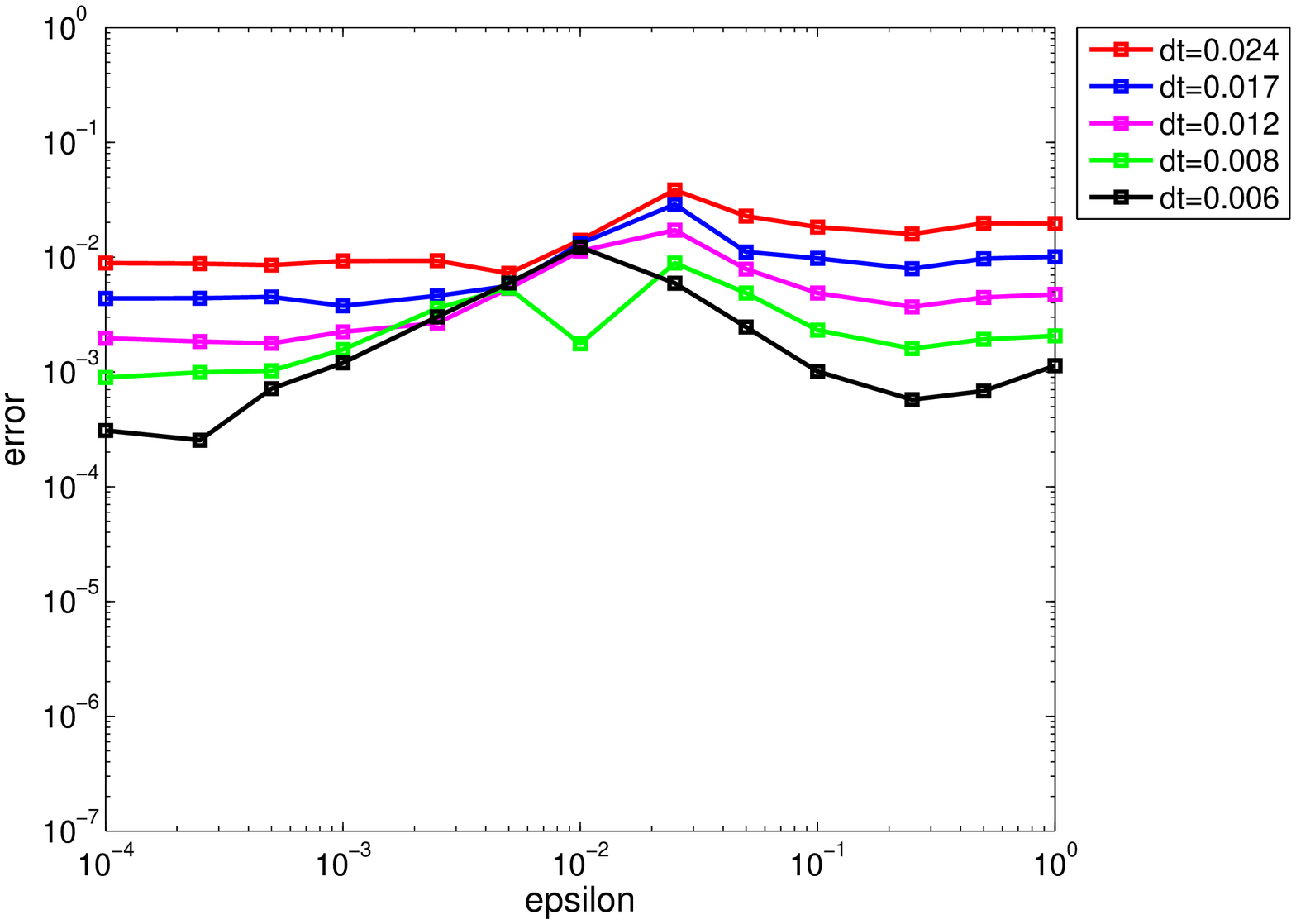}
\end{tabular}
\caption{Plot of the relative  $L^2$ error for the scheme {\sl AP without correction}. Left: error as a function of $\Delta t$ for different $\eps$. Right: error as a function of $\eps$ for different $\Delta t$.}
\label{figordresans}
\end{center}
\end{figure}

\bs
\ni
{\bf Numerical verification of the order 2 uniform accuracy with respect to $\eps$}: Figures \ref{figordrecorr},  \ref{figordresans} and  \ref{figordrespec}
\nopagebreak

\ms
\ni
Let us now proceed to more quantitative tests. We plot on the three next figures the  relative $L^2$ error between the numerical solutions computed with different schemes and a reference solution (computed with tiny time and space steps). The final time ($t=\pi/16$) is fixed. 

\ms
For Figure  \ref{figordrecorr}, the solution is computed with the scheme {\sl AP with correction}. On the left part, we represent (in logarithmic scales) the error as a function of the time step $\Delta t$, for different values of $\eps$ (from $\eps=1$ to $\eps=10^{-4}$): the slope is always close to 2 and the curves are very close together, indicating that the error is almost independent of $\eps$. This independence is confirmed on the right part of the figure, where we represent the error as a function of $\eps$, for different values of $\Delta t$: all the curves are nearly horizontal. These curves indicate that the error produced by the scheme {\sl AP with correction} is of the form $C\Delta t^2$, with $C$ independent of $\eps$. This proves experimentally the {\em second order} Asymptotic Preserving behavior of our scheme.

\ms
For Figure  \ref{figordresans}, the same tests are done for the scheme {\sl AP without correction}, i.e. for the scheme \fref{lf}, \fref{lw} with the initial data $F_0(\tau,\xi)=f_0(\xi)$. On the left part of the figure, we observe that the scheme behaves at an order 2 scheme for $\eps=\mathcal O(1)$ ($\eps=1$, 0.5 or 0.1) or for small values of $\eps$ (less than $10^{-3}$). But for intermediate regimes, the curves are more chaotic. On the right part of the figure, this feature is even more obvious: without the correction of the initial data, our scheme behaves well for $\eps=\mathcal O(1)$ and for $\eps$ very small (in fact, when the observed error is greater than $\eps$), but not for intermediate regimes. This shows that this initial correction is really needed and this validates numerically the analysis done in Section \ref{gene}.

\ms
For Figure  \ref{figordrespec}, the same tests are done for the {\sl splitting scheme} (well resolved in space, we only observe the error in the time step). On the left part of the figure, we observe that, for all fixed $\eps$, the Strang splitting scheme is of order 2 but the important fact is that the error strongly depends on $\eps$: the smaller is $\eps$, the smaller must be the time step to maintain a constant error. We also observe this feature on the right part of the figure. Experimentally, one can estimate that the error for the {\sl splitting scheme} is of the form $C(\Delta t/\eps)^2$.

\begin{figure}[!htbp]
\begin{center}
\hspace*{-15mm}
\begin{tabular}{@{}c@{}c@{}}
\includegraphics[width=10cm]{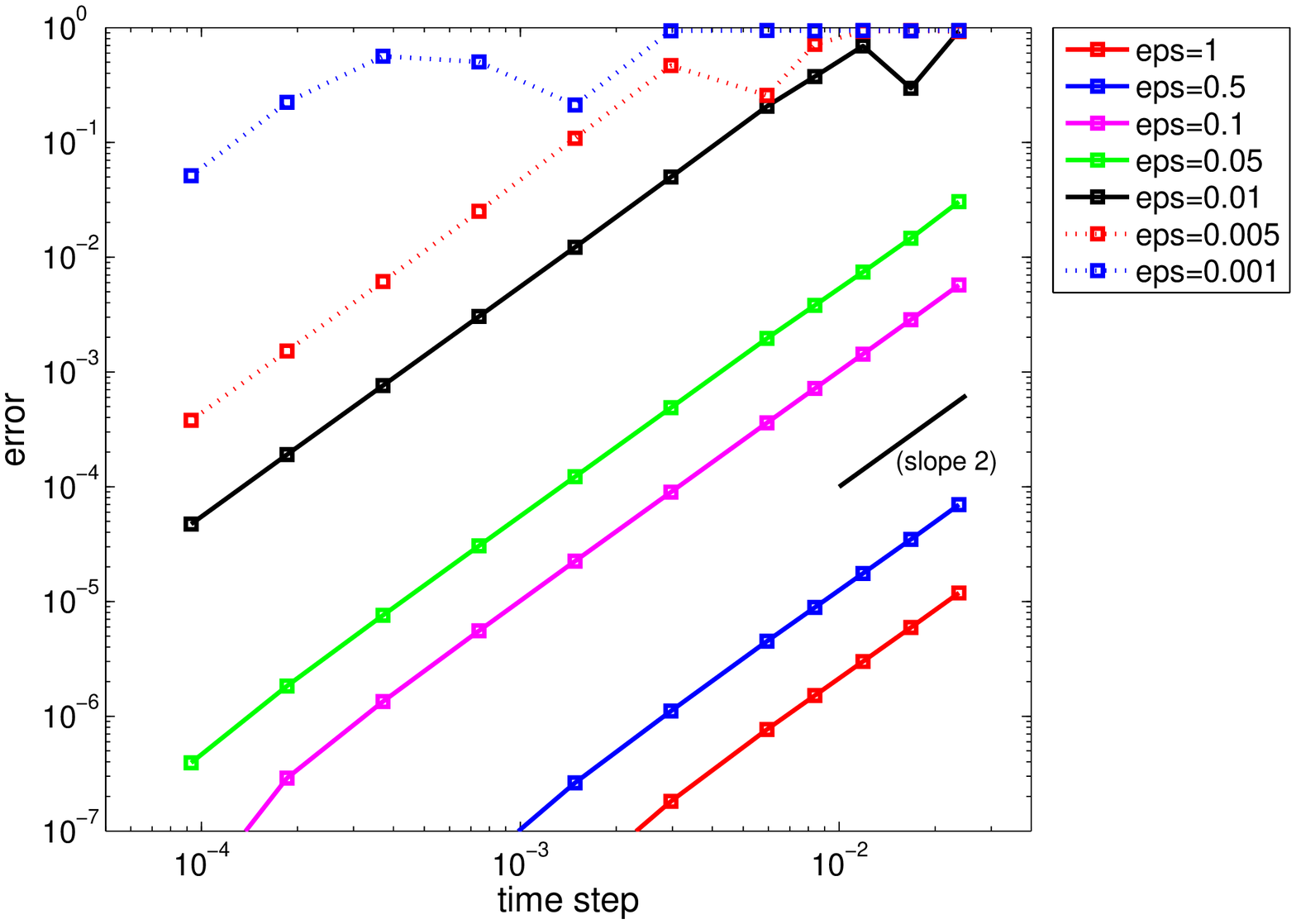}&
\hspace*{-7mm}\includegraphics[width=10cm]{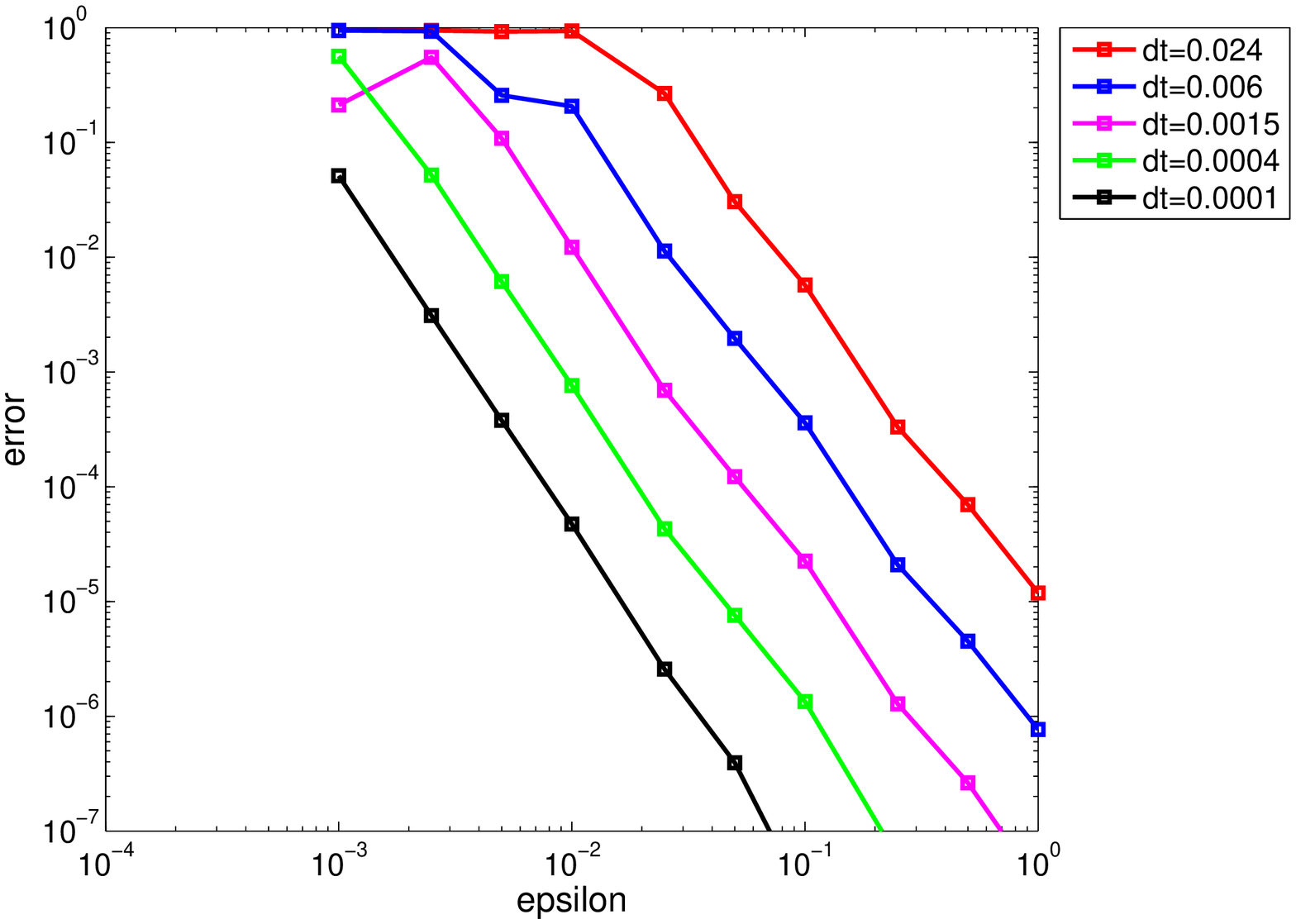}
\end{tabular}
\caption{Plot of the  relative $L^2$ error for the {\sl splitting scheme}. Left: error as a function of $\Delta t$ for different $\eps$. Right: error as a function of $\eps$ for different $\Delta t$. Conclusion: the error behaves like $C(\Delta t/\eps)^2$.}
\label{figordrespec}
\end{center}
\end{figure}

\begin{figure}[!htbp]
\begin{center}
\hspace*{-7mm}
\begin{tabular}{@{}c@{}c@{}}
\includegraphics[width=8cm]{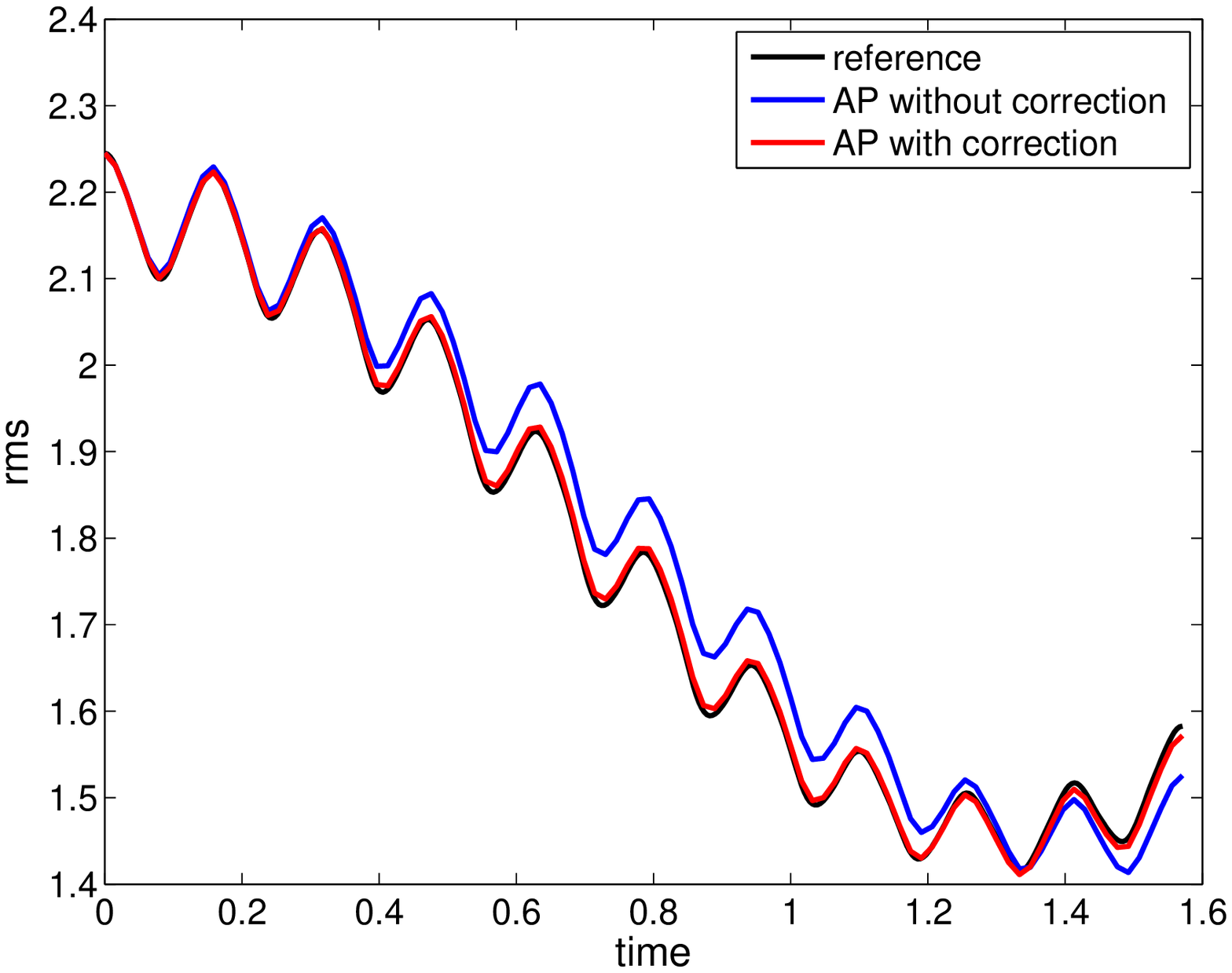} &     
\includegraphics[width=8cm]{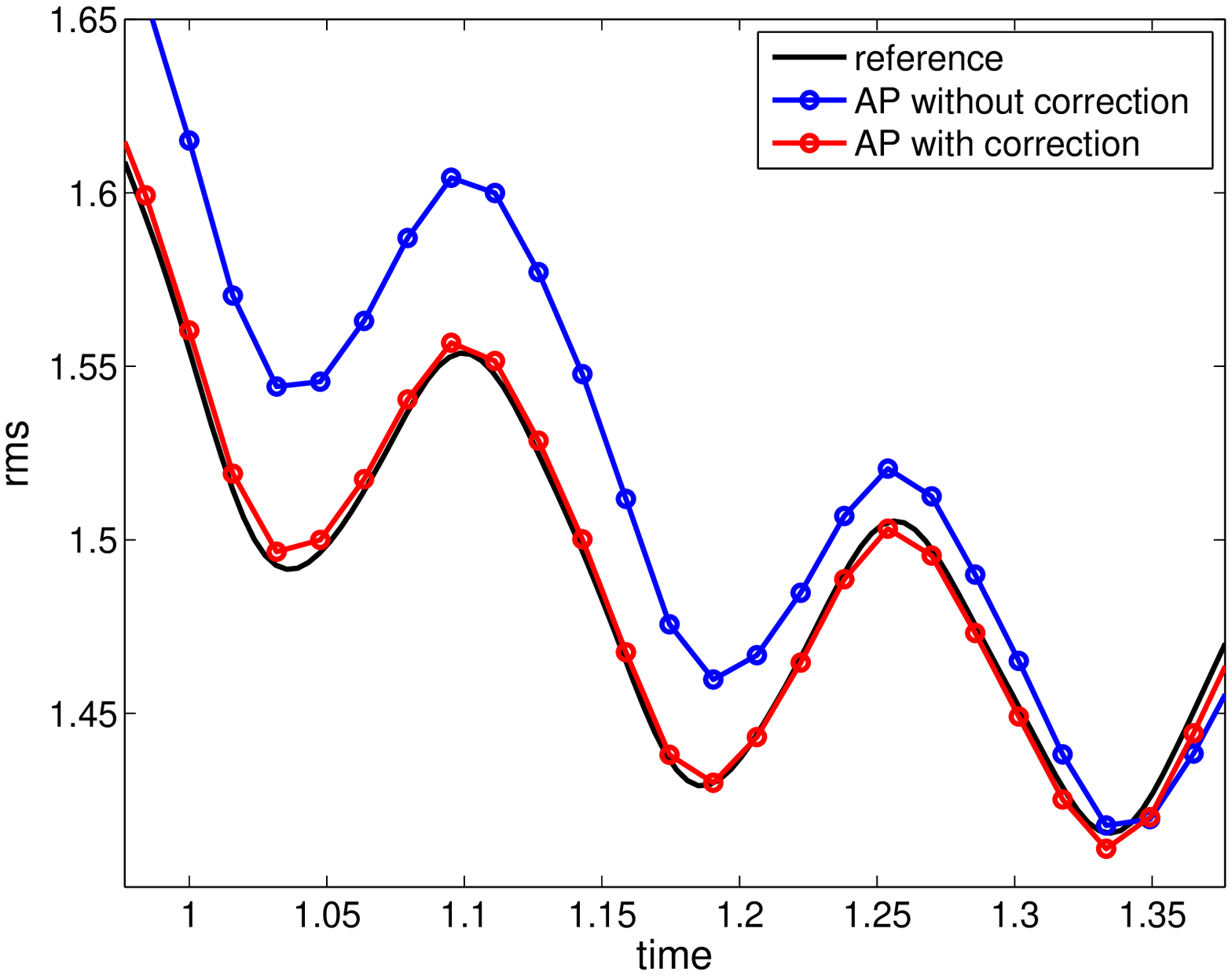}
\end{tabular}
\caption{Time history of $RMS(t)$ for $\varepsilon=0.05$, computed by {\sl AP with correction} and {\sl AP without correction}. On the right: zoom of the left figure.}
\label{figrms1}
\end{center}
\end{figure}

\begin{figure}[!htbp]
\begin{center}
\hspace*{-7mm}
\begin{tabular}{@{}c@{}c@{}}
\includegraphics[width=8cm]{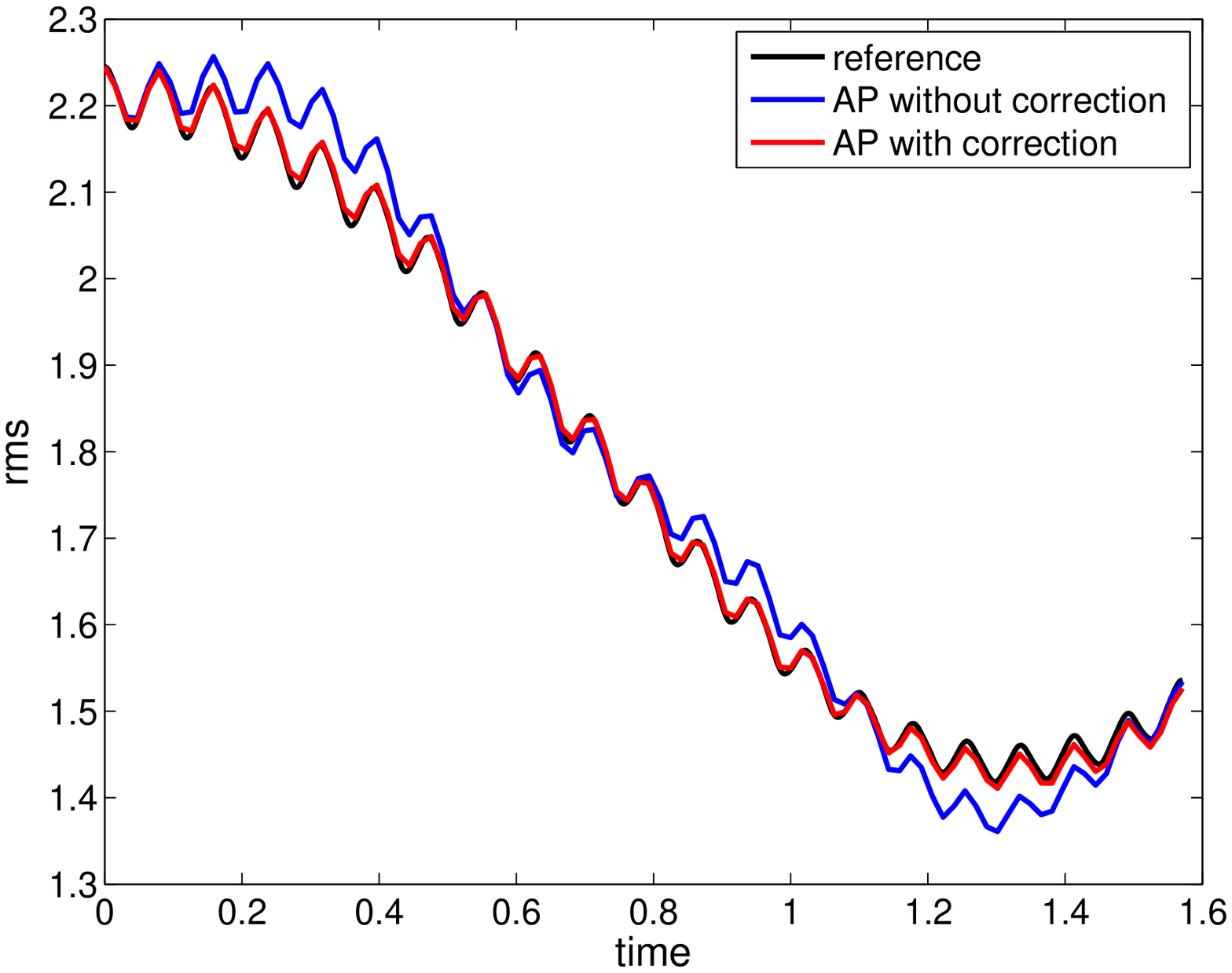} &     
\includegraphics[width=8cm]{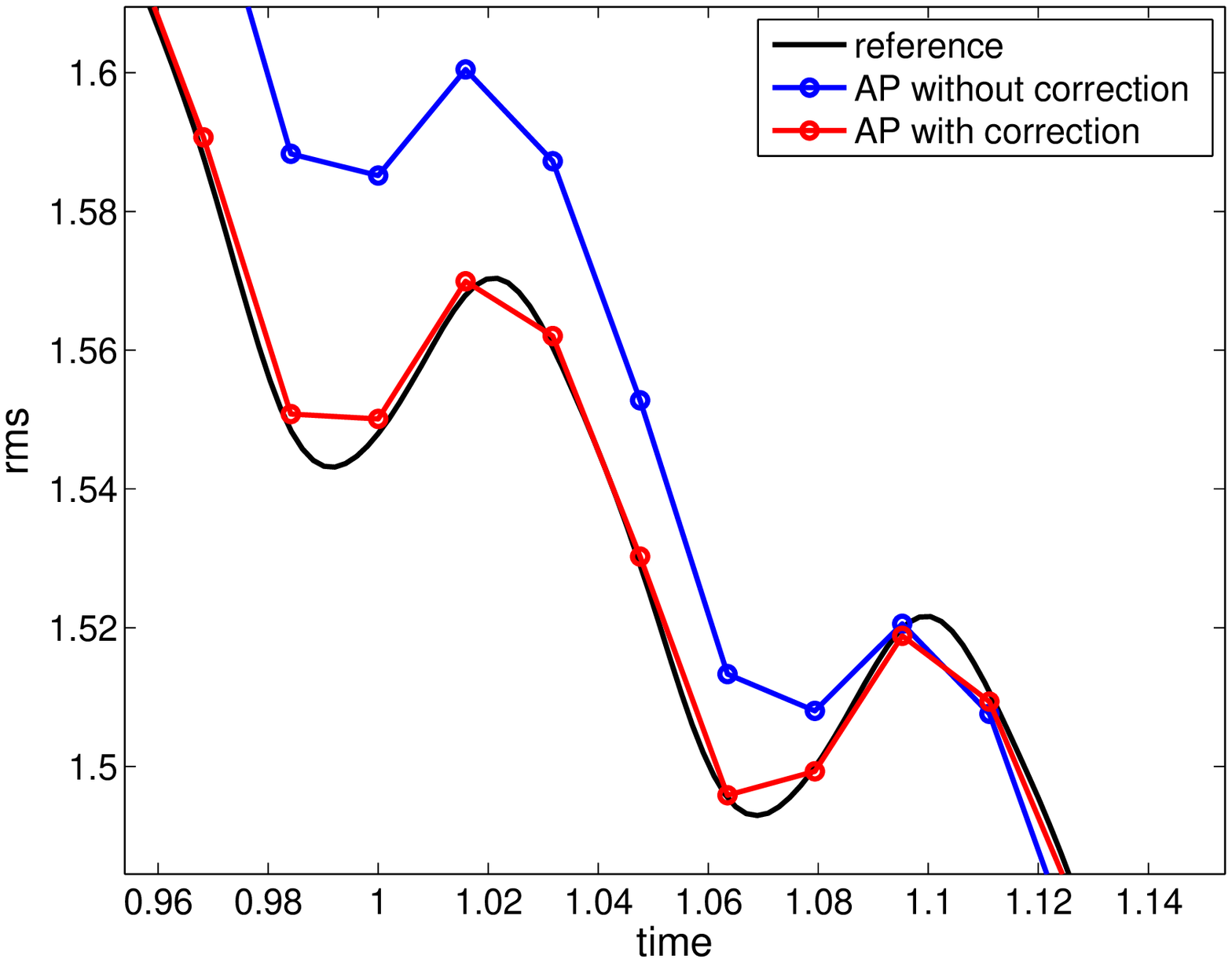}
\end{tabular}
\caption{Time history of $RMS(t)$ for $\varepsilon=0.025$, computed by {\sl AP with correction} and {\sl AP without correction}. On the right: zoom of the left figure.}
\label{figrms2}
\end{center}
\end{figure}

\begin{figure}[!htbp]
\begin{center}
\hspace*{-7mm}
\begin{tabular}{@{}c@{}c@{}}
\includegraphics[width=8cm]{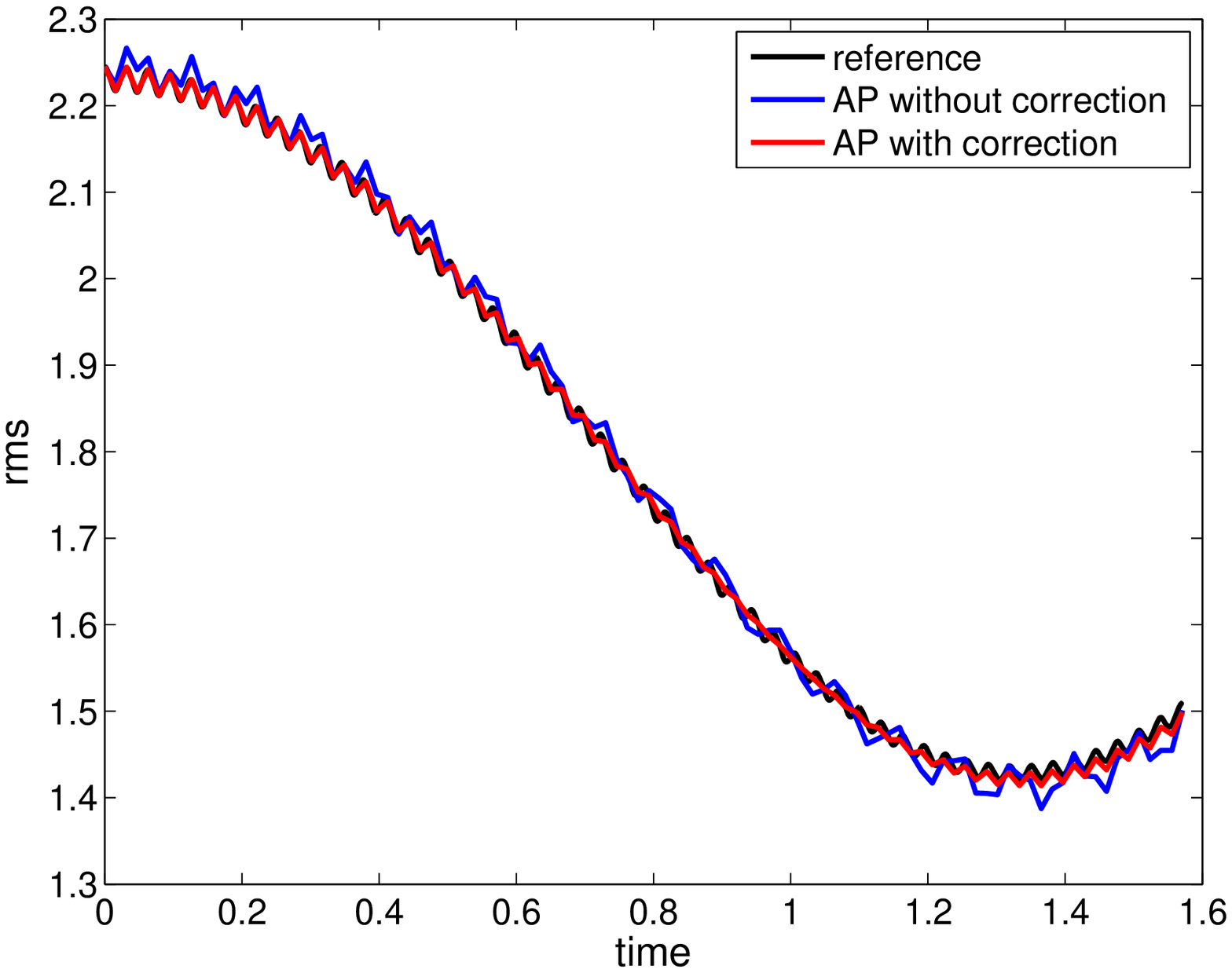} &     
\includegraphics[width=8cm]{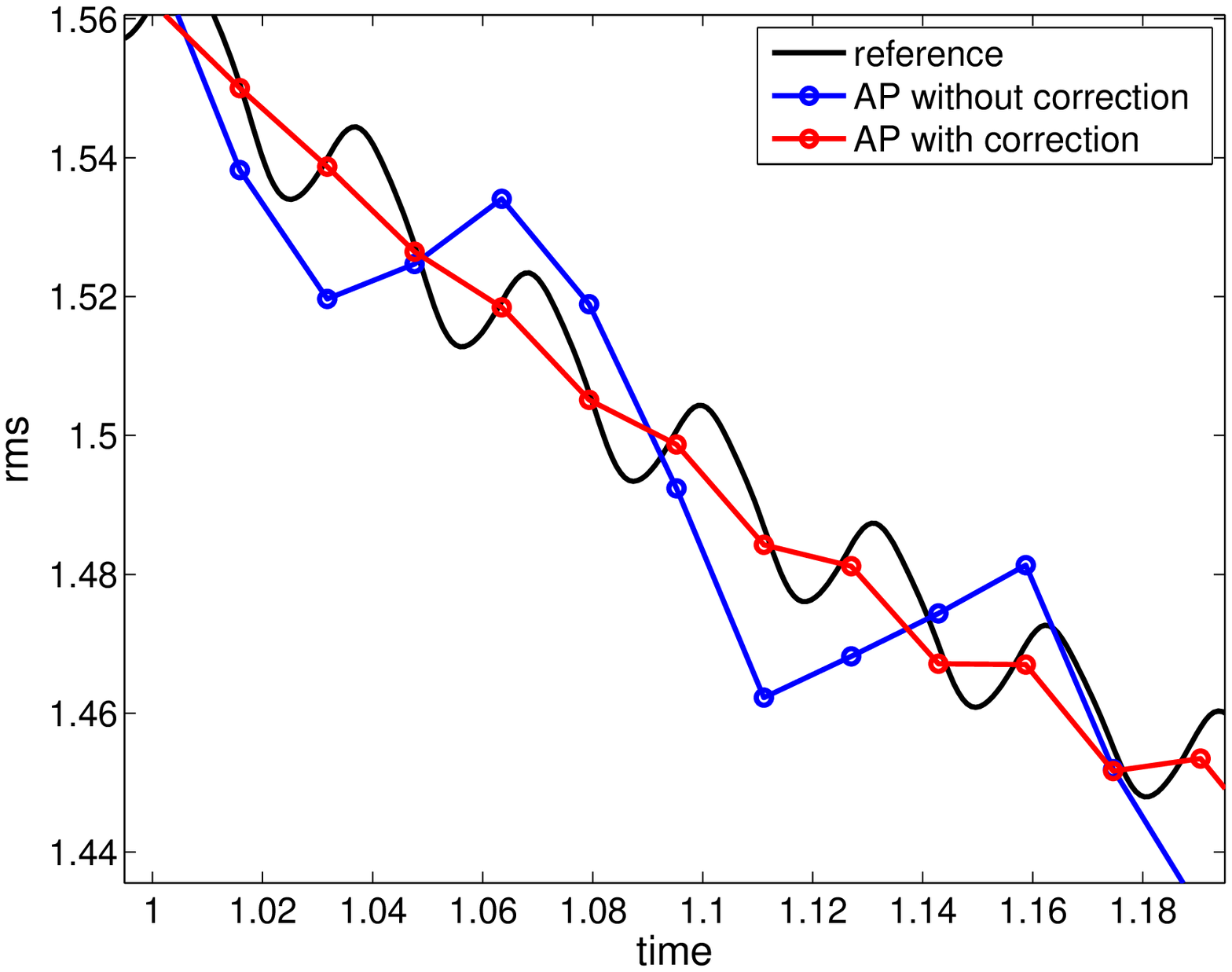}
\end{tabular}
\caption{Time history of $RMS(t)$ for $\varepsilon=0.01$, computed by {\sl AP with correction} and {\sl AP without correction}. On the right: zoom of the left figure.}
\label{figrms3}
\end{center}
\end{figure}

\begin{figure}[!htbp]
\begin{center}
\hspace*{-7mm}
\begin{tabular}{@{}c@{}c@{}}
\includegraphics[width=8cm]{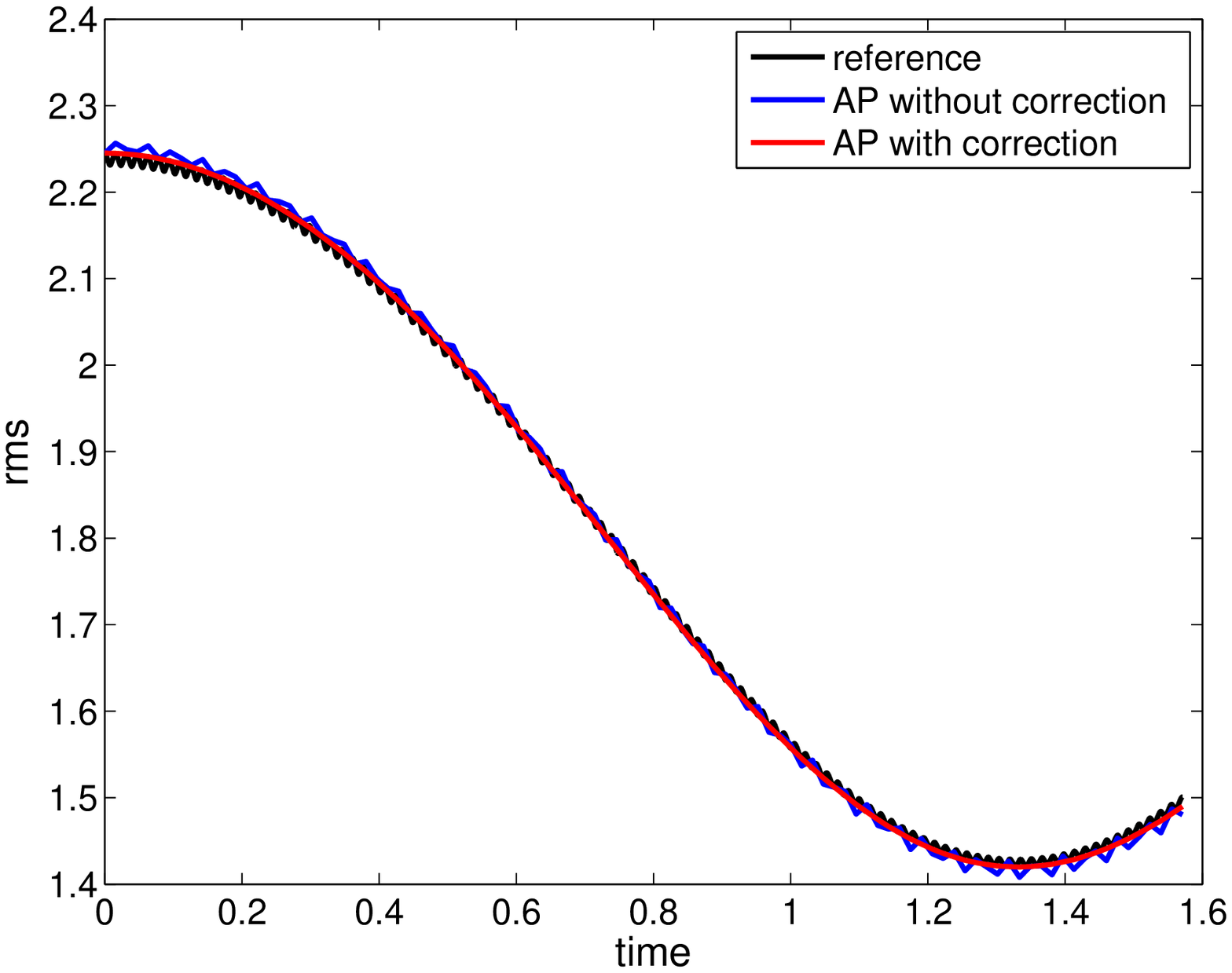} &     
\includegraphics[width=8cm]{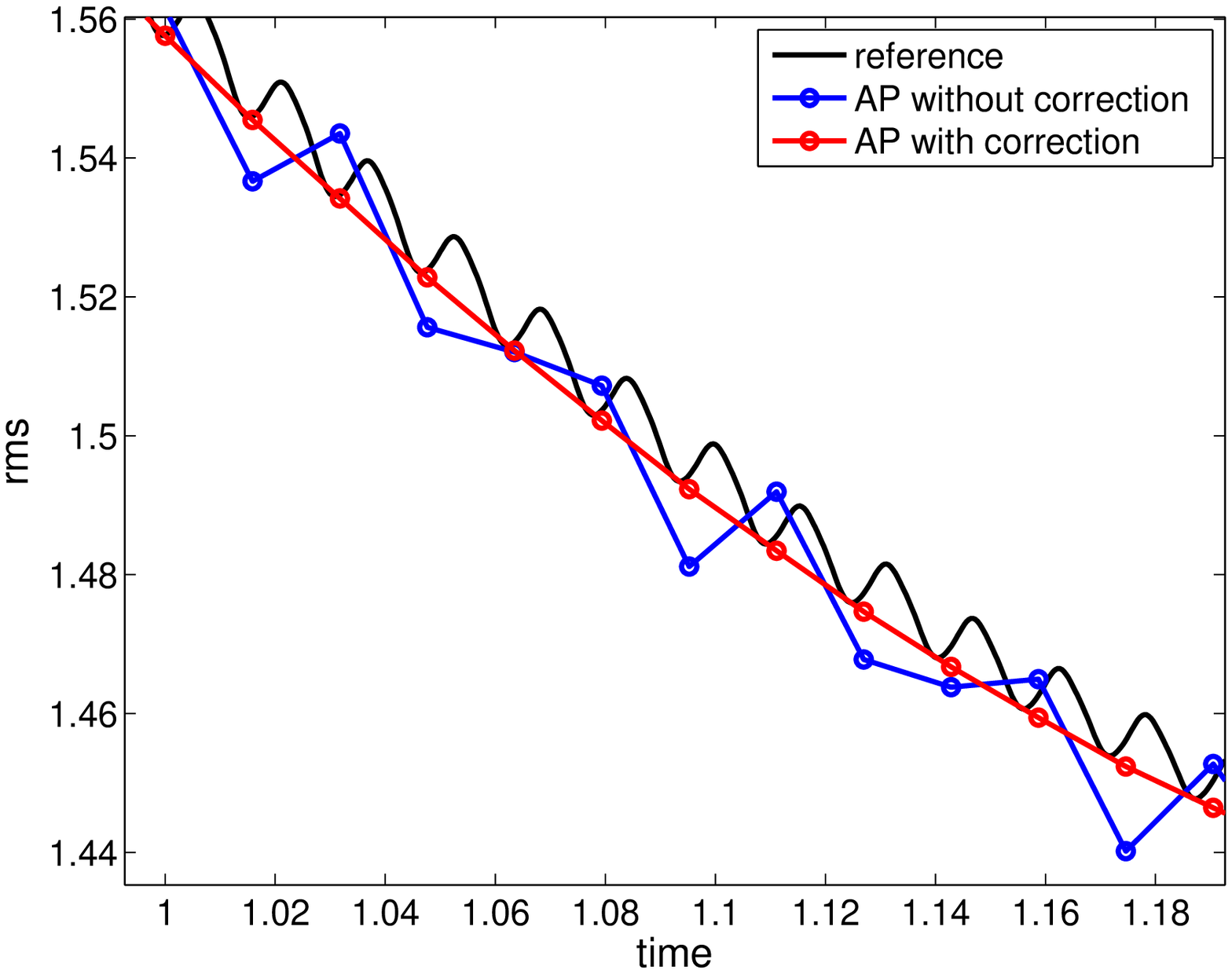}
\end{tabular}
\caption{Time history of $RMS(t)$ for $\varepsilon=0.005$, computed by {\sl AP with correction} and {\sl AP without correction}. On the right: zoom of the left figure.}
\label{figrms4}
\end{center}
\end{figure}

\bs
\ni
{\bf Evolution of an RMS quantity and observation of the oscillations in time}: Figures \ref{figrms1},  \ref{figrms2},  \ref{figrms3} and \ref{figrms4}
\nopagebreak

\ms
\ni
Let us now observe the evolution in time of a Root Mean Square (RMS) quantity associated to the filtered distribution function $\tilde f^\eps(t,\xi_1,\xi_2)$:
\be
\label{rms}
RMS(t)=\sqrt{\int_{\R^2}\xi_1^2\tilde f^\eps(t,\xi_1,\xi_2)d\xi }.
\ee
Note that, due to the filtering, this quantity does not oscillate at the limit $\eps=0$, and only the corrective terms for $\eps>0$ are rapidly oscillating. On Figures \ref{figrms1},  \ref{figrms2},  \ref{figrms3} and \ref{figrms4}, we represent respectively, for $\eps=0.05$, $\eps=0.025$, $\eps=0.01$ and $\eps=0.005$, the time history of $RMS(t)$ computed by the AP scheme with and without correction, and compare these numerical solutions to a reference solution. In all these simulations, we take $N=128$ and $\Delta t=0.02$. In particular, for $\eps=0.01$ and $\eps=0.005$, the time oscillation is not resolved by this time step. However, in all the cases, one can observe that the solution obtained by the scheme {\sl AP with correction} fits surprisingly well with the reference solution, see in particular the zooms on the right part of each figure: the red circles, which represent the only calculated points, are on the black reference curves, even when the oscillation is not resolved. This is another proof of the Asymptotic Preserving property of our scheme. On the contrary, one observes that the solution obtained with the scheme {\sl AP without correction} is less accurate: it converges to the right limit as $\eps\to 0$ but it is not able to correctly give the details of order $\mathcal O(\eps)$.

\subsection{The linear case}
\label{linear}

For all the numerical tests presented in this subsection, the self-consistent electric field is neglected in the model: we now set $E_f=0$ in \fref{eqf}. We are thus in the situation of the linear model of Section \ref{sectionlinear}, for which we have analytic expressions for the solution $F_{\rm limit}$ of the {\sl limit model} and for the solution $F_{\rm second order}$ of the {\sl second order model}, respectively given by \fref{lim1} and  \fref{lim2}.

\begin{figure}[!htbp]
\begin{center}
\begin{tabular}{@{}c@{}c@{}}
\includegraphics[width=7cm]{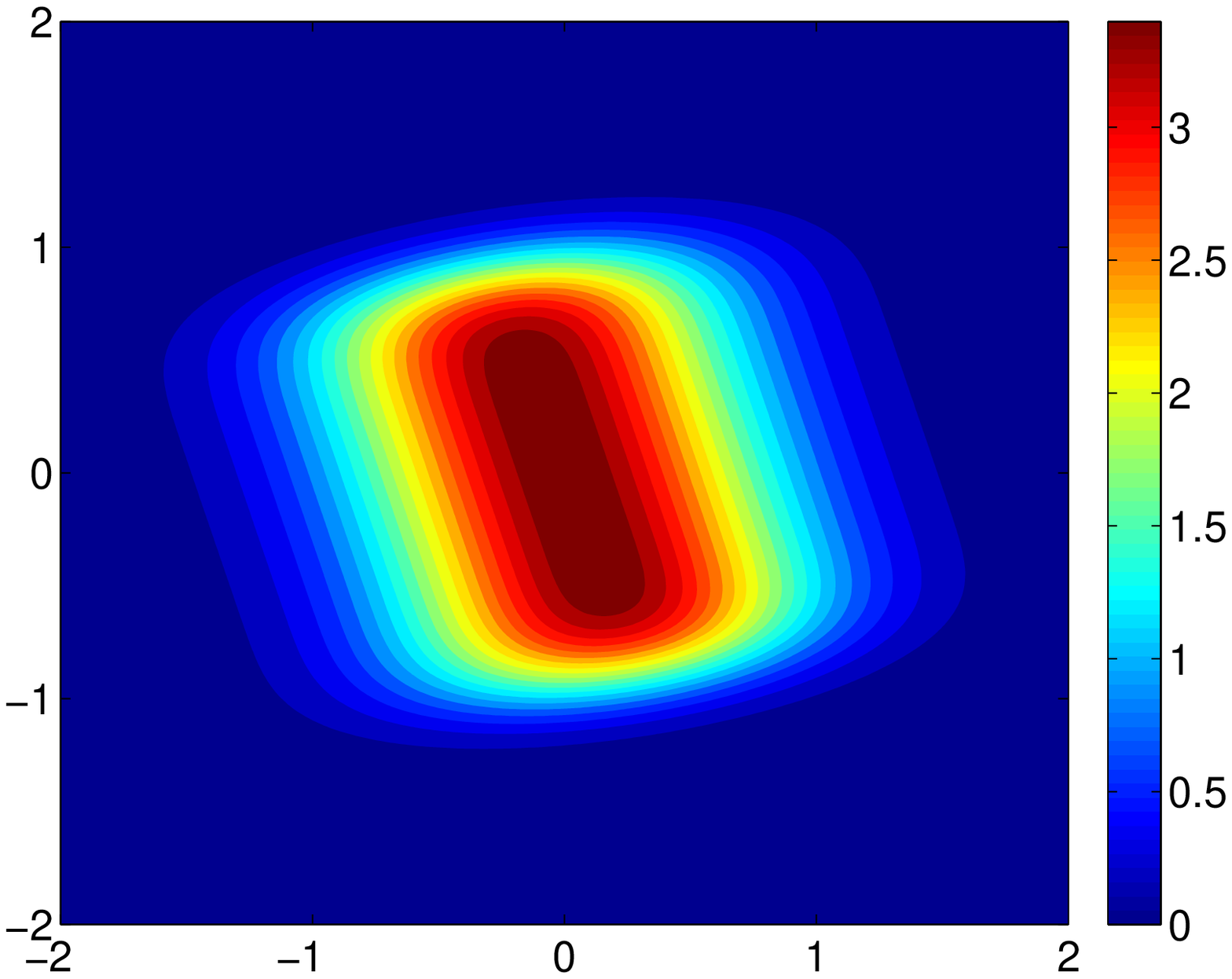} &     
\includegraphics[width=7cm]{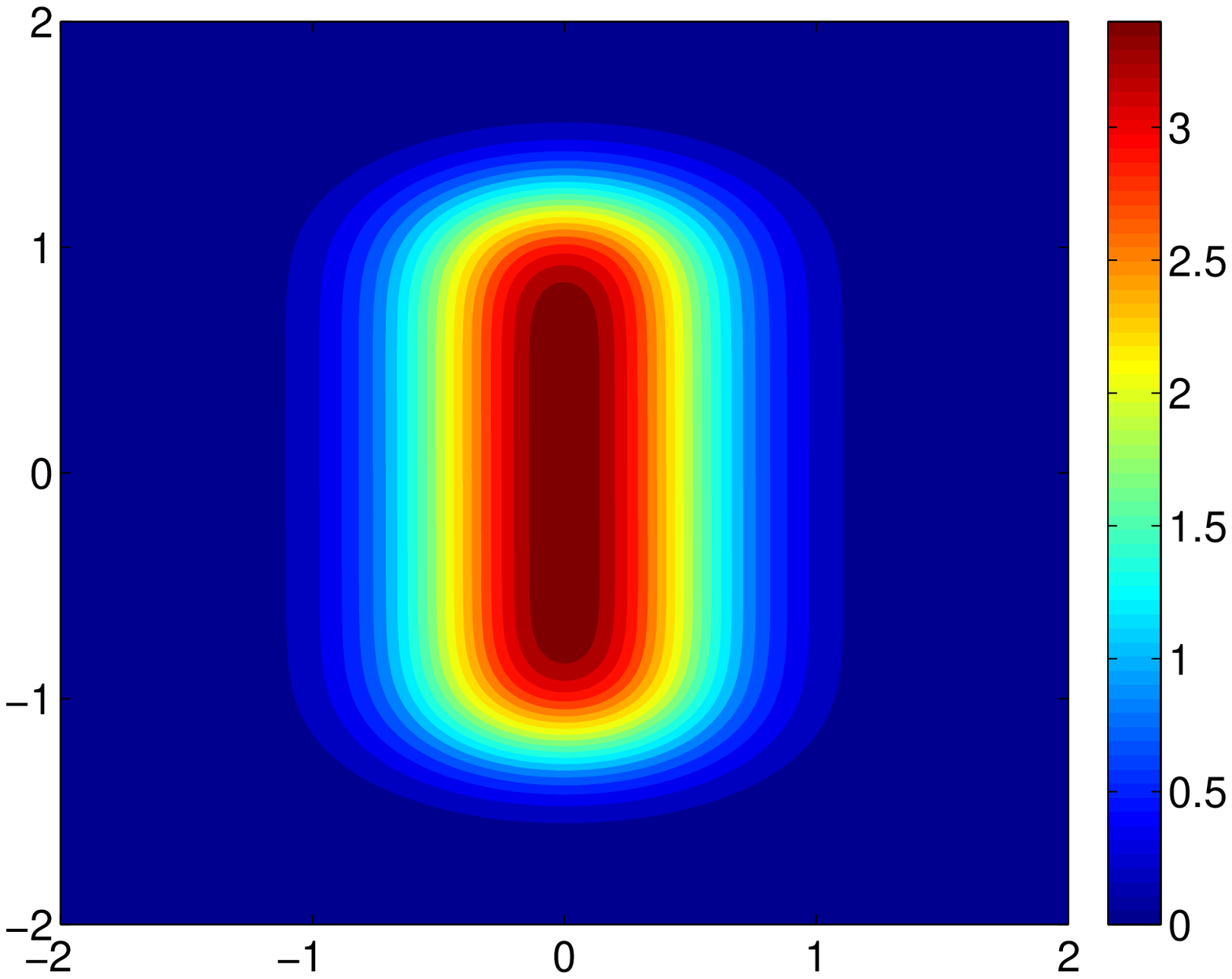}\\[-3mm]
\footnotesize $\eps=1$, {\sl splitting scheme}&\footnotesize $ \eps=0.01$, {\sl splitting scheme}\\
\includegraphics[width=7cm]{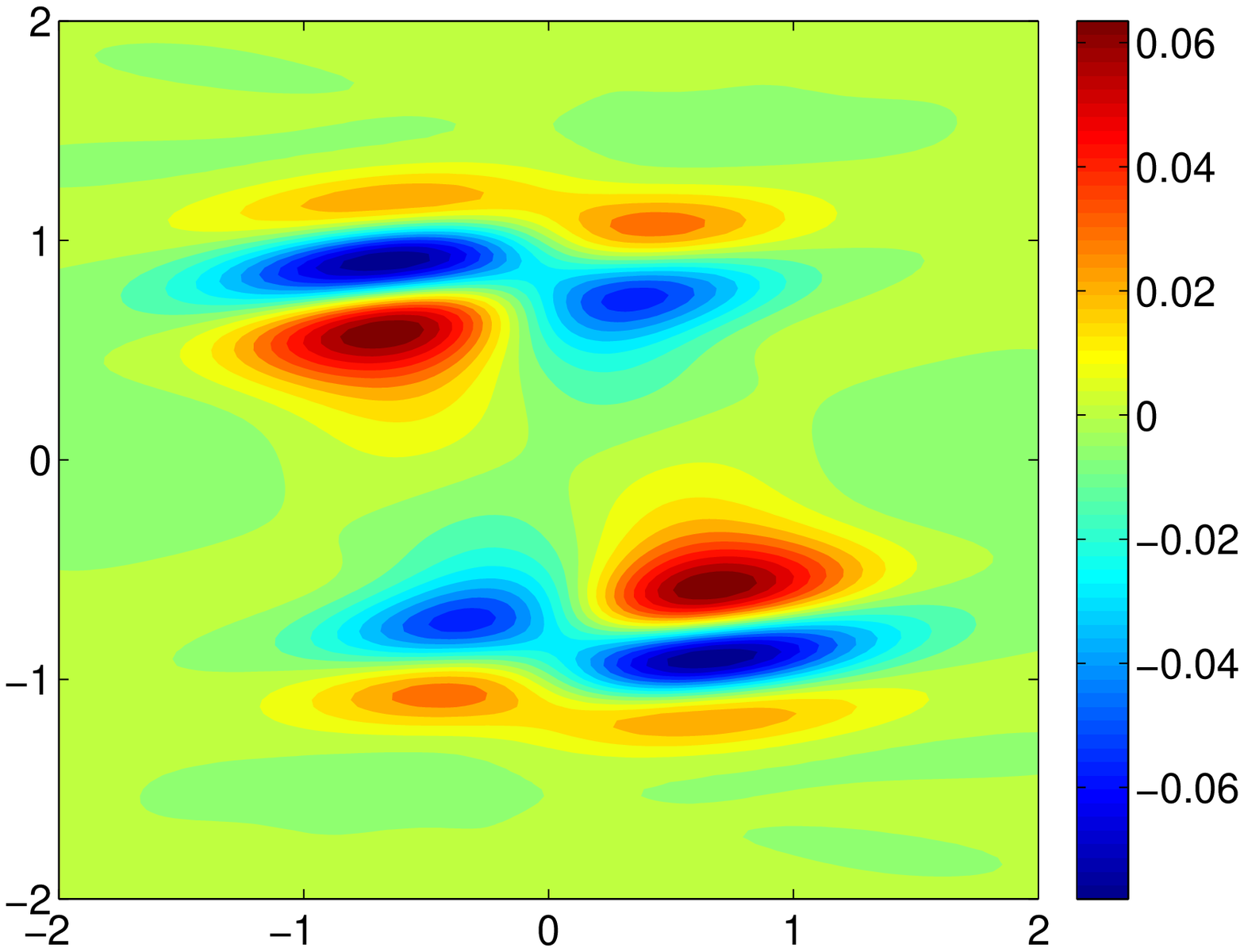} &     
\includegraphics[width=7cm]{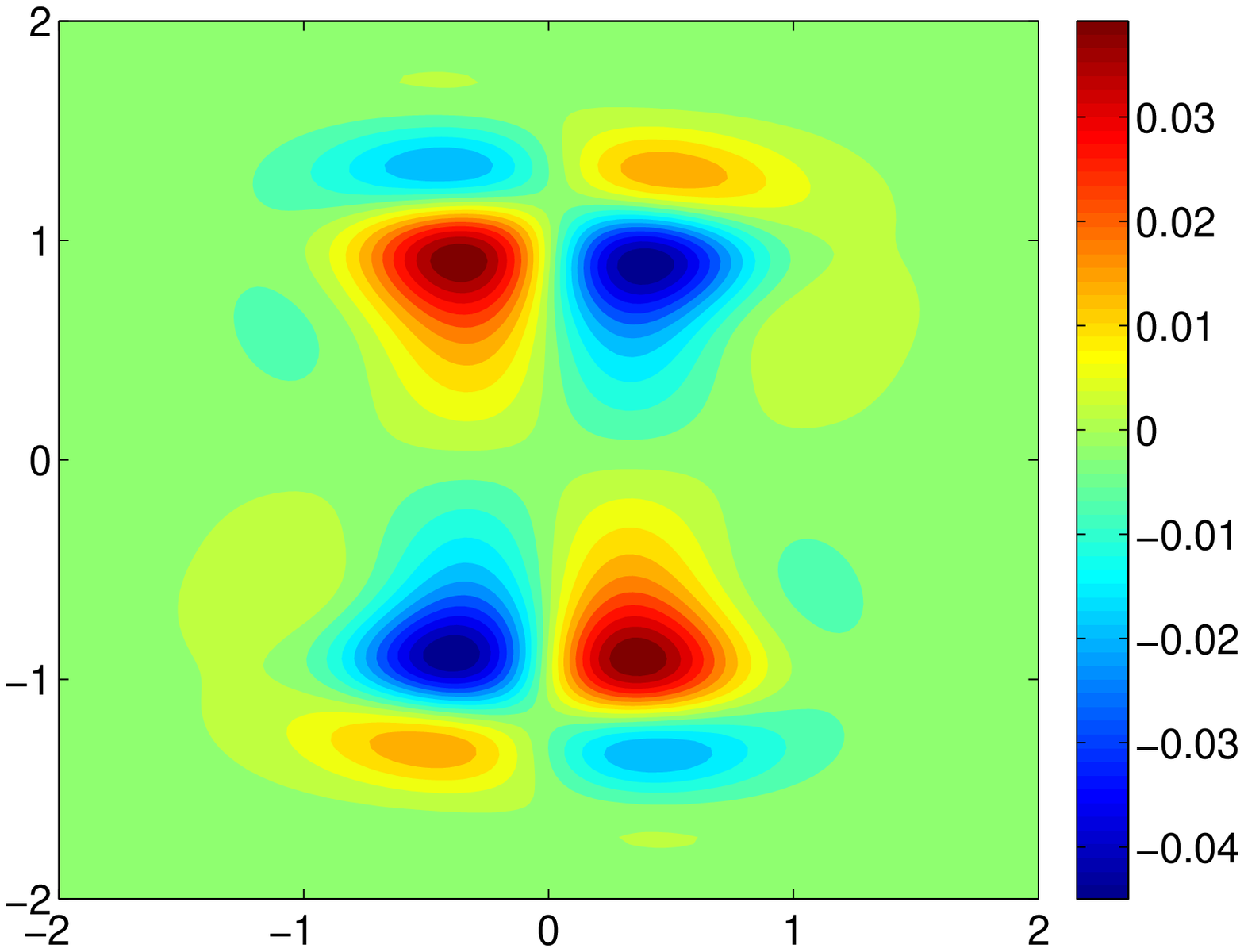}\\[-3mm]
\footnotesize $\tilde f^\eps_{\rm AP}-\tilde f^\eps_{\rm ref}$ for $\eps=1$&\footnotesize $\tilde f^\eps_{\rm AP}-\tilde f^\eps_{\rm ref}$ for $ \eps=0.01$\\
\includegraphics[width=7cm]{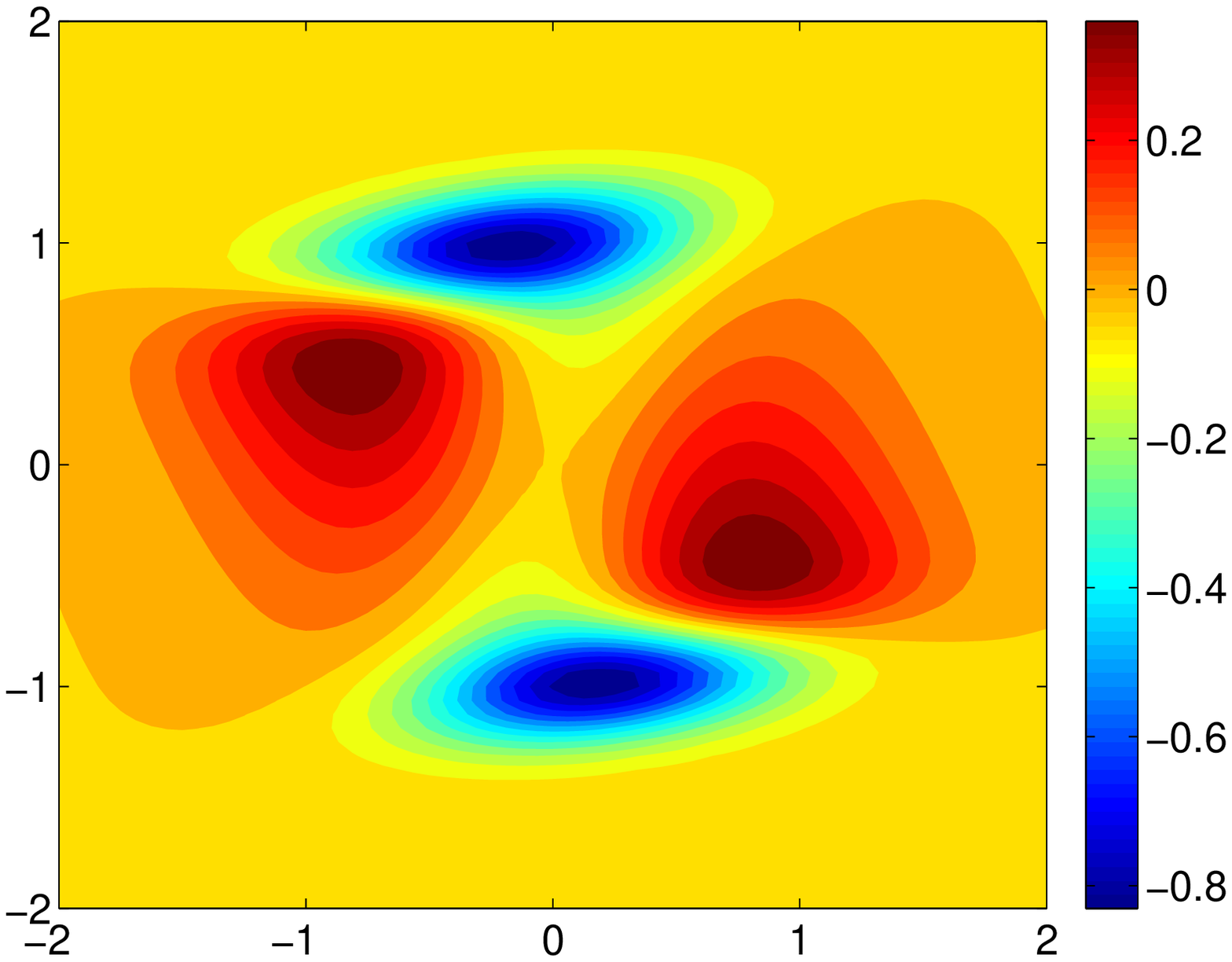} &     
\includegraphics[width=7cm]{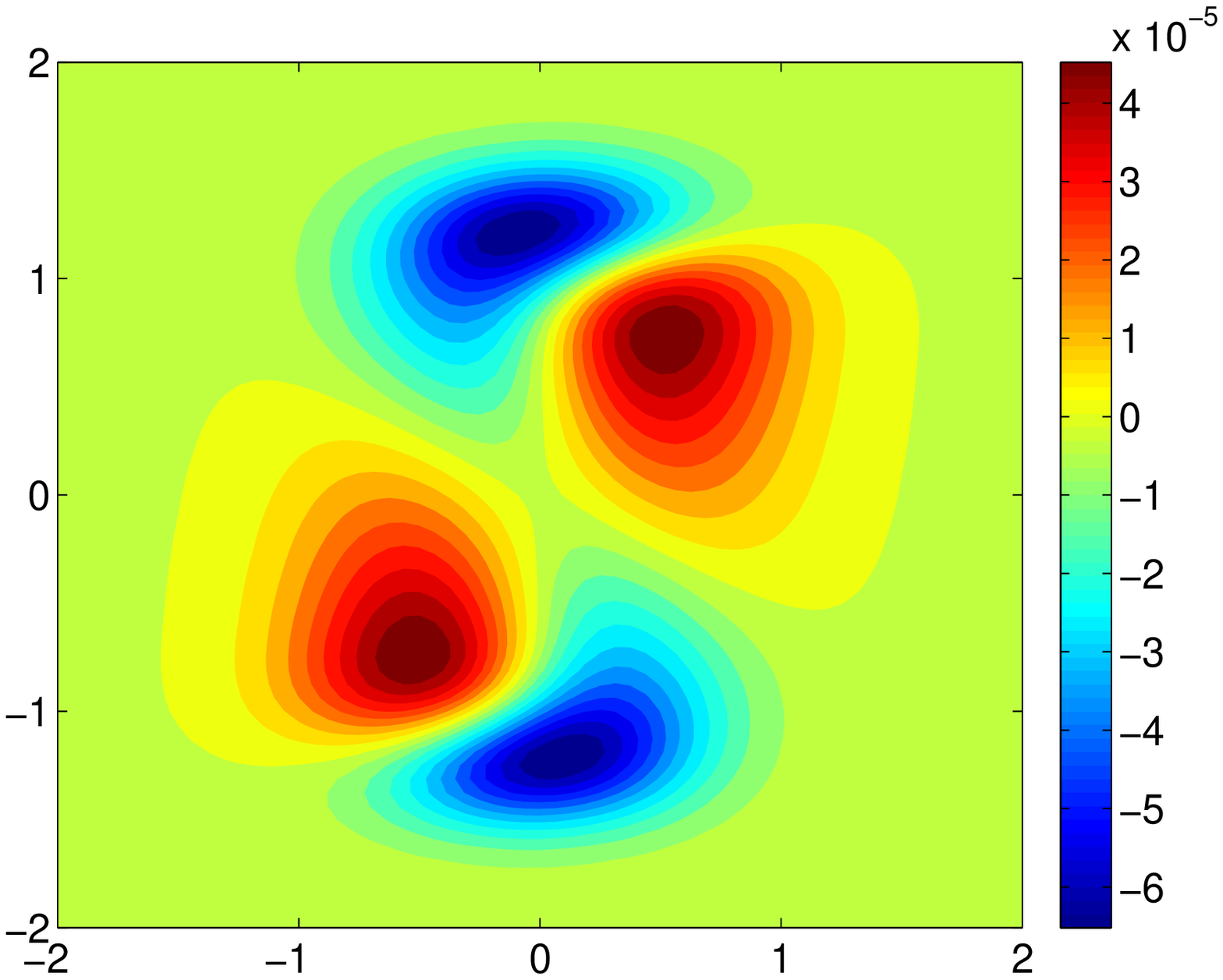}\\[-3mm]
\footnotesize $\tilde f^\eps_{\rm second}-\tilde f^\eps_{\rm ref}$ for $\eps=1$&\footnotesize $\tilde f^\eps_{\rm second}-\tilde f^\eps_{\rm ref}$ for $ \eps=0.01$
\end{tabular}
\caption{2D plots of $f^\eps(t,r,v)$ for $(r,v)\in [-2,2]^2$ for the linear beam model at time $t=2\pi$. Top line: reference solutions computed with the {\sl splitting scheme}, for $\eps=1$ and $\eps=0.01$. Middle line:  difference between the reference solutions and the numerical solutions with the scheme {\sl AP with correction}. Bottom line: difference between the reference solutions and the numerical solutions with the scheme {\sl second order model}.}
\label{figcontoursansPoi}
\end{center}
\end{figure}

\bs
\ni
{\bf Qualitative results for two regimes of $\eps$}: Figure \ref{figcontoursansPoi} 
\nopagebreak

\ms
\ni
As above for the nonlinear model, let us start with a few qualitative results. We first represent the 2D plot of the solution of the linear problem, at the final time $t_{final}=2\pi$, for the values $\eps=1$ and $\eps=0.01$. On the top line of Figure \ref{figcontoursansPoi}, we represent the plot of the reference solution $\tilde f^\eps_{\rm ref}$ computed with the {\sl splitting scheme}. Note that, for $\eps=0.01$ (top-right plot of the figure), the solution cannot be distinguished from the solution of the limit model, which is simply the initial data rotated of an angle $\pi/2$ (compare with Figure \ref{figf0}). Indeed, if $\omega$ is given by \fref{omega}, one has $\omega t_{final}=\frac{\pi}{2}+\frac{5\pi}{9600}\approx \frac{\pi}{2}$.

On the middle line of the same figure, we represent the 2D plot of the difference $\tilde f^\eps_{\rm AP}-\tilde f^\eps_{\rm ref}$, where $\tilde f^\eps_{\rm AP}$ is the numerical solution with the scheme {\sl AP with correction} ($N=256$), for $\eps=1$ and $\eps=0.01$ and, on the bottom line of the figure, we represent the difference $\tilde f^\eps_{\rm second}-\tilde f^\eps_{\rm ref}$, where $\tilde f^\eps_{\rm second}$ is the analytic solution of the {\sl second order model}, for the same values of $\eps$. We observe the following facts. The error for the AP scheme is almost the same (around $10^{-2}$) for the two values of $\eps$, whereas for the {\sl second order model}, the results are very dependent of $\eps$: for $\eps=1$, the error in $L^\infty$ norm is close to 1, whereas for $\eps=0.01$, this error is around $10^{-4}$. The  {\sl second order model} can be used only for small values of $\eps$ (its incapacity to predict the solution for $\eps=1$ is even clearer below on the RMS test).

\begin{figure}[!htbp]
\begin{center}
\includegraphics[width=10cm]{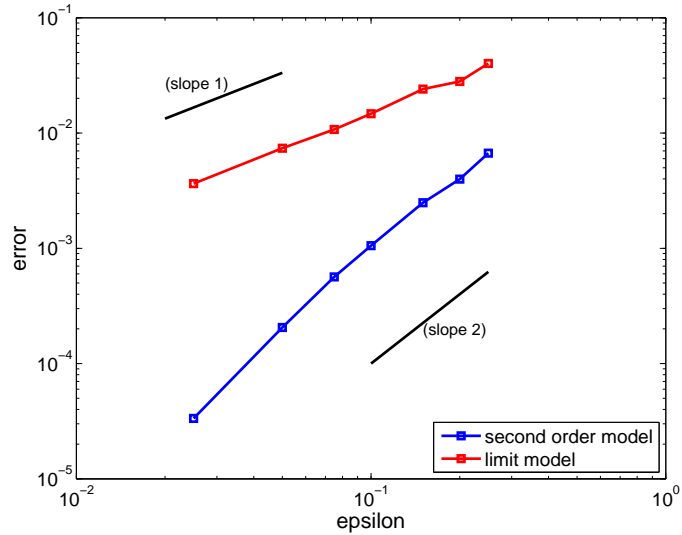}
\caption{Plot of the errors between the {\sl limit model} and the reference solution and between the {\sl second order model} and the reference solution, as functions of $\eps$.}
\label{figordrelim}
\end{center}
\end{figure}

\bs
\ni
{\bf Accuracy of the limit model and of the second order model}: Figure \ref{figordrelim}
\nopagebreak

\ms
\ni
Let us confirm more quantitatively the above observations. In the next table, we give the relative $L^\infty$ errors between the approximate solutions and the reference solution (still at time $t_{final}=2\pi$ and, for the AP scheme, we take $N=256$). This error is defined by $$\mbox{error}=\frac{\|\tilde f^\eps_{\rm approx}-\tilde f^\eps_{\rm ref}\|_{L^\infty}}{\|\tilde f^\eps_{\rm ref}\|_{L^\infty}}.$$
\begin{center}
\begin{tabular}{|l|l|l|l|l|l|l|l|l|c|c|r|} 
\hline
$\varepsilon$ & 1 & 0.5& 0.25  & 0.1 & 0.01 \\
\hline
error for \sl AP with correction     &  1.8\,\% & 1.5\,\%   & 1.5\,\%     & 1.4\,\%     & 1.3\,\%  \\
\hline
error for the \sl second order model  &  18\,\%& 4\,\%  & 1\,\%    & 0.15\,\%     & 0.001\,\%   \\
\hline
error for the \sl limit model &  37\,\% & 18\,\%   & 8.6\,\%     & 3.3\,\%     & 0.3\,\%   \\
\hline
\end{tabular} 
\end{center}

\ms
\ni
This table indicates that the error produced by the scheme {\sl AP with correction} is independent of $\eps$ (as we shown in the previous subsection for the nonlinear case), and that the {\sl limit model} and the {\sl second order model} seem respectively of orders 1 and 2 in $\eps$. On Figure \ref{figordrelim}, we illustrate numerically the accuracies of these two asymptotic models with respect to $\eps$ by plotting in logarithmic scales the $L^1([0,t_{final}])$ norm of the difference $RMS_{\rm approx}(t)-RMS_{\rm reference}(t)$, for these two models. One can check on this figure that the errors produced by these models are respectively $\mathcal O(\eps)$ and $\mathcal O(\eps^2)$.  In other terms, we confirm numerically the results given by Proposition \ref{theo2}. 

\begin{figure}[!htbp]
\begin{center}
\hspace*{-7mm}
\begin{tabular}{@{}c@{}c@{}}
\includegraphics[width=8cm]{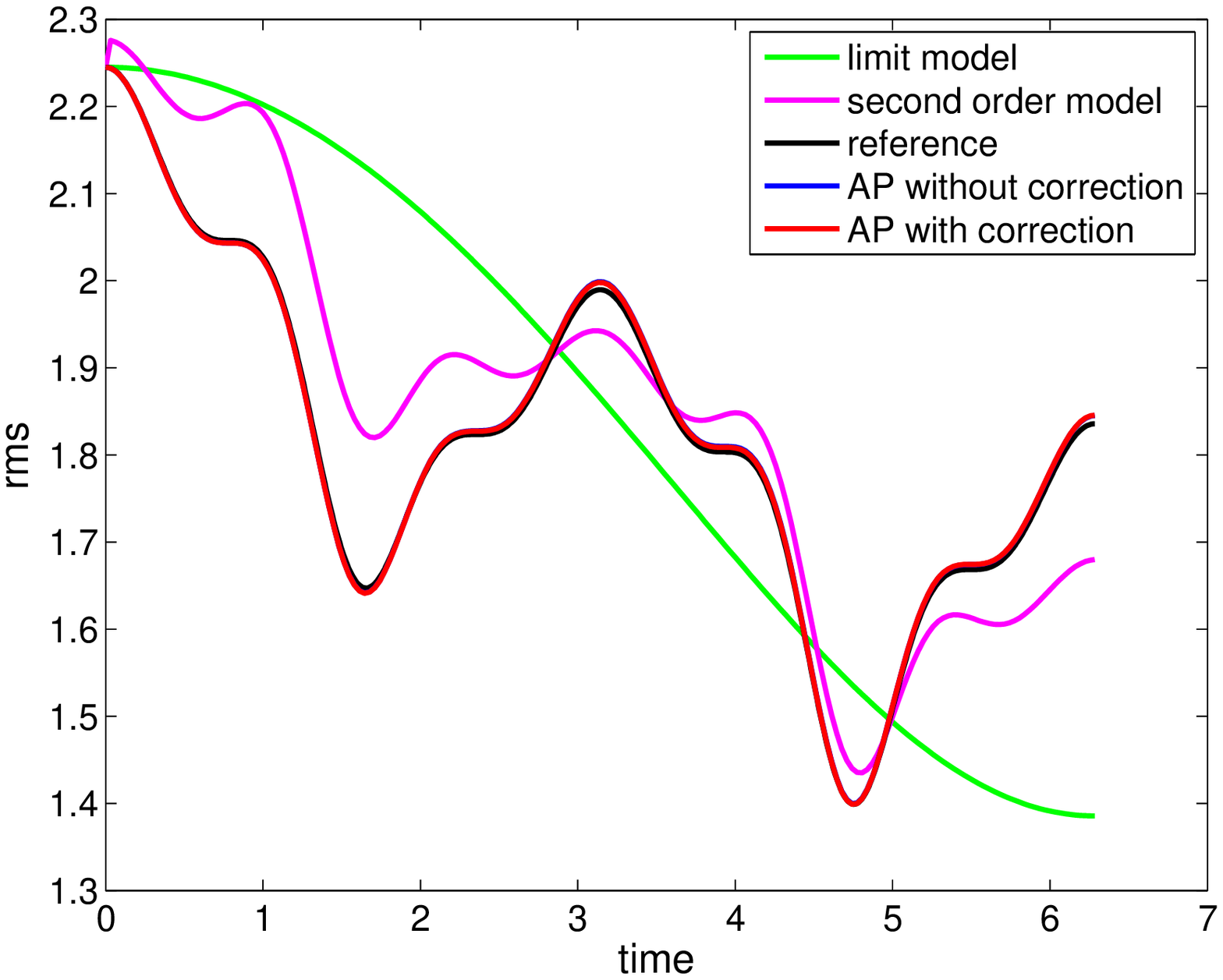} &     
\includegraphics[width=8cm]{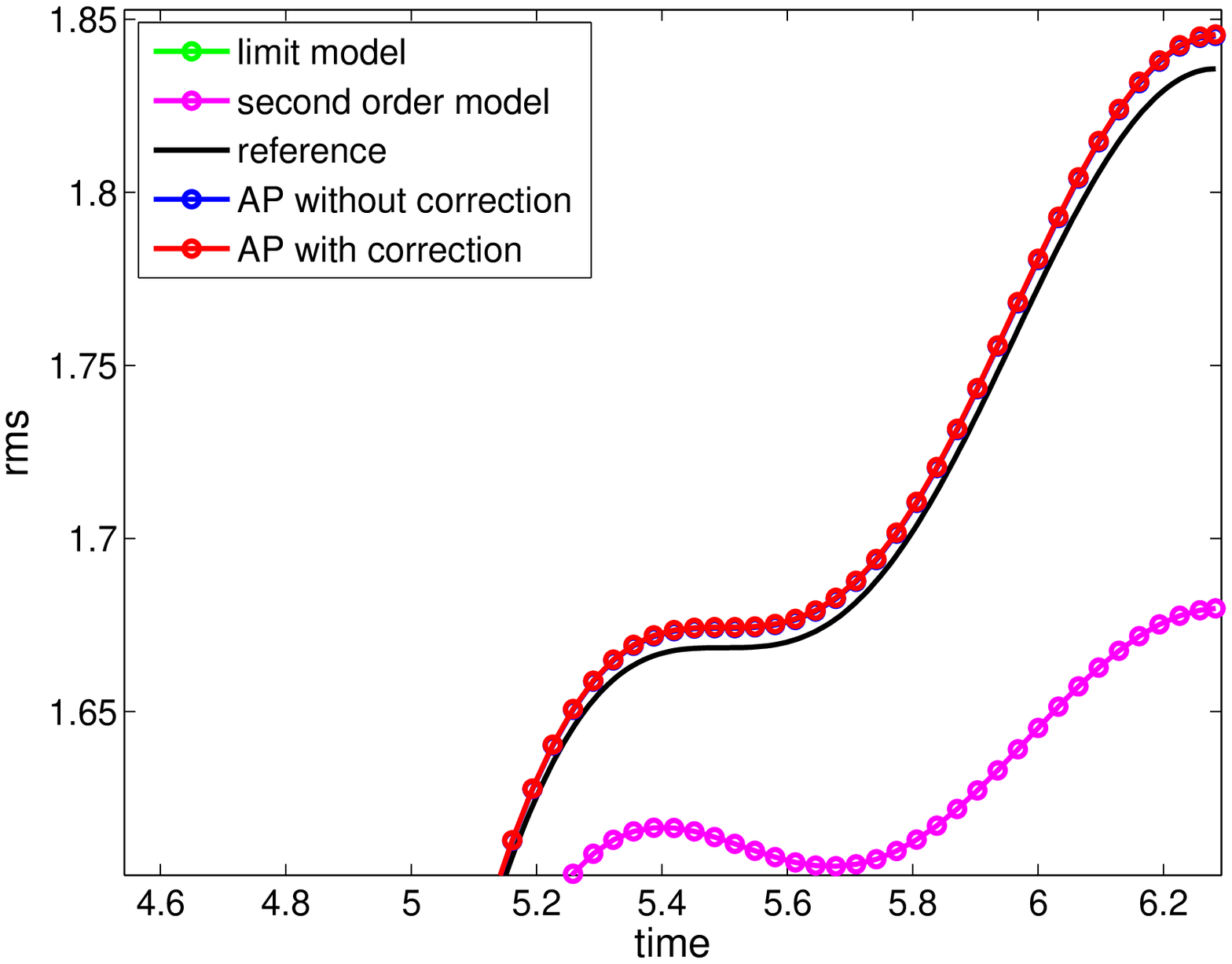}
\end{tabular}
\caption{Time history of $RMS(t)$ for $\varepsilon=1$, in the linear situation, computed with {\sl AP} (with or without correction; these two curves red and blue coincide), with the {\sl limit model} and with the {\sl second order model}. On the right: zoom of the left figure.}
\label{figrmssans1}
\end{center}
\end{figure}

\begin{figure}[!htbp]
\begin{center}
\hspace*{-7mm}
\begin{tabular}{@{}c@{}c@{}}
\includegraphics[width=8cm]{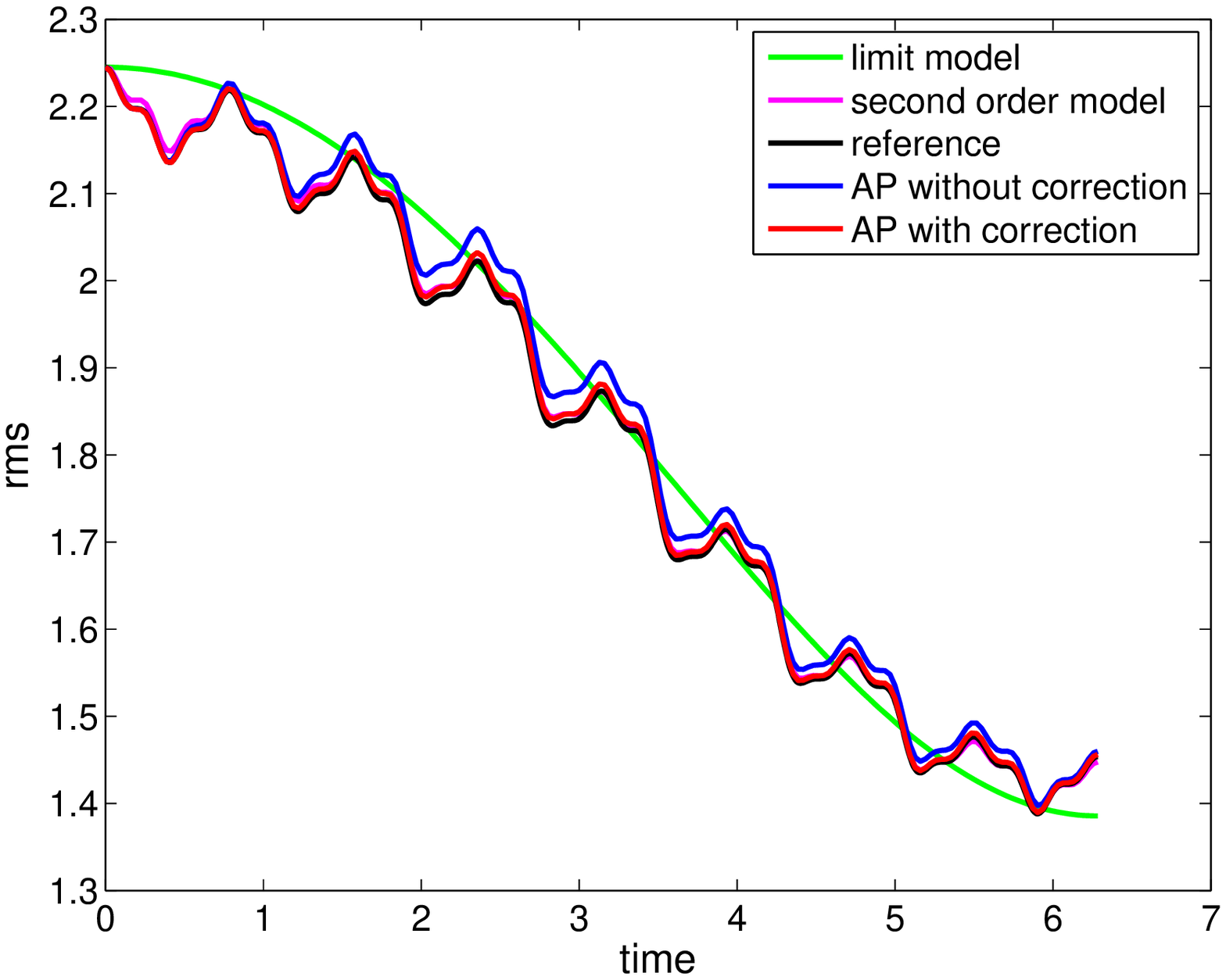} &     
\includegraphics[width=8cm]{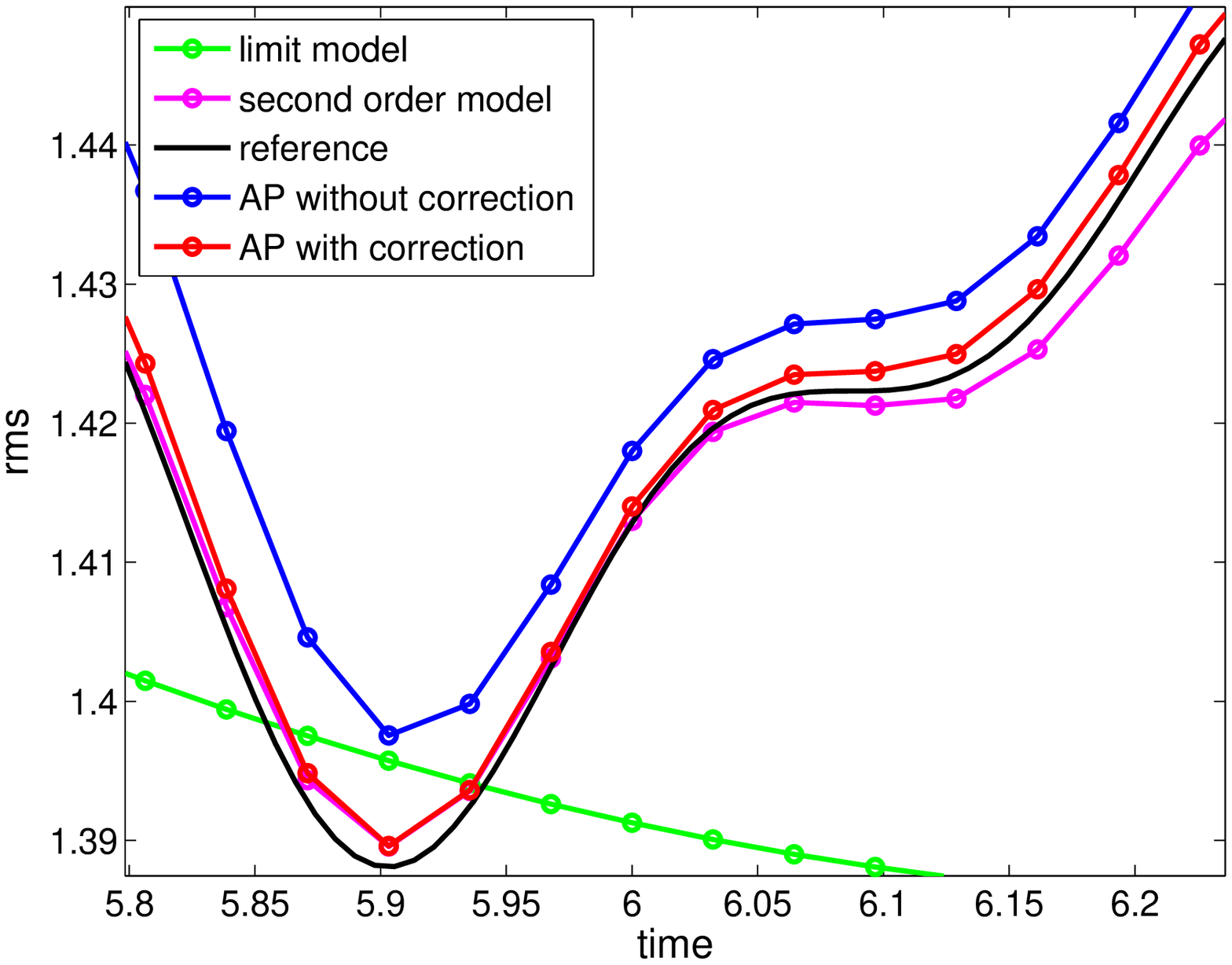}
\end{tabular}
\caption{Time history of $RMS(t)$ for $\varepsilon=0.25$, in the linear situation, computed with {\sl AP} (with or without correction), with the {\sl limit model} and with the {\sl second order model}. On the right: zoom of the left figure.}
\label{figrmssans2}
\end{center}
\end{figure}

\begin{figure}[!htbp]
\begin{center}
\hspace*{-7mm}
\begin{tabular}{@{}c@{}c@{}}
\includegraphics[width=8cm]{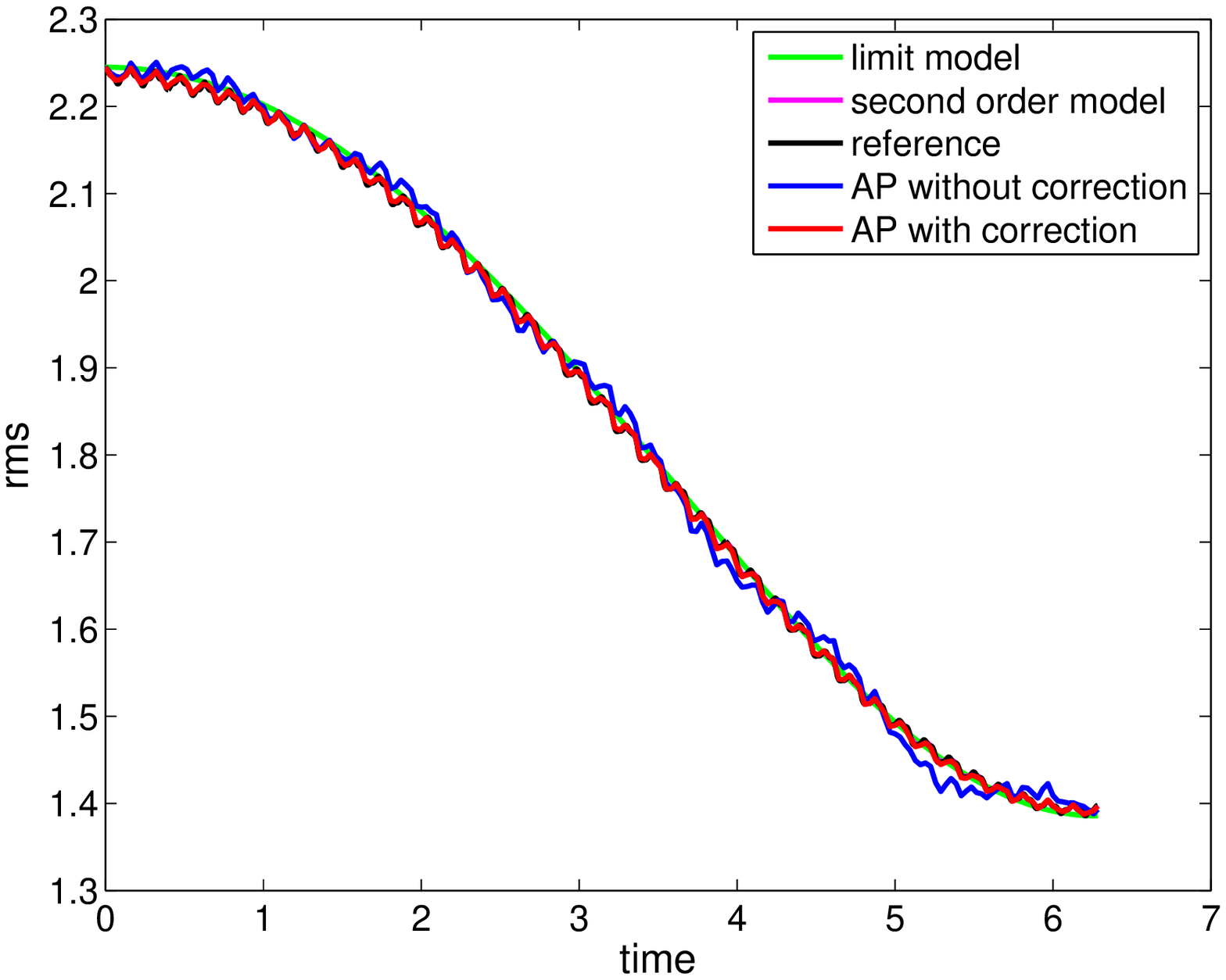} &     
\includegraphics[width=8cm]{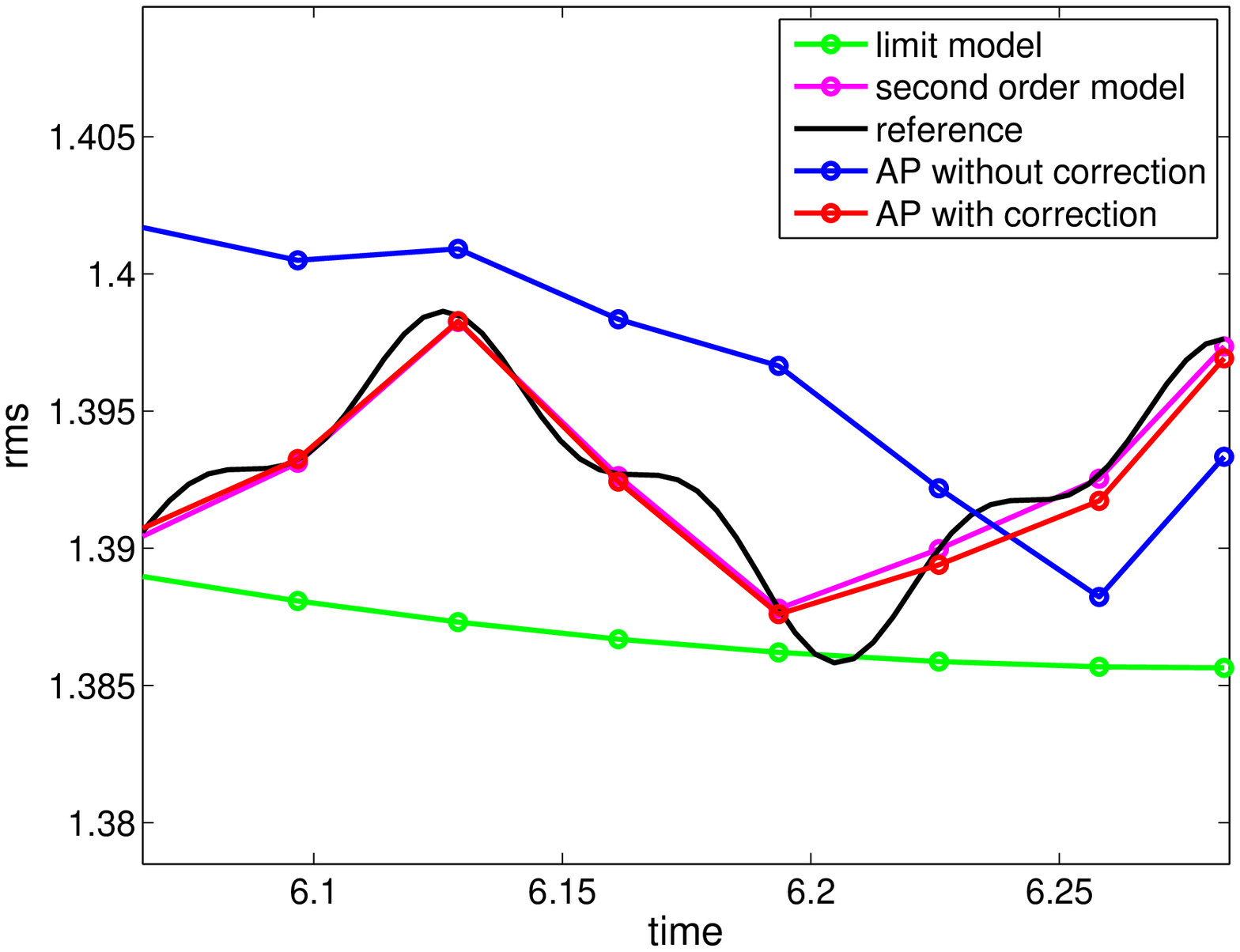}
\end{tabular}
\caption{Time history of $RMS(t)$ for $\varepsilon=0.05$, in the linear situation, computed with {\sl AP} (with or without correction), with the {\sl limit model} and with the {\sl second order model}. On the right: zoom of the left figure.}
\label{figrmssans3}
\end{center}
\end{figure}

\bs
\ni
{\bf Evolution of the RMS}: Figures \ref{figrmssans1},  \ref{figrmssans2} and \ref{figrmssans3}
\nopagebreak

\ms
\ni
We now observe the evolution in time of the RMS quantity defined by \fref{rms}. On Figures \ref{figrmssans1},  \ref{figrmssans2} and \ref{figrmssans3}, we represent respectively, for $\eps=1$, $\eps=0.25$ and $\eps=0.05$, the time history of $RMS(t)$ computed by the AP scheme in the linear case, with and without correction, and compare these numerical solutions to a reference solution and to the solutions of the {\sl limit model} and of the {\sl second order model}. In all these simulations, we take $N=64$ and $\Delta t=0.02$. In all the cases, one can observe that the solution obtained by {\sl AP with correction} fits very well with the reference solution (see the zooms on the right part of each figure) even when the oscillation is not well resolved. As for the nonlinear case, one observes that the solution obtained with {\sl AP without correction} is less accurate when $\eps$ is small and is not able to reproduce the details of order $\mathcal O(\eps)$.

One also observes that the {\sl limit model} is only able to give the averaged behavior of the curve. The {\sl second order model} is much better and follows the oscillations for small values of $\eps$. On Figure \ref{figrmssans3}, for $\eps=0.05$, its solution coincides with the reference solution and is more precise than {\sl AP with correction}. Recall indeed that the error made by the scheme {\sl AP with correction} is proportional to $\Delta t^2+\Delta \xi^2\approx 0.02$ whereas the error made by the {\sl second order model} is proportional to $\eps^2\approx 0.002$. Indeed, in this linear context, the  {\sl second order model} is analytic and does not produce any error in time or space. Obviously, 
in a more general case, the {\sl second order model} will also generate an error due to 
its space-time discretization. On Figure \ref{figrmssans2}, for $\eps=0.25$ ($\eps^2\approx 0.06$), the errors made by the two methods {\sl AP with correction} and {\sl second order model} are comparable. Finally, on Figure  \ref{figrmssans1}, for $\eps=1$, it appears again that the error made by the {\sl second order model} is of order $\mathcal O(1)$: this confirms that this averaged model is useless when $\eps$ is not small.

\section{Conclusion}

In this work, we have presented a general strategy to construct the so-called Asymptotic Preserving (AP) numerical schemes
for a family of highly oscillatory problems. Although we focus on the particular case of a charged particle beam to illustrate our strategy, the approach may be applied to many  known physical models belonging to this family of highly oscillatory problems.  Averaged models are usually used to approximate this type of problems but these models are not relevant in the intermediate regime since they miss important
informations from the original problem.

The starting idea in this construction is to write the oscillatory problem into a "double-scale" formulation where the rapid and slow time scales are separated, making the new distribution function more regular in some sense.  The new structure then  suggests to follow  a similar strategy as in the collisional case to develop AP schemes  on this formulation. 
However the completely different nature of  highly oscillatory problems (compared to collisional kinetic equations) induces new important difficulties. First, the double-scale formulation is overdetermined in the sense that a large family of initial data for this formulation is allowed.  We show in this paper that there is a suitable choice to make on this initial data in order to maintain the regularity of the distribution function at different orders of the oscillation parameter.  More precisely, the initial data is chosen to fit with a Chapman-Enskog like expansion which ensures a separation of the rapid and slow time scales at different orders of the expansion.
Based on this formulation, we then derive an Asymptotic Preserving scheme for the  original problem and show that  time-space discretizations of order $2$ are necessary to numerically observe the fine structures and filamentations that are generated by the coupling of Vlasov and Poisson equations.  Several numerical tests are performed to show the  efficiency of our strategy: uniform accuracy and ability to capture the oscillations of different magnitudes and the long time behavior. 

We emphasize that the AP property of our scheme is shown by making links with the so-called micro-macro decomposition, which is known to be a flexible tool to develop AP schemes in the context of collisional kinetic equations. In particular, this decomposition may be used to extend the present approach to other highly oscillatory problems such as the charged particle beam with a diffusion scaling, the guiding-center asymptotics and the finite Larmor radius approximation. This will be the subject of future works.

\end{document}